\documentclass{amsart}
\usepackage{amsmath}
\usepackage{amssymb}
\usepackage{extarrows}
\usepackage{mathabx}
\usepackage[all]{xy}%commutative graph
\usepackage{amsbsy}
\usepackage{graphicx}
\usepackage{appendix}
\usepackage{float}
\usepackage{subfigure}
\usepackage[numbers,sort&compress]{natbib}
\usepackage{amsthm}
\usepackage{geometry}
\usepackage{fancyhdr}%Heads
\usepackage{color} %colors
\usepackage{overpic}%pictures with text
\usepackage{enumerate}
\usepackage{mathrsfs}
\usepackage{colonequals}
\usepackage{microtype}
\usepackage{hyperref}
\usepackage{diagbox}

\newtheorem{thm}{Theorem}[section]
\newtheorem{cor}[thm]{Corollary}
\newtheorem{prop}[thm]{Proposition}
\newtheorem{lem}[thm]{Lemma}

\theoremstyle{definition}
\newtheorem{defn}[thm]{Definition}
\newtheorem{cons}[thm]{Construction}
\newtheorem{exmp}[thm]{Example}
\newtheorem{conj}[thm]{Conjecture}
\newtheorem*{fact}{Fact}

\newtheorem*{conv}{Convention}

\newtheorem*{ack}{Acknowledgement}

\theoremstyle{remark}
\newtheorem{rem}[thm]{Remark}

\numberwithin{equation}{section}

\newcommand{\beq}{\begin{equation*}\begin{aligned}}
\newcommand{\eeq}{\end{aligned}\end{equation*}}
\newcommand{\bpf}{\begin{proof}}
\newcommand{\epf}{\end{proof}}
\newcommand{\bthm}{\begin{thm}}
\newcommand{\ethm}{\end{thm}}
\newcommand{\bprop}{\begin{prop}}
\newcommand{\eprop}{\end{prop}}
\newcommand{\bcor}{\begin{cor}}
\newcommand{\ecor}{\end{cor}}
\newcommand{\blem}{\begin{lem}}
\newcommand{\elem}{\end{lem}}
\newcommand{\bdefn}{\begin{defn}}
\newcommand{\edefn}{\end{defn}}
\newcommand{\bcons}{\begin{cons}}
\newcommand{\econs}{\end{cons}}
\newcommand{\bexmp}{\begin{exmp}}
\newcommand{\eexmp}{\end{exmp}}
\newcommand{\brem}{\begin{rem}}
\newcommand{\erem}{\end{rem}}
\newcommand{\bfa}{\begin{fact}}
\newcommand{\efa}{\end{fact}}
\newcommand{\benu}{\begin{enumerate}[(1)]}
\newcommand{\eenu}{\end{enumerate}}

\newcommand{\bdia}{\begin{displaymath}\xymatrix}
\newcommand{\edia}{\end{displaymath}}

\newcommand{\shi}{\underline{\rm SHI}}

\newcommand{\khi}{\underline{\rm KHI}^-}
\newcommand{\khii}{\underline{\rm KHI}}

\newcommand{\deq}{\colonequals}

\newcommand{\al}{\alpha}
\newcommand{\be}{\beta}
\newcommand{\ga}{\gamma}
\newcommand{\Ga}{\Gamma}

\newcommand{\ot}{\otimes}%tensor product

\newcommand{\id}{\operatorname{Id}}

\newcommand{\p}{\prime}
\newcommand{\pp}{{\prime\prime}}

\newcommand{\aand}{~{\rm and}~}

\newcommand{\intg}{\mathbb{Z}}
\newcommand{\ft}{{\mathbb{F}_2}}

\newcommand{\ra}{\rightarrow}

\newcommand{\xra}{\xrightarrow}

\newcommand{\rgl}{\rangle}
\newcommand{\lgl}{\langle}

\DeclareMathOperator{\cok}{coker}
\DeclareMathOperator{\im}{Im}
\DeclareMathOperator{\ke}{ker}

\DeclareMathOperator{\cone}{Cone}

\newcommand{\sut}[1]{\mathbf{\Gamma}_{#1}}
\newcommand{\sutg}[2]{(\mathbf{\Gamma}_{#1},#2)}
\newcommand{\psp}[2]{\psi^{#1}_{+,#2}}
\newcommand{\psm}[2]{\psi^{#1}_{-,#2}}
\newcommand{\Psp}[2]{\Psi^{#1}_{+,#2}}
\newcommand{\Psm}[2]{\Psi^{#1}_{-,#2}}
\newcommand{\dehny}[1]{\mathbf{Y}_{#1}}

\newcommand{\sutb}[2]{\mathbf{\Gamma}^{#1}_{#2}}
\newcommand{\sutbg}[3]{(\mathbf{\Gamma}^{#1}_{#2},#3)}

\begin{document}

\title{Knot surgery formulae for instanton Floer homology II: applications}

%    Remove any unused author tags.

%    author one information
% \author{Zhenkun Li}
% \address{Department of Mathematics, Stanford University}
% \curraddr{}
% \email{zhenkun@stanford.edu}
% \thanks{}

\author{Zhenkun Li}
\address{Department of Mathematics and Statistics, University of South Florida}
\curraddr{}
\email{zhenkun@usf.edu}
\thanks{}

%    author two information
\author{Fan Ye}
\address{Department of Mathematics, Harvard University}
\curraddr{}
\email{fanye@math.harvard.edu}
\thanks{}

\date{}
\dedicatory{}
\begin{abstract}
This is a companion paper to earlier work of the authors, which proved an integral surgery formula for framed instanton homology. First, we present an enhancement of the large surgery formula, a rational surgery formula for null-homologous knots in any 3-manifold, and a formula encoding a large portion of $I^\sharp(S^3_0(K))$. Second, we use the integral surgery formula to study the framed instanton homology of many 3-manifolds: Seifert fibered spaces with nonzero orbifold degrees, especially nontrivial circle bundles over any orientable surface, surgeries on a family of alternating knots and all twisted Whitehead doubles, and splicings with twist knots. Finally, we use the previous techniques and computations to study almost L-space knots, \textit{i.e.}, the knots $K\subset S^3$ with $\dim I^\sharp(S_n^3(K))=n+2$ for some $n\in\mathbb{N}_+$. We show that an almost L-space knot of genus at least $2$ is fibered and strongly quasi-positive, and a genus-one almost L-space knot must be either the figure eight or the mirror of the $5_2$ knot in Rolfsen's knot table.
\end{abstract}
\maketitle

%\tableofcontents%table of contents
% \newpage

\tableofcontents
%————Start from here————

\section{Introduction}

Sutured instanton homology $SHI(M,\ga)$ for a balanced sutured manifold $(M,\ga)$ was introduced by Kronheimer-Mrowka \cite{kronheimer2011knot} and it leads to many important instanton invariants of $3$-manifolds and knots. Among them, the framed instanton homology $I^{\sharp}(Y)$ for a $3$-manifold $Y$ and the instanton knot invariant $KHI(Y,K)$ for a knot $K\subset Y$ are the two most important ones. It has been known that the framed instanton homology is closely related to the $SU(2)$-representations of the fundamental group $\pi_1(Y)$ and hence understanding its structural properties and computing its dimension are essential tasks in the study of instanton theory. However, the fact that instanton Floer homology is constructed based on a set of partial differential equations makes this task very difficult. Some previous computational results were obtained in \cite{scaduto2015instanton,lidman2020framed,baldwin2019lspace,baldwin2020concordance}.

Motivated by the conjecture proposed by Kronheimer-Mrowka \cite{kronheimer2010knots} that framed instanton homology and the hat version of Heegaard Floer homology are isomorphic to each other, and the known structural properties in the Heegaard Floer theory established by Ozsv\'ath-Szab\'o \cite{ozsvath2004holomorphicknot,Ozsvath2008integral,Ozsvath2011rational}, the authors of the current paper have established many structural properties that relate the framed instanton homology to instanton knot homology:

% was conjectured to be isomorphic to the sutured Floer homology $SFH(M,\ga)$. In particular, for a knot $K$ in a closed 3-manifold $Y$, it is conjectured that instanton knot homology $KHI(Y,K)$ is isomorphic to knot Floer homology $\widehat{HFK}(Y,K)$ and framed instanton homology $I^\sharp (Y)$ is isomorphic to Heegaard Floer homology $\widehat{HF}(Y)$. The construction of sutured instanton homology involves flat $SU(2)$-connections, which is equivalent to the homomorphisms from the fundamental group to $SU(2)$. In particular, from \cite[Theorem 4.6]{baldwin2018stein}, if $Y$ is a rational homology sphere with \begin{equation}\label{eq: inequality}
%   \dim I^\sharp(Y)>|H_1(Y;\mathbb{Z})|
%\end{equation}and $\pi_1(Y)$ is cyclically finite, then there exists a homomorphism $\rho:\pi_1(Y)\to SU(2)$ so that $\im \rho$ is nonabelian. From \cite[Proposition 4.9]{baldwin2018stein}, if $H_1(Y;\mathbb{Z})$ is cyclic and the order is a prime power, then $\pi_1(Y)$ is always cyclically finite. However, in general, it is hard to compute $I^\sharp(Y)$ by its definition; see  for some previous computations.

%Kronheimer-Mrowka's conjecture motivated the authors of this paper to establish many structural results for sutured instanton homology analogous to those in Heegaard Floer theory, which makes the computation feasible in many cases.
\benu
\item In \cite{LY2021,LY2021enhanced}, we established a decomposition of $SHI(M,\ga)$ along $H_1(M;\mathbb{Z})$, and showed that the enhanced Euler characteristic associated to this decomposition equals the Euler characteristic of $SFH(M,\ga)$ with respect to the spin$^c$ decomposition.
\item In \cite{LY2021large}, for a rationally null-homologous knot $K\subset Y$, we constructed two differentials, $d_+$ and $d_-$, on $KHI(Y,K)$ such that the homologies are isomorphic to $I^\sharp(Y)$. Using those differentials, we constructed some complexes called \textbf{bent complexes} whose homologies compute $I^\sharp(Y_n(K))$, where $Y_n(K)$ is obtained from $Y$ by Dehn surgery along $K$ with a large coefficient $n$. 
\item In \cite{LY2022integral1}, we established a formula based on the bent complexes that computes $I^\sharp(Y_m(K))$ for any nonzero integral $m$-surgery.
\eenu

Many applications have already been found based on this work: the proof that $\pi_1(S^3\backslash L)$ for almost all links $L$ admits an irreducible $SU(2)$-representation in \cite{yixiesu2}, the proof that $\pi_1(S^3_3(K))$ for any nontrivial knot admits an irreducible $SU(2)$-representation in \cite{BLSY21}, a strong restriction on the Alexander polynomial $\Delta_K(t)$ for any instanton L-space knot $K$ in \cite{LY2021large}, and the computation of $I^\sharp(S_r^3(K))$ for any genus-one quasi-alternating knot $K$ in \cite{LY2021large}, \textit{etc}.

In this paper, we present more applications of our previous work from (1) to (3), further demonstrating the power of these tools in dealing with the Dehn surgeries of knots: we upgrade the integral surgery formula proved in \cite{LY2022integral1} to a rational surgery formula; we study the $0$-surgery for knots inside $S^3$, which is a missing case in \cite{LY2022integral1}; we study almost L-space knots, which admit a surgery with next-to-minimal framed instanton homology, and we present the computations of many new families of the framed instanton homology of $3$-manifolds, including most Seifert fibered $3$-manifolds with non-zero orbifold degrees, the Dehn surgery along a large family of alternating knots and all twisted Whitehead doubles, and splicings of the complement of a twist knot with the complement of an arbitrary knot in $S^3$. Below, we give an outline of the contents of individual sections, providing more details of these results.

%The results in this paper fall into three parts, and most sections can be read separately. The first part (Section \ref{sec: preliminaries}-\ref{sec: Knots inside $S^3$}) is to generalize surgery formulae. This includes an enhancement of the large surgery formula, a rational surgery formula for null-homologous knots in any 3-manifold, and a formula encoding a large portion of $I^\sharp(S^3_0(K))$. The results are motivated by Ozsv\'ath-Szab\'o's work \cite{Ozsvath2008integral,Ozsvath2011rational}. The second part (Section \ref{sec: Surgeries on Borromean knots}-\ref{sec: twisted Whitehead doubles}) is to compute framed instanton homology by the integral surgery formula. The examples include Seifert fibered spaces with nonzero orbifold degrees, especially nontrivial circle bundles over surface, surgeries on a family of alternating knots and all twisted Whitehead doubles, and splicings with twist knots. The third part (Section \ref{sec: Almost L-space knots}) is to study almost L-space knots in $S^3$, which is based on the previous techniques and computations. Below, we give an outline of the contents of individual sections.

\quad

{\bf Section} \ref{sec: preliminaries}.
We review notations and results about surgery formulae in \cite{LY2021large,LY2022integral1}. We truncate the integral surgery formula to make it simpler for further usage. As a byproduct, we weaken the assumption on the large coefficient in the large surgery formula. In particular, when $K$ is null-homologous, the integer $2g(K)-1$ is large enough to apply the large surgery formula, while in \cite{LY2021large}, the minimal integer is $2g(K)+1$. 

\quad

{\bf Section} \ref{A rational surgery formula}. We establish a rational surgery formula for all null-homologous knots in instanton theory. The proof is similar to that in Heegaard Floer theory \cite{Ozsvath2011rational}. Suppose $K\subset Y$ is a null-homologous knot. The rational surgery along $K$ can be interpreted as the integral surgery along a knot $K_{\#}\subset Y\# L(p,q)$. The knot $K_{\#}$ is obtained by the connected sum of $K$ and a core knot in $L(p,q)$ (whose complement is a solid torus) for $(p,q)$ chosen according to the surgery slope. We establish a connected sum formula for the differentials in bent complexes in such cases and then apply the integral surgery formula to complete the proof.

\quad

{\bf Section} \ref{sec: Knots inside $S^3$}.
The statement of the integral surgery formula in \cite{LY2022integral1} excludes the case of $0$-surgery, \textit{i.e.}, the filling slope is the boundary of a Seifert surface. However, for a knot $K\subset S^3$, we can still understand a large portion of $I^\sharp(S_0^3(K))$ by examining an extra grading: after performing the $0$-surgery, the Seifert surface of $K$ is capped off by the meridian disk of the filling solid torus, which becomes an essential closed surface in $S^3_0(K)$. From \cite[Section 2.6]{baldwin2019lspace}, this surface induces a $\intg$-grading \begin{equation}\label{eq: z grading}
    I^{\sharp}(S^3_0(K))\cong \bigoplus_{s=1-g(K)}^{g(K)-1} I^{\sharp}(S^3_0(K),s).
\end{equation}
In this case, the integral surgery formula can be stated and proven grading-wise. As a result, we can understand $I^{\sharp}(S^3_0(K),s)$ for all $s$ but $0$.

\quad

The next three sections are about computations. To apply the integral surgery formula for a specific knot, there are two main tasks to solve:
\benu
    \item To compute differentials $d_\pm$ on $KHI(Y,K)$,
    \item To find the isomorphism $H(KHI(Y,K),d_+)\cong H(KHI(Y,K),d_-)$ in the statement of the surgery formula (\textit{cf.} Theorem \ref{thm: integral bent complex}).
\eenu

In the following three sections, we present many methods to deal with the above tasks (1) and (2).

\quad

\begin{figure}[ht]
	\begin{overpic}[width=0.5\textwidth]{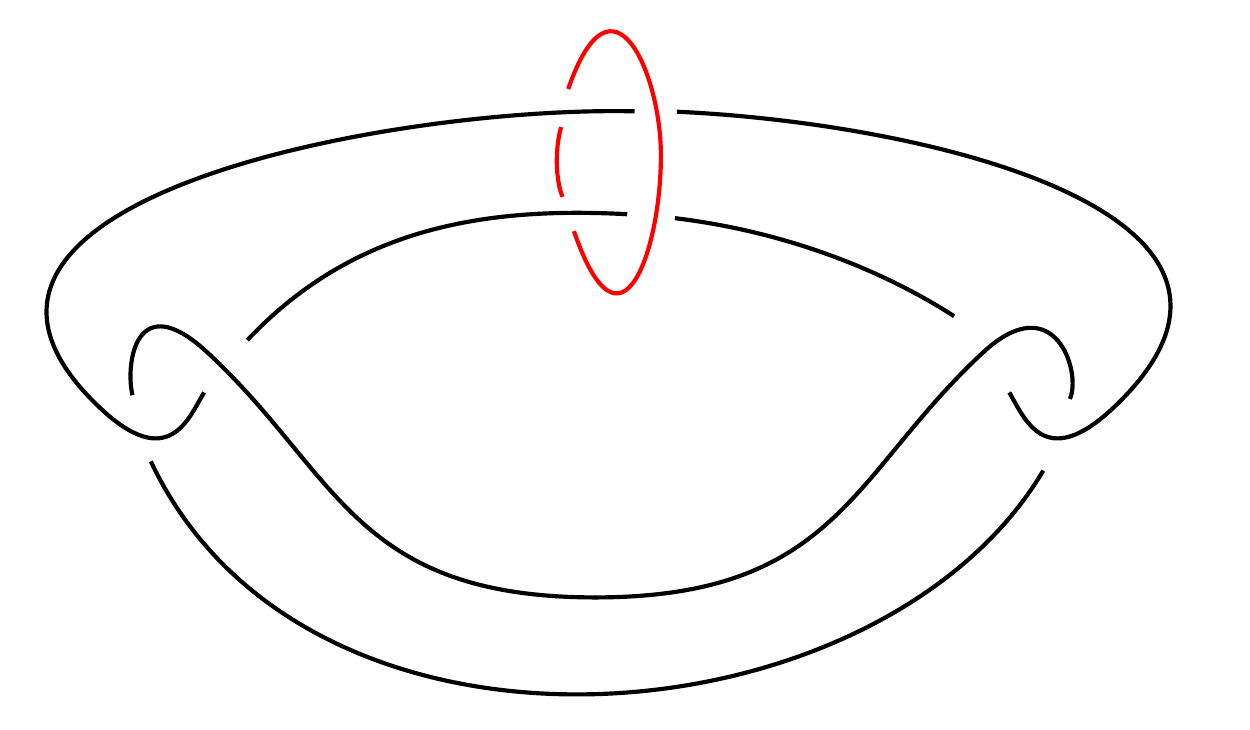}
		%\put(10,19){$\delta$}
		%\put(16,2){$K_+$}
		\put(47,30){\color{red} $K$}
		\put(80,46){$0$}
		\put(80,10){$0$}
	\end{overpic}
	\caption{The Borromean knot $K$ inside $S^1\times S^2\# S^1\times S^2$. The two copies of $S^1\times S^2$ come from the zero surgeries on the two (black) components of the Borromean link.\label{fig: borromean}}
\end{figure}

{\bf Section} \ref{sec: Surgeries on Borromean knots}. We deal with the Borromean knot as in Figure \ref{fig: borromean} and the connected sums of a few copies of them. Any such knot $K$ is inside a $3$-manifold $Y$ that is the connected sum of a few copies of $S^1\times S^2$. For this special family of knots, we have
$$KHI(Y,K)\cong I^{\sharp}(Y),$$
so task (1) is trivial. Moreover, the $H_1(Y)$-action in this case is essential. From \cite[Section 7.8]{scaduto2015instanton} and \cite[Theorem 7.16]{donaldson2002floer}, we have an identification
$$I^{\sharp}(Y)=\Lambda^*H_1(Y).$$
Hence we can regard all related instanton Floer homology groups as modules over the ring $\Lambda^*H_1(Y)$ and the task (2) can be done by the module structure.

It is worth mentioning that prior to the current paper, most computations of $I^\sharp(Y)$ are for rational homology spheres $Y$, while our computations for (connected sums of) Borromean knots, up to the author's knowledge, provide a first family of knots inside $3$-manifolds with positive first Betti number for which the framed instanton homology of their Dehn surgeries can be computed systematically. It is well known that the nonzero integral surgeries of connected sums of Borromean knots give nontrivial circle bundles over orientable surfaces. Hence, we obtain the following.
%surgery formula can compute surgeries on knots inside general $3$-manifolds. The second family of knots that is suitable to apply the integral surgery formula is the connected sums of Borromean knots inside connected sums of $S^1\times S^2$. The Borromean knot is the knot shown in Figure \ref{}. It is a knot in $\#^2 S^1\times S^2$ obtained from the three-component Borromean link in $S^3$ by performing $0$-surgeries along two components. We show that the differentials for the connected sums of Borromean knot are all trivial, which solves the first problem to apply the formula. For the second problem, though $I^\sharp (\# ^{2g}S^1\times S^2)$ is not $1$-dimensional, we can use the $\Lambda^* H_1(\# ^{2g}S^1\times S^2;\mathbb{C})$-action on $I^\sharp (\# ^{2g}S^1\times S^2)$ to give a strong restriction on the isomorphism between the complexes of $d_+$ and $d_-$. Indeed, we use the integral surgery formula in the sutured  setting (Theorem \ref{thm: integral surgery formula}), which also inherits the $\Lambda^*H_1$-action. Hence we can compute the isomorphism explicitly and then apply the integral surgery formula.

\bthm\label{thm: circle bundle computation}
For any $g>1$, $m\neq 0$, suppose $Y^g_m$ is the circle bundle over a surface of genus $g$ with Euler number $m$. We have the following.
\begin{enumerate}
	\item If $|m|\ge 2g-1$, then
	$$\dim I^{\sharp}(Y^g_{m})=2^{2g}\cdot |m|.$$
	\item If $|m|=2l$ with $l\leq g-1$, then
	$$\dim I^{\sharp}(Y_{m}^g)=2^{2g}\cdot |m|+4\cdot\sum_{j=1}^{g-l-1}\sum_{i=0}^{j-1}\binom{2g}{i}+2\cdot \sum_{i=0}^{g-l-1}\binom{2g}{i}.$$
	\item If $|m|=2l-1$ with $l\leq g-1$, then
	$$\dim I^{\sharp}(Y^g_{m})=2^{2g}\cdot |m|+4\cdot\sum_{j=1}^{g-l}\sum_{i=0}^{j-1}\binom{2g}{i}$$
\end{enumerate}
\ethm
\brem
In \cite[Theorem 5.5]{Ozsvath2008integral}, Ozsv\'ath-Szab\'o provided a formula for $HF^+_{red}(Y^g_m)$ using the integral surgery formula for $HF^+$.
\erem
Furthermore, we can recover any Seifert fibered space with nonzero orbifold degree by a non-zero integral surgery along the connected sum of Borromean knots and suitable core knots in lens spaces. We also use the $\Lambda^*H_1(Y)$-module structure to solve task (2). As a result, we prove the following theorem, which generalizes Alfieri-Baldwin-Dai-Sivek's result for Seifert fibered manifolds that are rational homology spheres \cite[Corollary 1.3]{alfieri2020framed}.

\bthm\label{thm: SF space}
Let $Y$ be a Seifert fibered space over a genus $g$ orbifold with Seifert invariants $(m,r_1/v_1,\dots,r_n/v_n)$. Suppose the orbifold degree is $$\deg Y=m+\sum_{i=1}^n\frac{r_1}{v_1}.$$If $\deg Y\neq 0$ and $\gcd(v_i,v_j)=1$ for any $i\neq j\in\{1,\dots,n\}$, then
$$\dim_\mathbb{C}I^{\sharp}(Y)=\dim_{\ft}\widehat{HF}(Y).$$
\ethm
\brem\label{rem: remove condition}
It is possible to compute $\dim_\mathbb{C} I^\sharp(Y)$ in Theorem \ref{thm: SF space} explicitly as in \cite[Theorem 10.1]{Ozsvath2011rational}. We need the condition $\deg Y\neq 0$ because we do not have a zero-surgery formula for knots inside general manifolds and $\deg Y=0$ corresponds to the zero-surgery on the connected sum. We need the condition $\gcd(v_i,v_j)=1$ because we want the first homology of the complement of the connected sum to be torsion-free, so that we can use the grading from the Seifert surface to capture all information in the spin$^c$ decomposition of Heegaard Floer theory. This condition could be removed if we utilize the work in \cite{LY2021enhanced} to obtain a further composition of our integral surgery formula.
\erem

\quad

{\bf Section} \ref{sec: Surgeries on some alternating knots}. We also study more families of knots inside $S^3$. Since there are isomorphisms
$$H(KHI(S^3,K),d_+)\cong H(KHI(S^3,K),d_-)\cong I^{\sharp}(S^3)\cong \mathbb{C},$$
the choice of the isomorphism between them is only up to a scalar. Hence task (2) is trivial, and all we need is to deal with task (1). 

It is well known that alternating knots are thin in the Heegaard Floer theory \cite{Ozsvath2003}. From Petkova's classification of thin complexes \cite[Section 3.1]{Petkova2009thin}, the knot Floer complex of an alternating knot is fully determined by its Alexander polynomial and the tau invariant (which is related to the signature for alternating knots). Since there is no known integral Maslov grading in instanton theory, we do not have a proper definition of thin knots in the instanton setting. 

Instead, we can consider knots whose two spectral sequences from $KHI(S^3,K)$ to $I^\sharp(S^3)$ collapse on the second pages, \textit{i.e.}, only differentials $d_{1,\pm}$ are nontrivial. We say such knots have \textbf{torsion order one} (\textit{cf.} Definition \ref{defn: torsion order one}). For knots having torsion order one, we have a similar classification of complexes involving $d_\pm$ as the thin complexes, and hence the complexes are again fully determined by the Alexander polynomial and the tau invariant in instanton theory. 

In order to prove that some families of knots have torsion order one, we make use of the oriented skein relation in instanton theory studied in \cite{Lim2009,kronheimer2010instanton}. Unlike the original setup, where we have an oriented smoothing of the crossing to derive a link in $S^3$, we consider its knotification, or equivalently a knot inside $S^1\times S^2$. 

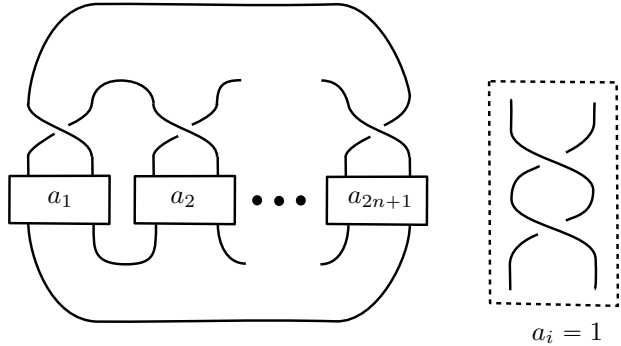
\begin{figure}[ht]
	\begin{overpic}[width=0.6\textwidth]{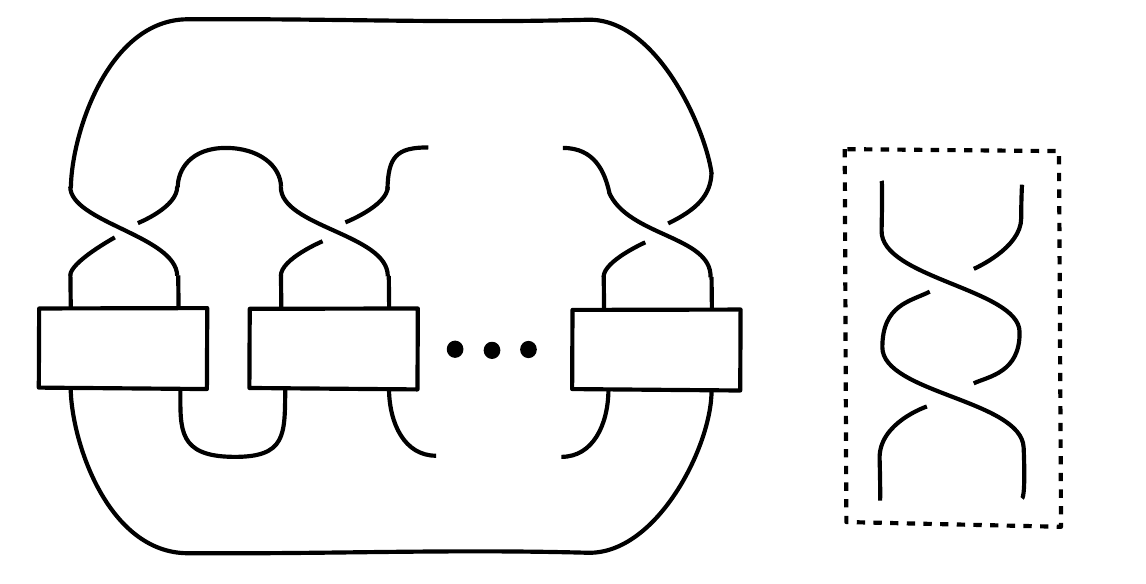}
		\put(9,20){$a_1$}
		\put(27,20){$a_2$}
		\put(53,20){$a_{2n+1}$}
		%\put(56,39){$\delta$}
		\put(80,0){$a_i=1$}
	\end{overpic}
	\caption{The knot $K(a_1,\dots,a_{2n+1})$.}\label{fig: alternating1}
\end{figure}

This idea of using oriented skein relation works for a large family of alternating knots. In particular, we can deal with the family of pretzel knots as shown in Figure \ref{fig: alternating1}. Note that all the crossings in this family of diagrams are positive, since the induction starts with the torus knots $T(2,2n+1)$ (\textit{i.e.} $a_i=0$ for all $i$), whose crossings are all positive. We prove that those knots have torsion order one, and then we can compute the framed instanton homologies on their surgeries.

%We use the oriented skein relation for $KHI$ studied in  has torsion order one. We start with knots with known differentials (the torus knots $T(2,2n+1)$) and compute the differentials for new knots by induction. During the induction, we need to deal with knots in $S^1\times S^2$, for which we can also discuss the torsion order one property. Finally, we have the following theorem.

%We also compute differentials for a family of alternating knots so that we can apply the integral/rational surgery formula or the large surgery formula plus the dimension formula in \cite{baldwin2020concordance} to compute $I^\sharp$ of any nonzero surgery. The family of knots is shown in Figure \ref{} and denoted by $K(a_1,\dots,a_{2n+1})$, where $a_1,\dots, a_{2n+1}$ are all odd numbers of half-twists. 

\bthm\label{thm: alternating knot}
Suppose $K\subset S^3$ is a knot as shown in Figure \ref{fig: alternating1} so that $$a_{i}\geq 0\text{ for all }i=1,\dots,2n+1, \text{ and }$$
$$\#\{i~|~a_i\ge 1\}\leq n+1.$$
Then we have $g(K)=n$ and for any $r=p/q\in\mathbb{Q}\backslash\{0\}$ with $q\ge 1$, we have
$$\dim I^{\sharp}(S^3_r(K))=\dim_{\ft}\widehat{HF}(S^3_r(K))=
\frac{(||\Delta_K(t)||-2n-3)\cdot q}{2}+ |p-q\cdot (2n-1)|,$$where $||\cdot||$ is the sum of the absolute values of the coefficients.
\ethm

\brem\label{rem: remove condition2}
The one-dimensional argument that task (2) is trivial for knots inside $S^3$ can also be generalized to a knot $K$ in any instanton L-space $Y$. If $H_1(Y\backslash N(K);\mathbb{Z})$ is torsion-free, then we may use the grading from the Seifert surface to decompose our integral surgery formula, so that the one-dimensional argument can be applied to each summand. If $H_1(Y\backslash N(K);\mathbb{Z})$ is not torsion-free, we could utilize the work in \cite{LY2021enhanced} to obtain a decomposition, but that needs further study concerning the interaction of the decomposition and the construction of the integral surgery formula (\textit{cf.} Remark \ref{rem: remove condition}). 
\erem

\quad

{\bf Section} \ref{sec: twisted Whitehead doubles}. We also use the techniques involving oriented skein relation to study twisted Whitehead doubles.

\begin{figure}[ht]
	\begin{overpic}[width=0.7\textwidth]{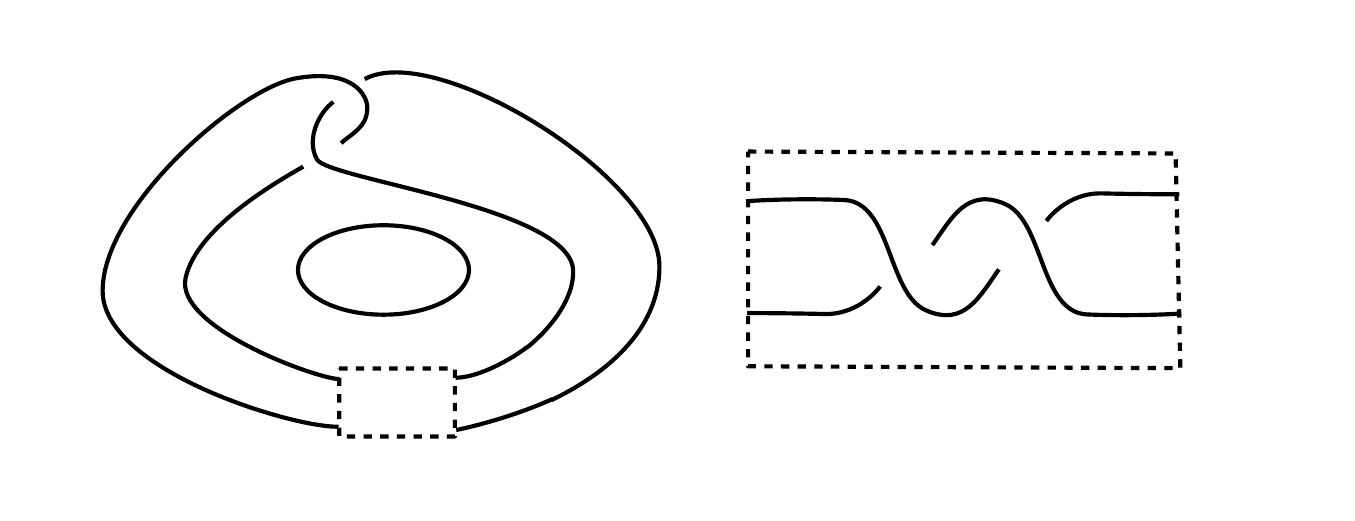}
		%\put(10,16){$\widehat{K}$}
		%\put(-2,16){$\widehat{V}$}
		\put(28.5,8){$n$}
		\put(68,7){$n=1$}
	\end{overpic}
	\caption{Whitehead double.\label{fig: pattern of WHD 2}}
\end{figure}

\bthm\label{thm: Whitehead double}
Suppose $J\subset S^3$ is a knot and $K=D^+_t(J)$ is the $t$-twisted Whitehead double of $J$ with positive clasp; see Figure \ref{fig: pattern of WHD 2}. Suppose $\tau_I$ is the instanton tau invariant and $\Ga_{n}\subset \partial (S^3\backslash N(J))$ consists of two copies of curves with slope $-n$. Then we have the following.
\begin{enumerate}
	\item $KHI(S^3,K,1)\cong SHI(S^3\backslash N(J),\Gamma_{-t}).$
	\item $\tau_I(K)=\begin{cases}
		1&t<2\cdot \tau_I(J)\\
		0&t\geq 2\cdot \tau_I(J) \end{cases}$
	\item $
	\dim I^{\sharp}(S^3_{\pm1}(K))=\begin{cases}
		2\cdot\dim SHI(S^3\backslash N(J),\Ga_{-t})\mp1&t<2\cdot \tau_I(J)\\
		2\cdot\dim SHI(S^3\backslash N(J),\Ga_{-t})+1&t\geq 2\cdot \tau_I(J)
	\end{cases}$
% 	\item If $m<2\cdot \tau_I(J)$, then we have
% 	$$\tau_I(K)=\nu^{\sharp}(K)=1,~{\rm and}\dim I^{\sharp}(S^3_1(K))=2 \cdot \dim SHI(S^3\backslash N(J),\Gamma_{-t})-1.$$
% 	\item If $m\geq 2\cdot \tau_I(J)$, then we have
% 	$$\tau_I(K)=\nu^{\sharp}(K)=0,~{\rm and}\dim I^{\sharp}(S^3_1(K))=2\cdot\dim SHI(S^3\backslash N(J),\Gamma_{-t})+1.$$
\end{enumerate} 
\ethm

\brem
According to \cite[Theorem 1.1]{baldwin2020concordance}, the data provided in Theorem \ref{thm: Whitehead double} part (3) is enough to compute the framed instanton homology of all nonzero rational surgeries of the twisted Whitehead doubles with positive clasps. Also, note that we have
$$\overline{D^-_m(K)}=D^+_{-m}(\widebar{K}),$$
where $\widebar{K}$ is the mirror of $K$. So we also know all the information for twisted Whitehead doubles with negative clasps. 
\erem

Theorem \ref{thm: Whitehead double} can also be applied to study splicings with knot complements of twist knots. Note that twist knots $K_n$ are the positively clasped $n$-twisted Whitehead doubles of the unknot.

\bthm\label{thm: splicing with twist knot}
Suppose $K_n$ is the twist knot. Suppose $J\subset S^3$ is a non-trivial knot. Let $Y$ be obtained by gluing the complement of $K_n$ with the complement of $J$ so that the gluing map sends the meridian of one knot to the longitude of the other and vice versa. Let $\Ga_0\subset\partial (S^3\backslash N(J))$ consist of two Seifert longitudes. Then
\begin{equation*}
	\dim I^{\sharp}(Y)=\begin{cases}
		2\cdot|n|\cdot\dim SHI(S^3\backslash N(J),\Ga_0)+1&\tau_{I}(J)\leq 0\\
		|n|\cdot(2\cdot\dim SHI(S^3\backslash N(J),\Ga_0)-1)+|1+n|&\tau_I(J)>0
	\end{cases}
\end{equation*}
\ethm

\brem
From \cite[Section 5.2]{li2019tau}, we have the following equality for $n\in\intg$ ({\it cf.} Lemma \ref{lem: symmetry and dimension of Gamma_n})
$$\dim SHI(S^3\backslash N(J),\Gamma_n)=\dim SHI(S^3\backslash N(J),\Gamma_{-2\tau_I(K)})+|n+2\tau_I(K)|.$$
So, for a knot $K\subset S^3$, as long as we know its $\tau_I$ and  $\dim SHI(S^3\backslash N(J),\Gamma_n)$ for any one $n\in\intg$, we can obtain the dimensions for all $n\in\intg$. Furthermore, from Theorem \ref{thm: Whitehead double} and Theorem \ref{thm: splicing with twist knot}, we obtain the framed instanton homology of Dehn surgeries along all of their twisted Whitehead doubles as well as the splicing with the complements of the twist knots. Here is the list of knots where all such data are known.
\begin{itemize}
	\item Genus-one quasi-alternating knots ({\it cf.} \cite[Section 6]{LY2021large}).
	\item Instanton L-space knots ({\it cf.} \cite[Section 5]{LY2021large}).
	\item Knots described in Theorem \ref{thm: alternating knot} ({\it cf.} Section \ref{sec: Surgeries on some alternating knots}).
\end{itemize}
\erem

\quad

{\bf Section} \ref{sec: Almost L-space knots}.
Finally, we study almost L-space knots in $S^3$. A knot $K\subset S^3$ is called an {\bf almost (instanton) L-space knot} if it is not an instanton L-space knot and there exists $n\in\mathbb{N}_+$ such that
\begin{equation}\label{eq: defn almost l}
    \dim I^{\sharp}(S^3_n(K))=n+2.
\end{equation}
Note that $n+2$ is the second minimal value of the dimension since the Euler characteristic is $n$ \cite[Corollary 1.4]{scaduto2015instanton}. See \cite{baldwin22knot52} for the results in Heegaard Floer theory.

Similar to the previous work on instanton L-space knots \cite{LY2021large}, we can impose strong restrictions on almost L-space knots. Moreover, we can classify all genus-one almost L-space knots.

\bthm\label{thm: almost L space knots, intro}
Suppose $K\subset S^3$ is an almost L-space knot. Then we have the following.
\begin{enumerate}
	\item If $g(K)\geq 2$, then $K$ is fibered, strongly quasi-positive, and $\tau_I(K)=g(K)$.
	\item If $g(K)=1$, then $K$ is either the figure eight or the mirror of the $5_2$ knot in Rolfsen's knot table (with signature $-4$, denoted by $\widebar{5}_2$).
\end{enumerate}
\ethm
A direct corollary of Theorem \ref{thm: almost L space knots, intro} is the following.
\bcor\label{cor: 1 surgery having dimension 3}
Suppose $K\subset S^3$ is a knot. Suppose further that 
$$\dim I^{\sharp}(S^3_1(K))=3.$$ 
Then $K$ is either the left-handed trefoil, the figure-eight, or the knot $\widebar{5}_2$.
\ecor

The proof for $g(K)\ge 2$ largely depends on our previous work in \cite[Section 5]{LY2021large}. The classification of genus-one almost L-space knots is more complicated. We first proved that $KHI(S^3,K)$ is $1$- or $2$-dimensional in the top Alexander grading, for which we know a list of all possible knots. If the top grading is $1$-dimensional, then the knot is fibered \cite[Corollary 7.19]{kronheimer2010knots}. It is well known that the trefoil and the figure eight are the only genus-one fibered knots. If the top grading is $2$-dimensional, Baldwin-Sivek \cite{BS2022nearly} recently classified all such knots in the Heegaard Floer setting. According to \cite{LY2022nearly}, the same results apply to the instanton setting. This also leads to the following theorem, which is a complete classification of genus-one nearly fibered knots in terms of instanton knot homology.
\bthm\label{thm: genus-one nearly fibred}
Suppose $K\subset S^3$ is a genus-one knot with
$$\dim KHI(S^3,K,1)=2.$$
Let $J$ be the right-handed trefoil. Then we know the following.
\begin{enumerate}
	\item $K$ is $5_2$ or its mirror if and only if 
	$$\dim KHI(S^3,K)=7.$$
	\item $K$ is, up to mirror, either $15n_{43522}$ or $D^-_2({J})$ if and only if
	$$\dim KHI(S^3,K)=9~{\rm and}~\Delta_{K}(t)=2t-3+2t^{-1}.$$
	\item $K$ is one of the pretzel knots $P(-3,3,2n+1)$ for some $n\in\intg$, $D^+_2({J})$, or their mirrors if and only if
	$$\dim KHI(S^3,K)=9~{\rm and}~\Delta_{K}(t)=-2t+5-2t^{-1}.$$
\end{enumerate}
\ethm

\brem
Prior to the computation in this paper, due to Baldwin-Sivek's work \cite{BS2022nearly}, we know that if $K$ is genus-one and
$\dim KHI(S^3,K,1)=2,$
then $K$ must be one of the knots listed in Theorem \ref{thm: genus-one nearly fibred}.
Furthermore, we already know that $\dim KHI(S^3,K)=7$ for $K=5_2$ and $\dim KHI(S^3,K)=9$ for $K=P(-3,3,2n+1)$ and $K=D^+_2({J})$. The last piece for the above complete classification is the computations for $D^-_2({J})$ and $15n_{43522}$. This is finished in Section \ref{sec: twisted Whitehead doubles} and Section \ref{sec: Almost L-space knots}, respectively, by studying their Dehn surgeries.
\erem

\begin{ack}
The authors thank John A. Baldwin and Steven Sivek for introducing us to almost L-space knots. The authors also thank Ciprian Manolescu, Thomas E. Mark, Tomasz Mrowka, and Jacob Rasmussen for helpful comments and valuable discussions.
\end{ack}

\section{Preliminaries on surgery formulae}\label{sec: preliminaries}

\subsection{Conventions}

If it is not mentioned, all manifolds are smooth, oriented, and connected. Homology groups and cohomology groups are defined with $\mathbb{Z}$ coefficients. We write $\mathbb{Z}_n$ for $\mathbb{Z}/n\mathbb{Z}$ and $\mathbb{F}_2$ for the field with two elements. If there is no subscript for $\dim$, then it means $\dim_\mathbb{C}$.

A knot $K\subset Y$ is called \textbf{null-homologous} if it represents the trivial homology class in $H_1(Y;\mathbb{Z})$, while it is called \textbf{rationally null-homologous} if it represents the trivial homology class in $H_1(Y;\mathbb{Q})$.

For any oriented 3-manifold $M$, we write $-M$ for the manifold obtained from $M$ by reversing the orientation. For any surface $S$ in $M$ and any suture $\ga\subset \partial M$, we write $S$ and $\ga$ for the same surface and suture in $-M$, without reversing their orientations. For a knot $K$ in a 3-manifold $Y$, we write $(-Y,K)$ for the induced knot in $-Y$ with induced orientation, called the \textbf{mirror knot} of $K$. The corresponding balanced sutured manifold is $(-Y\backslash N(K),-\ga_K)$.

\subsection{Sutured instanton homology for knot complements}\label{sec: SHI, 1}

For any \textbf{balanced sutured manifold} $(M,\ga)$ \cite[Definition 2.2]{juhasz2006holomorphic}, Kronheimer-Mrowka \cite[Section 7]{kronheimer2010knots} constructed an isomorphism class of $\mathbb{C}$-vector spaces $SHI(M,\ga)$. Later, Baldwin-Sivek \cite[Section 9]{baldwin2015naturality} dealt with the naturality issue and constructed (untwisted and twisted versions of) projectively transitive systems related to $SHI(M,\ga)$. We will use the twisted version, which we write as $\shi(M,\ga)$ and call \textbf{sutured instanton homology}. 

Moreover, there is a relative $\mathbb{Z}_2$-grading on $\shi(M,\ga)$ obtained from the construction of sutured instanton homology, which we consider as a \textbf{homological grading} and use to take the Euler characteristic.

\bdefn
Suppose $K$ is a knot in a closed 3-manifold $Y$. Let $Y(1)\deq Y\backslash B^3$ and let $\delta$ be a simple closed curve on $\partial Y(1)\cong S^2$. Let $Y\backslash N(K)$ be the knot complement and let $\Ga_{\mu}$ be two oppositely oriented meridians of $K$ on $\partial (Y\backslash N(K))\cong T^2$. Define\[I^\sharp(Y)\deq \shi(Y(1),\delta)\aand \khii(Y,K)\deq \shi(Y\backslash N(K),\Ga_{\mu}).\]
\edefn
% \brem\label{rem: basept}
% By the naturality results, we should specify the places of the removing ball, the neighborhood of the knot, and the sutures to define $I^\sharp(Y)$ and $\khii(Y,K)$. These data can be fixed by choosing a basepoint in $Y$ or $K$. For simplicity, we omit those choices in the notations.
% \erem
From now on, we will suppose $K\subset Y$ is a rationally null-homologous knot and fix some notations. Let $\mu$ be the meridian of $K$ and pick a longitude $\lambda$ (such that $\lambda\cdot \mu=1$) to fix a framing of $K$.  We will always assume $Y\backslash N(K)$ is irreducible, but many results still hold due to the following connected sum formula of sutured instanton homology \cite[Section 1.8]{li2018contact}:
$$\shi(Y^\p\# Y\backslash N(K),\ga)\cong I^\sharp(Y^\p)\otimes \shi(Y\backslash N(K),\ga).$$

Given coprime integers $r$ and $s$, let $\Gamma_{r/s}$ be the suture on $\partial (Y\backslash N(K))$ that consists of two oppositely oriented, simple closed curves of slope $-r/s$, with respect to the chosen framing (\textit{i.e.} the homology of the curves are $\pm(-r\mu+s\lambda)\in H_1(\partial N(K))$). Pick $S$ to be a minimal genus Seifert surface of $K$.
\begin{conv}
We will use $p$ to denote the order of $[K]\in H_1(Y)$, \textit{i.e.}, $p$ is the minimal positive integer satisfying $p[K]=0\in H_1(Y)$. Let $q=\partial S\cdot \lambda$ and let $g=g(S)$ be the genus of $S$. When $K$ is null-homologous, we always choose the Seifert framing $\lambda=\partial S$. In such a case, we have $(p,q)=(1,0)$. 
\end{conv}
\brem\label{rem: notations}
The meanings of $p$ and $q$ follow from \cite{LY2022integral1}, but are different from our previous papers \cite{LY2020,LY2021large}. Before, we used $\hat{\mu}$ and $\hat{\lambda}$ to denote the meridian of the knot $K$ and its preferred framing. When $\partial S$ is connected, it is the same as the homological longitude $\lambda$ in previous papers. Hence $p$ and $q$ in this paper should be $q$ and $q_0$ in previous papers.
\erem

For simplicity, we use the bold symbol of the suture to represent the sutured instanton homology of the balanced sutured manifold with the reversed orientation:
$$\mathbf{\Gamma}_{r/s}\deq \shi(-(Y\backslash N(K)),-\Gamma_{r/s}).$$
When $(r,s)=(\pm 1,0)$, we write $\Ga_{r/s}=\Ga_{\mu}$. When $s=\pm 1$, we write $\Ga_{n}=\Ga_{n/1}=\Ga_{(-n)/(-1)}$. We also write $\mathbf{\Gamma}_{\mu}$ and $\mathbf{\Gamma}_{n}$ for the corresponding sutured instanton homologies.

Also, we write
$$\mathbf{Y}_{r/s}\deq I^{\sharp}(-Y_{-r/s}(K)),$$
and in particular
$$\mathbf{Y}_n\deq I^{\sharp}(-Y_{-n}(K))\aand \dehny{}\deq I^{\sharp}(-Y).$$
% \brem\label{rem: ambiguity}
% Strictly speaking, the sutures corresponding to $(r,s)=(1,0)$ and $(-1,0)$ are not identical because the orientations are opposite. Since both sutures are on $\partial (Y\backslash N(K))$ of the same slope, they are isotopic. Moreover, we can choose a canonical isotopy by rotating the suture along the direction specified by the orientation of the knot. Due to discussion in Heegaard Floer theory \cite{sarkar15moving,Zemke2019} and the conjectural relation between Heegaard Floer theory and instanton theory \cite{kronheimer2010knots}, it is expected that rotating the suture back to the original position induces a nontrivial isomorphism of the sutured instanton homology. So we pick the canonical isotopy to be the minimal rotation of the suture. Hence we can abuse notations and write $\Ga_\mu$ for both sutures. The same discussion also applies to the relation between $\Ga_{n/1}$ and $\Ga_{(-n)/(-1)}$.
% \erem

We always assume that $S$ has minimal intersection with $\Ga_{r/s}$. By work of \cite{li2019direct}, the Seifert surface $S$ induces either a $\mathbb{Z}$-grading or a $(\mathbb{Z}+\frac{1}{2})$-grading on $\mathbf{\Gamma}_{r/s}$, depending on the parity of the intersection number $\partial S\cdot(s\lambda-r\mu)$. We write the graded part of $\sut{r/s}$ as
$$\sutg{r/s}{i}\deq \shi(-(Y\backslash N(K)),-\Gamma_{r/s},S,i)$$
with $i\in\intg$ or $i\in\intg+\frac{1}{2}$, depending on the parity of the intersection number.

For simplicity, we omit the definitions  of bypass maps $\psi_{\pm,*}^*$ and surgery maps $F_n$, $G_n$, $H_n$, $A_{n-1}$, $B_{n-1}$, $C_n$ in \cite[Section 2.2]{LY2022integral1} and only list their properties as follows. The proofs and references can be found in \cite[Section 2.2]{LY2022integral1}.

\blem\label{lem: vanishing grading}
We have $\sutg{r/s}{i}=0$ when $$|i|>g+\frac{|rp-sq|-1}{2}.$$
\elem

\blem\label{lem: bypass n,n+1,mu}
For any $n\in\intg$, there are two graded bypass exact triangles
\begin{equation*}
\xymatrix@R6ex{
\sutg{n}{i+\frac{p}{2}}\ar[rr]^{\psp{n}{n+1}}&&\sutg{n+1}{i}\ar[dl]^{\psp{n+1}{\mu}}\\
&\sutg{\mu}{i-\frac{np-q}{2}}\ar[ul]^{\psp{\mu}{n}}&
}	
\end{equation*}
\begin{equation*}
\xymatrix@R6ex{
\sutg{n}{i-\frac{p}{2}}\ar[rr]^{\psm{n}{n+1}}&&\sutg{n+1}{i}\ar[dl]^{\psm{n+1}{\mu}}\\
&\sutg{\mu}{i+\frac{np-q}{2}}\ar[ul]^{\psm{\mu}{n}}&
}	
\end{equation*}
where the maps are homogeneous with respect to the homological $\mathbb{Z}_2$-gradings.  
\elem
% \brem
% The reason to use balanced sutured manifold with reversed orientation is because of the above bypass exact triangles.
% \erem
% \brem
%  Note that if we do not mention gradings, the above triangles and the results in the rest of this subsection also hold without the assumption that $K$ is rationally-null-homologous since the proofs only involve the neighborhood of $\partial(-Y\backslash N(K))$.
% \erem

% \bcor\label{cor: psi^n_+,n+1 is an iso}
% For any large enough integer $n$, we have the following properties for maps.
% \begin{enumerate}
% 	\item The map $\psp{n}{n+1}|{\sutg{n}{i}}$ is an isomorphism when $i\le \frac{np-q-1}{2}-g$.
% 	\item The map $\psm{n}{n+1}|{\sutg{n}{i}}$ is an isomorphism when $i\ge g-\frac{np-q-1}{2}$.
% 	\item For any $g-\frac{np-q-1}{2}\le i\le\frac{np-q-1}{2}-g-p$, there is an isomorphism
% 	$$(\psp{n}{n+1})^{-1}\circ\psm{n}{n+1}:\sutg{n}{i}\xra{\cong}\sutg{n}{i+p}.$$
% 	\item The map $\psm{-n}{1-n}|{\sutg{-n}{i}}$ is an isomorphism when $i\le \frac{(n-1)p+q-1}{2}-g$.
% 	\item The map $\psp{-n}{1-n}|{\sutg{-n}{i}}$ is an isomorphism when $i\ge g-\frac{(n-1)p+q-1}{2}$.
% 	\item For any $g-\frac{np+q-1}{2}\le i\le\frac{np+q+1}{2}-g-p$, there is an isomorphism
% 	$$(\psp{-n}{1-n})^{-1}\circ\psm{-n}{1-n}:\sutg{-n}{i}\xra{\cong}\sutg{-n}{i+p}.$$

% \end{enumerate}
% \ecor
% \bpf
% It is a combination of Lemma \ref{lem: vanishing grading} and Lemma \ref{lem: bypass n,n+1,mu}.
% \epf
\bdefn
The maps in Lemma \ref{lem: bypass n,n+1,mu} are called \textbf{bypass maps}. The ones with subscripts $+$ and $-$ are called \textbf{positive} and \textbf{negative bypass maps}, respectively. We will use $\pm$ to denote either of the bypass maps. For any integer $n$ and any positive integer $k$, define $$\Psi_{\pm,n+k}^n\deq \psi_{\pm,n+k}^{n+k-1}\circ\cdots\circ \psi_{\pm,n+1}^n:\mathbf{\Gamma}_{n}\rightarrow \mathbf{\Gamma}_{n+k}.$$
\edefn

\begin{lem}\label{lem: com diag for n,n+1,n+2}
For any $n\in\intg$, we have the following commutative diagrams up to scalars.
\begin{equation*}
	\xymatrix{
	\sut{n}\ar[rr]^{\psp{n}{n+1}}\ar[dd]_{\psm{n}{n+1}}&&\sut{n+1}\ar[dd]^{\psm{n+1}{n+2}}\\
	&&\\
	\sut{n+1}\ar[rr]^{\psp{n+1}{n+2}}&&\sut{n+2}
	}
	\xymatrix{
	\sut{n+2}\ar[rr]^{\psp{n+2}{\mu}}\ar[dd]_{\psm{n+2}{\mu}}&&\sut{\mu}\ar[dd]^{\psp{\mu}{n}}\\
	&&\\
	\sut{\mu}\ar[rr]^{\psm{\mu}{n}}&&\sut{n}
	}
\end{equation*}
\end{lem}

\begin{lem}\label{lem: comm diag for n,n+1,mu}
For any $n\in\intg$, we have the following commutative diagrams up to scalars
\begin{center}
\begin{minipage}{0.4\textwidth}
	\begin{equation*}
\xymatrix{
\sut{n}\ar[rr]^{\psp{n}{n+1}}&&\sut{n+1}\\
&\sut{\mu}\ar[lu]^{\psm{\mu}{n}}\ar[ru]_{\psm{\mu}{n+1}}&
}	
\end{equation*}
\end{minipage}
\begin{minipage}{0.4\textwidth}
	\begin{equation*}
\xymatrix{
\sut{n}\ar[rr]^{\psm{n}{n+1}}&&\sut{n+1}\\
&\sut{\mu}\ar[lu]^{\psp{\mu}{n}}\ar[ru]_{\psp{\mu}{n+1}}&
}	
\end{equation*}
\end{minipage}
\end{center}

\begin{center}
\begin{minipage}{0.4\textwidth}
	\begin{equation*}
\xymatrix{
\sut{n}\ar[rr]^{\psp{n}{n+1}}\ar[dr]_{\psm{n}{\mu}}&&\sut{n+1}\ar[dl]^{\psm{n+1}{\mu}}\\
&\sut{\mu}&
}	
\end{equation*}
\end{minipage}
\begin{minipage}{0.4\textwidth}
		\begin{equation*}
\xymatrix{
\sut{n}\ar[rr]^{\psm{n}{n+1}}\ar[dr]_{\psp{n}{\mu}}&&\sut{n+1}\ar[dl]^{\psp{n+1}{\mu}}\\
&\sut{\mu}&
}	
\end{equation*}\end{minipage}
\end{center}
\end{lem}

\blem\label{lem: bypass n+1,n,2n+1/2}
For a knot $K\subset Y$ and $n\in\intg$, there are two graded bypass exact triangles
\begin{equation*}
\xymatrix{
\sutg{n-1}{i+\frac{np-q}{2}}\ar[rr]^{\psp{n-1}{\frac{2n-1}{2}}}&&\sutg{\frac{2n-1}{2}}{i}\ar[dl]^{\psp{\frac{2n-1}{2}}{n}}\\
&\sutg{n}{i-\frac{(n-1)p-q}{2}}\ar[ul]^{\psp{n}{n-1}}&
}	
\end{equation*}
\begin{equation*}
\xymatrix{
\sutg{n-1}{i-\frac{np-q}{2}}\ar[rr]^{\psm{n-1}{\frac{2n-1}{2}}}&&\sutg{\frac{2n-1}{2}}{i}\ar[dl]^{\psm{\frac{2n-1}{2}}{n}}\\
&\sutg{n}{i+\frac{(n-1)p-q}{2}}\ar[ul]^{\psm{n}{n-1}}&
}	
\end{equation*}
\elem

% \bcor\label{cor: n+1,n,2n+1/2 is an iso}
% For any large enough integer $n$, we have the following properties for restrictions of maps.
% \begin{enumerate}
% 	\item The map $\psp{n}{n+1}|{\sutg{n}{i}}$ is an isomorphism when $i\le \frac{np-q-1}{2}-g$.
% 	\item The map $\psm{n}{n+1}|{\sutg{n}{i}}$ is an isomorphism when $i\ge g-\frac{np-q-1}{2}$.
% 	\item For any $g-\frac{np-q-1}{2}\le i\le\frac{np-q-1}{2}-g-p$, there is an isomorphism
% 	$$(\psp{n}{n+1})^{-1}\circ\psm{n}{n+1}:\sutg{n}{i}\xra{\cong}\sutg{n}{i+p}.$$
% 	\item The map $\psm{-n}{1-n}|{\sutg{-n}{i}}$ is an isomorphism when $i\le \frac{(n-1)p+q-1}{2}-g$.
% 	\item The map $\psp{-n}{1-n}|{\sutg{-n}{i}}$ is an isomorphism when $i\ge g-\frac{(n-1)p+q-1}{2}$.
% 	\item For any $g-\frac{np+q-1}{2}\le i\le\frac{np+q+1}{2}-g-p$, there is an isomorphism
% 	$$(\psp{-n}{1-n})^{-1}\circ\psm{-n}{1-n}:\sutg{-n}{i}\xra{\cong}\sutg{-n}{i+p}.$$

% \end{enumerate}
% \ecor

\blem\label{lem: comm diag for n+1 to n}
For a knot $K\subset Y$ and $n\in\intg$, there are commutative diagrams up to scalars
\begin{center}
\begin{minipage}{0.4\textwidth}
	\begin{equation*}
\xymatrix{
	\sut{\mu}\ar[rr]^{\psp{\mu}{n-1}}\ar[dd]_{\psm{\mu}{n-1}}&&\sut{n-1}\ar[dd]^{\psp{n-1}{\frac{2n-1}{2}}}\\
	&&\\
	\sut{n-1}\ar[rr]^{\psm{n-1}{\frac{2n-1}{2}}}&&\sut{\frac{2n-1}{2}}
	}
	\end{equation*}
\end{minipage}
\begin{minipage}{0.4\textwidth}
	\begin{equation*}
\xymatrix{
	\sut{\frac{2n-1}{2}}\ar[rr]^{\psp{\frac{2n-1}{2}}{n}}\ar[dd]_{\psm{\frac{2n-1}{2}}{n}}&&\sut{n}\ar[dd]^{\psm{n}{\mu}}\\
	&&\\
	\sut{n}\ar[rr]^{\psp{n}{\mu}}&&\sut{\mu}
	}
	\end{equation*}
\end{minipage}

\begin{minipage}{0.4\textwidth}
	\begin{equation*}
	\xymatrix{
\sut{\mu}\ar[rr]^{\psp{\mu}{n-1}}&&\sut{n-1}\\
&\sut{n}\ar[lu]^{\psm{n}{\mu}}\ar[ru]_{\psp{n}{n-1}}&
}
\end{equation*}
\end{minipage}
\begin{minipage}{0.4\textwidth}
	\begin{equation*}
\xymatrix{
\sut{\mu}\ar[rr]^{\psm{\mu}{n-1}}&&\sut{n-1}\\
&\sut{n}\ar[lu]^{\psp{n}{\mu}}\ar[ru]_{\psm{n}{n-1}}&
}
\end{equation*}
\end{minipage}
\end{center}

\begin{center}
\begin{minipage}{0.4\textwidth}
	\begin{equation*}
\xymatrix{
\sut{n-1}\ar[rr]^{\psp{\frac{2n-1}{2}}{n-1}}\ar[dr]_{\psp{n-1}{n}}&&\sut{\frac{2n-1}{2}}\ar[dl]^{\psm{\frac{2n-1}{2}}{n}}\\
&\sut{n}&
}	
\end{equation*}
\end{minipage}
\begin{minipage}{0.4\textwidth}
		\begin{equation*}
\xymatrix{
\sut{n-1}\ar[rr]^{\psm{\frac{2n-1}{2}}{n-1}}\ar[dr]_{\psm{n-1}{n}}&&\sut{\frac{2n-1}{2}}\ar[dl]^{\psp{\frac{2n-1}{2}}{n}}\\
&\sut{n}&
}		
\end{equation*}\end{minipage}
\end{center}
\elem

% Suppose $\alpha$ is a connected non-separating simple closed curve on $\partial (Y\backslash N(K))$, we can push $\alpha$ into the interior of $Y\backslash N(K)$ and apply the surgery exact triangle associated to the surgeries along $\alpha$ with respect to the framing induced by $\partial (Y\backslash N(K))$. According to \cite[Section 4]{baldwin2016contact}, when $\alpha$ intersects the suture at two points, the $0$-surgery along $\alpha$ corresponds to a $2$-handle attachment along $\alpha$ and hence leads to the Dehn filling of $Y\backslash N(K)$ along $\alpha$.
\begin{lem}\label{lem: surgery triangles}
	For any $n\in\intg$, we have the following exact triangles.
\begin{center}
\begin{minipage}{0.4\textwidth}
	\begin{equation*}
\xymatrix{
\sut{n}\ar[rr]^{H_n}&&\sut{n+1}\ar[dl]^{F_{n+1}}\\
&\dehny{}\ar[ul]^{G_n}&
}	
\end{equation*}
\end{minipage}
\begin{minipage}{0.4\textwidth}
		\begin{equation*}
\xymatrix{
\sut{\mu}\ar[rr]^{A_{n-1}}&&\sut{n-1}\ar[dl]^{B_{n-1}}\\
&\dehny{n}\ar[ul]^{C_n}&
}	
\end{equation*}\end{minipage}
\end{center}
\end{lem}

\blem\label{lem: comm diag for n,n+1,dehn}
For any $n\in\intg$, we have the following commutative diagrams up to scalars
\begin{center}
\begin{minipage}{0.4\textwidth}
	\begin{equation*}
\xymatrix{
\sut{n}\ar[rr]^{\psp{n}{n+1}}&&\sut{n+1}\\
&\dehny{}\ar[lu]^{G_n}\ar[ru]_{G_{n+1}}&
}	
\end{equation*}
\end{minipage}
\begin{minipage}{0.4\textwidth}
		\begin{equation*}
\xymatrix{
\sut{n}\ar[rr]^{\psm{n}{n+1}}&&\sut{n+1}\\
&\dehny{}\ar[lu]^{G_n}\ar[ru]_{G_{n+1}}&
}	
\end{equation*}
\end{minipage}
\end{center}

\begin{center}
\begin{minipage}{0.4\textwidth}
	\begin{equation*}
\xymatrix{
\sut{n}\ar[rr]^{\psp{n}{n+1}}\ar[dr]_{F_{n}}&&\sut{n+1}\ar[dl]^{F_{n+1}}\\
&\dehny{}&
}	
\end{equation*}
\end{minipage}
\begin{minipage}{0.4\textwidth}
		\begin{equation*}
\xymatrix{
\sut{n}\ar[rr]^{\psm{n}{n+1}}\ar[dr]_{F_{n}}&&\sut{n+1}\ar[dl]^{F_{n+1}}\\
&\dehny{}&
}	
\end{equation*}\end{minipage}
\end{center}
\elem

% \blem\label{lem: F_n and G_n are iso when n large}
% Let $F_n$ and $G_n$ be defined as in Lemma \ref{lem: surgery triangles}. Then for any large enough integer $n$, we have the following properties
% \begin{enumerate}
% 	\item The map $G_{n-1}$ is zero and $F_n$ is surjective. Moreover, for any grading $$g-\frac{np-q-1}{2}\le i_0\le\frac{np-q-1}{2}-g-p+1,$$the restriction of the map
% 	$$F_n:\bigoplus_{i=0}^{p-1}\sutg{n}{i_0+i}\to \dehny{}$$
% 	is an isomorphism.
% 	\item The map $F_{-n+1}$ is zero and $G_{-n}$ is injective. Moreover, for any grading $$g-\frac{np+q-1}{2}\le i_0\le \frac{np+q-1}{2}-g-p+1,$$the map
% 	$${\rm Proj}\circ G_{-n}:\dehny{}\to \bigoplus_{i=0}^{p-1}\sutg{-n}{i_0+i},$$
% 	is an isomorphism, where 	$${\rm Proj}:\sut{-n}\ra\bigoplus_{i=0}^{p-1}\sutg{-n}{i_0+i}$$
% 	is the projection.
% \end{enumerate}
% \elem

\blem[{\cite[Lemma 4.17, Proposition 4.26, Lemma 4.29 and Proposition 4.32]{LY2020}}]\label{lem: F_n and G_n are iso when n large}
Let $F_n$ and $G_n$ be as in Lemma \ref{lem: surgery triangles}. Then, for any large enough integer $n$, the following properties hold:
\begin{enumerate}
	\item The map $G_{n-1}$ is zero and $F_n$ is surjective. Moreover, for any grading $$g-\frac{np-q-1}{2}\le i_0\le\frac{np-q-1}{2}-g-p+1,$$the restriction of the map
	$$F_n:\bigoplus_{i=0}^{p-1}\sutg{n}{i_0+i}\to \dehny{}$$
	is an isomorphism.
	\item The map $F_{-n+1}$ is zero and $G_{-n}$ is injective. Moreover, for any grading $$g-\frac{np+q-1}{2}\le i_0\le \frac{np+q-1}{2}-g-p+1,$$the map
	$${\rm Proj}\circ G_{-n}:\dehny{}\to \bigoplus_{i=0}^{p-1}\sutg{-n}{i_0+i},$$
	is an isomorphism, where 	$${\rm Proj}:\sut{-n}\ra\bigoplus_{i=0}^{p-1}\sutg{-n}{i_0+i}$$
	is the projection.
\end{enumerate}
\elem

\blem\label{lem: -1 surgery and bypass}
For any $n\in\intg$, let the maps $H_{n}$ and $\psi^{n}_{\pm,n+1}$ be as in Lemma \ref{lem: surgery triangles} and Lemma \ref{lem: bypass n,n+1,mu} respectively. Then there exist $c_1,c_2\in\mathbb{C}\backslash\{0\}$ such that
$$H_{n}=c_1\psp{n}{n+1}+c_2\psm{n}{n+1}$$
\elem

\begin{conv}
Though maps between projectively transitive systems are only well-defined up to scalars in $\mathbb{C}\backslash\{0\}$, in \cite[Section 2.3]{LY2022integral1}, we introduced a way to fix the representatives of the systems and the scalars of maps between them. Hence, we can consider the sutured instanton homologies used in this paper as actual vector spaces, and all commutative diagrams above hold without introducing scalars, except for the second commutative diagram in Lemma \ref{lem: com diag for n,n+1,n+2}. Moreover, we can set the scalars $c_1=1 $ and $c_2=-1$ for any $n\in\mathbb{Z}$ in Lemma \ref{lem: -1 surgery and bypass}.
\end{conv}

\subsection{Integral surgery formulae}
Suppose $K\subset Y$ is a rationally null-homologous knot with a Seifert surface $S$. Suppose $(\lambda,\mu)$ is the chosen framing for $K$ and $(p,q)$ is defined as in Section \ref{sec: SHI, 1}. Then we state two versions of integral surgery formulae, one from the sutured theory and the other from the bent complex.

\bthm[{\cite[Theorem 3.1]{LY2022integral1}}]\label{thm: integral surgery formula} Suppose $m$ is a fixed integer such that $mp-q\neq 0$, \textit{i.e.}, the suture $\Ga_m$ is not parallel to $\partial S$. Then, for any large enough integer $k$, there exists an exact triangle
\begin{equation*}\label{eq: generalized mapping cone}
	\xymatrix{
	\mathbf{\Ga}_{\frac{2m+2k-1}{2}}\ar[rr]^{\pi}&&\mathbf{\Ga}_{m+2k-1}\ar[dl]\\
	&\mathbf{Y}_m\ar[ul]&
	}
\end{equation*}
	where $\pi=\pi^+_{m,k}+\pi^-_{m,k}$ and
	$$\pi_{m,k}^{\pm}=\Psi^{m+k}_{\pm,m+2k-1}\circ\psi^{\frac{2m+2k-1}{2}}_{\mp,m+k}$$are compositions of bypass maps. Hence we have an isomorphism$$\mathbf{Y}_m\cong H(\cone(\pi_{m,k}^++\pi_{m,k}^-))\cong \ker\pi\oplus \cok \pi.$$
\ethm

In \cite[Section 3.4]{LY2021large}, for any rationally null-homologous knot $K\subset Y$, we constructed  two spectral sequences $\{E_{r,+},d_{r,+}\}_{r\ge 1}$ and $\{E_{r,-},d_{r,-}\}_{r\ge 1}$ from $\sut{\mu}$ to $\dehny{}$, where the $\mathbb{Z}$-grading shift of $d_{r,\pm}$ is $\pm rp$. In short, we obtain two spectral sequences from the following unrolled exact couples about bypass maps
\begin{equation}\label{eq: useful triangles}
\xymatrix@R=6ex{
\cdots&\ar[l]\sut{n+1}\ar[dr]_{\psi_{\pm,\mu}^{n+1}}&&\sut{n}\ar[ll]_{\psi_{\pm,n+1}^{n}}\ar[dr]_{\psi_{\pm,\mu}^n}&&\sut{n-1}\ar[ll]_{\psi_{\pm,n}^{n-1}}\ar[dr]_{\psi_{\pm,\mu}^{n-1}}&&\sut{n-2}\ar[ll]_{\psi_{\pm,n-1}^{n-2}}&\ar[l]\cdots
\\
&\cdots&\sut{\mu}\ar[ur]^{\psi_{\pm,n}^\mu}&&\sut{\mu}\ar[ur]^{\psi_{\pm,n-1}^\mu}&&\sut{\mu}\ar[ur]^{\psi_{\pm,n-2}^\mu}&\cdots&
}
\end{equation}
The spectral sequences are independent of the choice of $n$. Then we lift the spectral sequences to filtered chain complexes with differentials $d_+$ and $d_-$ by fixing an inner product on $\sut{\mu}$. By construction we have $$H(\sut{\mu},d_+)\cong H(\sut{\mu},d_-)\cong \dehny{}.$$
\bdefn[{\cite[Construction 3.27 and Definition 5.12]{LY2021large}}]\label{defn: complex A}
For any rationally null-homologous knot $K\subset Y$, let $B^\pm(K)$ be the complexes $(\sut{\mu},d_\pm)$. For any integer $s$, define the \textbf{bent complex}\[A(K,s)\deq (\bigoplus_{k\in\mathbb{Z}}(\sut{\mu},s+kp),d_s),\]where for any element $x\in (\sut{\mu},s+kp)$, 
\[
d_s(x)=\begin{cases}
d_+(x)&k>0,\\
d_+(x)+d_-(x)&k=0,\\
d_-(x)&k<0.
\end{cases}\]
Let $B^\pm(K,s)$ be copies of $B^\pm(K)$. Define\[\pi^+(K,{s}):A(K,s)\to B^+(K,s)\aand \pi^-(K,{s}):A(K,s)\to B^-(K,s)\]by
\[\pi^+(K,s)(x)=\begin{cases}
x&k\ge 0,\\
0&k< 0,\end{cases}\aand \pi^-(K,s)(x)=\begin{cases}
x&k\le 0,\\
0&k>0,\end{cases}\]where $x\in (\sut{\mu},s+kp)$. Define $$\pi^\pm(K):\bigoplus_{s\in\mathbb{Z}}A(K,s)\to \bigoplus_{s\in\mathbb{Z}}B^\pm(K,s)$$by putting $\pi^\pm(K,s)$ together for all $s$ . We also use the same notation for the induced map on homology. If $K$ is fixed, we will omit it in $A(K,s),B^\pm(K,s)$ and $\pi^\pm(K,s)$.
\edefn
\bthm[{\cite[Theorem 3.18]{LY2022integral1}}]\label{thm: integral bent complex}
Suppose $m$ is a fixed integer such that $mp-q\neq 0$. Then there exists an isomorphism $$\Xi_{m}:\bigoplus_{s\in\mathbb{Z}}H(B^+(s))\xra{\cong}\bigoplus_{s\in\mathbb{Z}}H(B^-(s+mp-q))$$as the direct sum of isomorphisms$$\Xi_{m,s}:H(B^+(s))\xra{\cong}H(B^-(s+mp-q))$$ so that $$\dehny{m}\cong H\bigg(\cone(\pi^-+\Xi_m\circ \pi^+:\bigoplus_{s\in\mathbb{Z}}H(A(s))\to  \bigoplus_{s\in\mathbb{Z}}H(B^-(s)))\bigg).$$
\ethm
\brem
Theorem \ref{thm: integral surgery formula} is a little stronger than Theorem \ref{thm: integral bent complex} when we consider the $\Lambda^* H_1(Y;\mathbb{C})$-action on the sutured instanton homology. From Corollary \ref{cor: bypass maps for Borromean knot}, the action is trivial on $\sut{\mu}$ of the Borromean knot, and hence is trivial on the bent complex. But it is nontrivial on $\sut{n}$ by Lemma \ref{lem: gamma_n for Borromean knot}. In this paper, we will use both versions of surgery formulae.
\erem

\subsection{Truncation of the integral surgery formulae}

In this subsection, we will use the following algebraic lemma to truncate the integral surgery formula.

\blem\label{lem: reduce the chain complex}
Suppose $(C,d_C)$ is a chain complex and suppose $C=D\oplus E\oplus F$. For $A,B\in\{D,E,F\}$, we write $d^A_B:A\to B$ for the restriction of $d_C$. We write elements in $C$ as column vectors. Suppose $d_C$ has the form \begin{equation*}
d_C = 
\begin{pmatrix}
0 & d^D_E & 0 \\
0 & 0 & 0 \\
0  & d^F_E  & d^F_F
\end{pmatrix}
\end{equation*}where $d^D_E$ is an isomorphism. Then we have an isomorphism
$$H(C,d_C)\cong H(F,d_{F}^F).$$
\elem

\bpf
We have a short exact sequence$$0\to D\oplus E\to C\to F\to 0$$which induces a long exact sequence. The assumption on $d_C$ implies $H(D\oplus E)=0$ and hence $H(C)\cong H(F)$. 
\epf
We also need some structural lemmas for sutured instanton homologies.

\blem\label{lem: structure of Gamma_n}
Suppose $K\subset Y$ is a framed rationally null-homologous knot. Suppose $n\in\intg$ such that $(n-1)p-q\geq 2g$, then we have the following.
\begin{enumerate}
	\item When $|i|>\frac{np-q-1}{2}+g$, we have $\sutg{n}{i}=0$.
	\item When $\frac{np-q-1}{2}+g\geq i\ge  g-\frac{np-q-1}{2}$, we have an isomorphism
	$$\psi^{n}_{\mp,n+1}:\sutg{n}{\pm i}\xra{\cong }\sutg{n+1}{\pm i\pm \frac{p}{2}}.$$
	\item When $\frac{np-q-1}{2}-g\geq i,j\geq g-\frac{np-q-1}{2}$ and $i-j=p$, we have an isomorphism
	$$(\psm{n}{n+1})^{-1}\circ\psp{n}{n+1}:\sutg{n}{i}\xra{\cong}\sutg{n}{j}.$$
	\item When $g-\frac{np-q-1}{2}\le i_0\le \frac{np-q-1}{2}-g-p+1$, the restriction of the map
	$$F_n:\bigoplus_{i=0}^{p-1}\sutg{n}{i_0+i}\to \dehny{}$$
	is an isomorphism.
\end{enumerate}
\elem
\bpf
Part (1) follows directly from Lemma \ref{lem: vanishing grading}. Part (2) and (3) follow from Lemma \ref{lem: vanishing grading} and Lemma \ref{lem: bypass n,n+1,mu}. Part (4) follows from \cite[Lemma 2.19 part (1)]{LY2022integral1}.
\epf

\blem\label{lem: structure of Gamma_2n-1/2}
Suppose $K\subset Y$ is a framed rationally null-homologous knot. Suppose $n\in\intg$ such that $(n-1)p-q\geq 2g$. Then we have the following.
\begin{enumerate}
	\item When $|i|>\frac{(2n-1)p-2q-1}{2}+g$, we have
	$$\sutg{\frac{2n-1}{2}}{i}=0.$$
	\item When $\frac{(2n-1)p-2q-1}{2}+g\geq i\ge g-\frac{p-1}{2}$, we have an isomorphism
	$$\psi^{\frac{2n-1}{2}}_{\pm,n}:\sutg{\frac{2n-1}{2}}{\pm i}\xra{\cong}\sutg{n}{\pm i\mp\frac{(n-1)p+q}{2}}.$$
	\item When $\frac{(2n-1)p-2q-1}{2}-g\geq i\ge g-\frac{(2n-1)p-2q-1}{2}$, we have an isomorphism
	$$\sutg{\frac{2n-1}{2}}{i}\cong H(A(i)),$$
	where $A(i)$ is the bent complex defined as in Definition \ref{defn: complex A}.
\end{enumerate}
\elem

\bpf
Part (1) follows directly from Lemma \ref{lem: vanishing grading}. Part (2) follows from Lemma \ref{lem: bypass n+1,n,2n+1/2} and Lemma \ref{lem: structure of Gamma_n} part (1). Part (3) follows from \cite[Theorem 3.23]{LY2021large}.
\epf

\blem\label{lem: the map pi^pm_m,k}
Suppose $K\subset Y$ is a framed rationally null-homologous knot. Suppose $\pi^{\pm}_{m,k}$ is defined as in Theorem \ref{thm: integral surgery formula}. Let $\pi^{\pm,i}_{m,k}$ be the restriction of $\pi^{\pm}_{m,k}$ on $\sutg{\frac{2m+2k-1}{2}}{i}$. Then we have the following. 
\begin{enumerate}
	\item We have
	$$\im\pi^{\pm,i}_{m,k}\subset\sutg{m+2k-1}{i\pm\frac{mp-q}{2}}.$$
	\item When $i>\frac{p-1}{2}+g$, we have $\pi^{\pm,\pm i}_{m,k}=0$.
	\item When $i\ge \frac{p-1}{2}+g$, the map $\pi_{m,k}^{\mp,\pm i}$ is an isomorphism. 
\end{enumerate}
\elem
\bpf
Part (1) follows directly from the grading shifts in Lemma \ref{lem: bypass n,n+1,mu} and Lemma \ref{lem: bypass n+1,n,2n+1/2}. For the grading $i$ in part (2), by Lemma \ref{lem: bypass n+1,n,2n+1/2} and Lemma \ref{lem: structure of Gamma_n}, we have  $$\psi^{\frac{2m+2k-1}{2}}_{\mp,m+k}=0$$and hence $\pi^{\pm,\pm i}_{m,k}=0$. Part (3) follows from Lemma \ref{lem: structure of Gamma_2n-1/2} part (2) and Lemma \ref{lem: structure of Gamma_n} part (2).
\epf
\bprop\label{prop: truncated integral surgery formula}
Suppose $m\in\intg$ such that $mp-q\neq 0$. Then, for any large enough integer $k$, we have the following.

\benu

\item If $(m-1)p-q+2<2g$, then
% 	$$\dehny{m}\cong H\bigg({\rm Cone}(\pi':\bigoplus_{i=-\frac{(2m+1)p-2q-1}{2}-g}^{\frac{(2m+1)p-2q-1}{2}+g}\sutg{\frac{2m+2k-1}{2}}{i}\to\bigoplus_{i=-\frac{(m+1)p-q-1}{2}-g}^{\frac{(m+1)p-q-1}{2}+g}\sutg{m+2k-1}{i})\bigg),$$
	
	$$\begin{aligned}\dehny{m}&\cong H\bigg({\rm Cone}(\pi^\p:\bigoplus_{|i|<\frac{p-1}{2}+g}\sutg{\frac{2m+2k-1}{2}}{i}\to\bigoplus_{|i|<\frac{(1-m)p+q-1}{2}+g}\sutg{m+2k-1}{i})\bigg)\\&\cong H\bigg({\rm Cone}(\pi^\pp:\bigoplus_{|s|<\frac{p-1}{2}+g}H(A(s))\to\bigoplus_{s=-\frac{p-1}{2}+1-mp+q-g}^{\frac{p-1}{2}+g-1}H(B^-(s)))\bigg).\end{aligned}$$
	
\item If $(m-1)p-q+2\ge 2g$, then $$\dehny{m}\cong\bigoplus_{i=\frac{p-1}{2}-mp+q+g}^{\frac{p-1}{2}+g-1}\sutg{\frac{2m+2k-1}{2}}{i}\cong \bigoplus_{s=\frac{p-1}{2}-mp+q+g}^{\frac{p-1}{2}+g-1} H(A(s)).$$
\eenu
	Here $\pi^\p$ and $\pi^\pp$ are the restrictions of $\pi$ and $\pi^-+\Xi_m\circ \pi^+$ is defined as in Theorem \ref{thm: integral surgery formula} and Theorem \ref{thm: integral bent complex}.
\eprop
\bpf
This is a reduction of Theorem \ref{thm: integral surgery formula} and Theorem \ref{thm: integral bent complex}. We only prove the first isomorphism in each case. The proof of the second isomorphism follows directly from the reformulation of the integral surgery formula by bent complexes in \cite[Section 3.3]{LY2022integral1}.

From Lemma \ref{lem: the map pi^pm_m,k}, the grading shift of $\pi^\pm_{m,k}$ is $\frac{mp-q}{2}$. When $i> \frac{p-1}{2}+g$, we have $\pi^{+, i}_{m,k}=0$ and $\pi^{-, i}_{m,k}$ is an isomorphism. Let $C_1$ be the total mapping cone in Theorem \ref{thm: integral surgery formula} and let $$D_1=\bigoplus_{i>\frac{p-1}{2}+g} \sutg{\frac{2m+2k-1}{2}}{i}\aand E_1=\bigoplus_{i> \frac{p-1}{2}-\frac{mp-q}{2}+g}\sutg{m+2k-1}{i}.$$Then $\pi$ restricts to $\pi^-_{m,k}$ on $D$ and induces an isomorphism of $D\cong E$. Then we apply Lemma \ref{lem: reduce the chain complex} to remove $D_1\oplus E_1$ from $C_1$. Let $C_2$ be the quotient $C_1/(D_1\oplus E_1)$. Since $\pi^{-,i}_{m,k}$ is also an isomorphism for $i=\frac{p-1}{2}+g$, we can apply Lemma \ref{lem: reduce the chain complex} again to remove $$D_2=\sutg{\frac{2m+2k-1}{2}}{\frac{p-1}{2}+g}\aand E_2=\sutg{\frac{2m+2k-1}{2}}{\frac{p-1}{2}-\frac{mp-q}{2}+g}$$from $C_2$. Let $C_3$ be the quotient $C_2/(D_2\oplus E_2)$. Note that the grading induced by the Seifert surface is either a $\mathbb{Z}$-grading or a $(\mathbb{Z}+\frac{1}{2})$-grading. If $$\frac{p-1}{2}-\frac{mp-q}{2}+g>\frac{1}{2},$$then we can similarly apply Lemma \ref{lem: reduce the chain complex} to $$D_3=\bigoplus_{i<-\frac{p-1}{2}-g} \sutg{\frac{2m+2k-1}{2}}{i}\aand E_3=\bigoplus_{i< -\frac{p-1}{2}+\frac{mp-q}{2}-g}\sutg{m+2k-1}{i}$$and then also $$D_4=\sutg{\frac{2m+2k-1}{2}}{-\frac{p-1}{2}-g}\aand E_4=\sutg{\frac{2m+2k-1}{2}}{-\frac{p-1}{2}+\frac{mp-q}{2}-g}.$$We conclude that $$H(C_1)\cong H\bigg({\rm Cone}(\pi':\bigoplus_{|i|<\frac{p-1}{2}+g}\sutg{\frac{2m+2k-1}{2}}{i}\to\bigoplus_{|i|<\frac{(1-m)p+q-1}{2}+g}\sutg{m+2k-1}{i})\bigg).$$If $$\frac{p-1}{2}-\frac{mp-q}{2}+g\le \frac{1}{2},$$then we can apply Lemma \ref{lem: reduce the chain complex} to $$D_3^\p=\bigoplus_{i<\frac{p-1}{2}-mp+q+g} \sutg{\frac{2m+2k-1}{2}}{i}\aand E_3^\p=\bigoplus_{i< \frac{p-1}{2}-\frac{mp-q}{2}+g}\sutg{m+2k-1}{i}.$$We conclude that $$H(C_1)\cong \bigoplus_{i=\frac{p-1}{2}-mp+q+g}^{\frac{p-1}{2}+g-1}\sutg{\frac{2m+2k-1}{2}}{i}$$
\epf

If $K$ is null-homologous, then $(p,q)=(1,0)$. The inequality $(m-1)p-q+2\ge 2g$ reduces to $m\ge 2g-1$. In such a case, the result in Proposition \ref{prop: truncated integral surgery formula} is indeed stronger than the large surgery formula in \cite[Theorem 1.22]{LY2021large} because the assumption in that paper is $m\ge 2g+1$. This difference is essential when $g$ is small (\textit{e.g.} $g=1$).

\begin{prop}\label{prop: genus one knot}
	Suppose $K\subset Y$ is a null-homologous knot bounding a Seifert surface of genus $1$. Then for any $m\in\mathbb{N}_+$, we have
	$$\dim I^{\sharp}(Y_{m}(K))=\dim I^{\sharp}(Y_{1}(K))+(m-1)\cdot \dim I^{\sharp}(Y).$$
\end{prop}
\bpf
Since $g=1$, we apply Proposition \ref{prop: truncated integral surgery formula} part (2) to any $m>0$. In particular, we have $$\dehny{1}\cong H(A(0))\aand \dehny{m}\cong \bigoplus_{s=-m+1}^{0}H(A(s)).$$By the construction of $A(s)$ in Definition \ref{defn: complex A}, we know $$H(A(s))\cong H(B^+(s))\cong \dehny{}$$for any $s<0$. Since $\dim I^\sharp(-Y)=\dim I^\sharp(Y)$ for any closed 3-manifold $Y$, we conclude the dimension equality.
\epf
\brem
When $Y$ is a rational homology sphere, this corollary follows directly from the adjunction inequality for the instanton cobordism map; see for example \cite[Theorem 1.16]{baldwin2019lspace}. However, for technical reasons such an adjunction inequality relies on the assumption that the first Betti number of the cobordism vanishes. So when $b_1(Y)>0$, the existing adjunction inequality does not apply.
\erem

\subsection{Instanton tau invariant}\label{sec: SHI, 2}
We present some results from \cite{li2019tau} for knots inside $S^3$.
\bdefn[{\cite[Definition 5.4]{li2019direct}}]\label{defn: khi minus}
Suppose $K\subset Y$ is a rationally null-homologous knot. Let $\khi(-Y,K)$ be the direct limit of $$\cdots\to \sut{n}\xra{\psi_{-,n+1}^n}\sut{n+1}\xra{\psi_{-,n+2}^{n}}\sut{n+2}\to \cdots.$$Let $U$ be the action on $\khi(-Y,K)$ defined by $\{\psi_{+,n+1}^n\}_{n\in\mathbb{N}_+}$. It is well-defined due to the commutativity in Lemma \ref{lem: com diag for n,n+1,n+2}.
\edefn
\bdefn[{\cite[Definition 5.7]{li2019direct}}]\label{defn: tau}
Suppose $K\subset S^3$ is a knot. We define
$$\tau_I(K)=\max\{i~|~\exists~x\in \khi(-S^3,K,i)~s.t.~U^j\cdot x\neq0~{\rm for~any~}j\geq0\}.$$
\edefn

We have the following basic properties for $\tau_I$.
\blem\label{lem: basic properties of tau}
Suppose $K\subset S^3$ is a knot. Then we have the following.
\begin{enumerate}
	\item (\cite[Proposition 3.17 and Corollary 5.3]{li2019tau}) For $n\in\intg$ large enough, we have
	\beq
	\tau_I(K)&=\max\{i~|~\exists~x\in\sutg{n}{i}~s.t.~F_{n}(x)\neq0\in I^{\sharp}(-S^3)\}-\frac{n-1}{2}\\
	&=\min\{i~|~\exists~x\in\sutg{n}{i}~s.t.~F_{n}(x)\neq0\in I^{\sharp}(-S^3)\}+\frac{n-1}{2}.
	\eeq
	\item (\cite[Proposition 1.12 and Proposition 1.14]{li2019tau}) We have $\tau_I(K)=-\tau_I(\widebar{K})$ where $\widebar{K}$ is the mirror of $K$.
\end{enumerate}
\elem

\blem[{\cite[Section 5]{li2019tau}}]\label{lem: symmetry and dimension of Gamma_n}
Suppose $K\subset S^3$ is a knot. Then we have the following.
\begin{enumerate}
	\item For any $*\in\mathbb{Q}\cup\{\mu\}$, we have $\sutg{*}{i}\cong\sutg{*}{-i}$.
	\item For any $n\in\mathbb{Z}$, we have
	$$\dim\sut{n}=\dim\sut{-2\tau_I(K)}+|n+2\tau_I(K)|$$
\end{enumerate}
\elem

\section{A rational surgery formula}\label{A rational surgery formula}

Suppose $K$ is a null-homologous knot in a 3-manifold $Y$. In this section, we will study the $u/v$-surgery on $K$. The integral surgery formula Theorem \ref{thm: integral bent complex} is an analog of the Morse (integral) surgery formula for Heegaard Floer homology in \cite[Section 6]{Ozsvath2011rational}. To generalize the formula to rational surgeries, we use the same strategy as in \cite[Section 7]{Ozsvath2011rational}. For simplicity, we use similar notations as in Ozsv\'ath-Szab\'o's work \cite{Ozsvath2011rational}. The symbols $(p,q,r,a)$ in \cite{Ozsvath2011rational} are replaced by $(u,v,r,m)$ since we define $p$ and $q$ in Section \ref{sec: SHI, 1} (indeed $(p,q)=(1,0)$ because $K$ is null-homologous). Suppose $$m=\lfloor \frac{u}{v}\rfloor$$ is the greatest integer smaller than or equal to $u/v$, and $$\frac{u}{v}=m+\frac{r}{v}.$$Let $O_{v/r}$ be the knot obtained by viewing one component of the Hopf link as a knot inside the lens space $L(v,-r)$ thought of as $-v/r$ surgery on the
other component of the Hopf link, which is framed by the Seifert framing of the unknot in $S^3$. Note that $O_{v/r}$ is a \textbf{core knot} of the lens space, \textit{i.e.}, the complement is a solid torus. Since $Y_{u/v}(K)$ can be obtained by $m$-surgery on the connected sum $$K\# O_{v/r}\subset Y\# L(v,-r),$$
it suffices to understand the bent complex associated to  $K\# O_{v/r}$ in terms of the bent complex of $K$. 

\subsection{The connected sum with a core knot}\label{subsec: The connected sum with a core knot}
Given two knots $K_i\subset Y_i$ for $i=1,2$, let $K^\p\subset Y^\p$ be the connected sum of $K_1$ and $K_2$. Note that $Y^\p\backslash N(K^\p)$ is obtained from gluing $Y_i\backslash N(K_i)$ by an annulus along the meridians of $K_i$ for $i=1,2$. Conversely, the disjoint union of $Y_i\backslash N(K_i)$ is obtained from $Y^\p\backslash N(K^\p)$ by a product annulus decomposition in the sense of \cite[Proposition 6.7]{kronheimer2010knots}. The instanton version of that proposition implies \begin{equation}\label{eq: connected sum}
    \khii(Y^\p,K^\p)\cong \khii(Y_1,K_1)\otimes \khii(Y_2,K_2).
\end{equation}Moreover, if $K_i$ are rationally null-homologous, in our previous work \cite[Proposition 5.15]{LY2021large}, we generalized the above isomorphism to a graded version with respect to the gradings associated to Seifert surfaces.  (Note the result in \cite[Proposition 5.15]{LY2021large} is stated for knots inside rational homology spheres but the proof works for rationally null-homologous knots inside arbitrary $3$-manifolds.) However, we need a stronger version of the connected sum formula that encodes the information in bent complexes. Inspired by the formula in Heegaard Floer theory \cite[Lemma 7.1]{ozsvath2004holomorphicknot}, we have the following conjecture. 

\begin{conj}\label{conj: connected sum}

Suppose $K_i\subset Y_i$ for $i=1,2$ are rationally null-homologous knots. Then there exist chain homotopy equivalences$$B^\pm(K_1\# K_2)\simeq B^\pm(K_1)\otimes B^\pm(K_2),$$where $B^\pm(K)$ is defined in Definition \ref{defn: complex A}.
\end{conj}
% The proof of the above conjecture is an ongoing project with Ghosh \cite{GLY2022}. 
In this subsection, we only prove the special case where $K_2$ is a core knot in a lens space. First, we present some results for core knots.

\blem[{\cite[Proposition 4.10]{li2019direct}}]\label{lem: dim}
Suppose $K$ is a core knot in a lens space $Y$. Then we have
\begin{equation*}\label{eq: dim p}
    \sutg{r/s}{i}\cong \mathbb{C}\text{ for any } |i|\le \frac{|rp-sq|-1}{2}.
\end{equation*}For other grading $i$, we have $\sutg{r/s}{i}=0$.
\elem

\bcor\label{cor: split}
Suppose $K$ is a core knot in a lens space $Y$. Then the bypass exact triangles in Lemma \ref{lem: bypass n,n+1,mu} are always split, and there are two canonical isomorphisms induced by bypass maps between the direct sum of two spaces with smaller dimensions and the third space.
\ecor
\bpf
From Lemma \ref{lem: dim}, we know the dimensions of $\sut{n},\sut{n+1}, \sut{\mu}$ are $|np-q|,|(n+1)p-q|,|p|$, respectively. Since the sum of two smaller integers equals the third integer, we know the triangle always splits. Since each nontrivial grading summand $\sut{n},\sut{n+1}, \sut{\mu}$ are $1$-dimensional, the restrictions of the bypass maps induce the canonical isomorphisms.
\epf

From Lemma \ref{lem: dim}, for a core knot $K\subset Y$, we have $$\dim \khii(-Y,K)=\dim I^\sharp (-Y)\aand d_{\pm}=0.$$
Then Conjecture \ref{conj: connected sum} reduces to the following proposition.
\bprop[]\label{prop: Kunneth diff}
Suppose $K_1\subset Y_1$ is a rationally null-homologous knot and $K_2\subset Y_2$ is a core knot in a lens space. Then there exist identifications$$B^\pm(K_1\# K_2)= B^\pm(K_1)\ot \khii(-Y_2,K_2).$$
\eprop

\begin{conv}
To distinguish sutures for different knot complements, we write $\sut{r/s}^{\bullet}$, $\psi_{\pm,n+1}^{n,\bullet}$, and $F_{n}^\bullet$ with $\bullet\in\{1,2,\#\}$ for the sutured instanton homology, bypass maps, and the cobordism maps in Lemma \ref{lem: surgery triangles} associated to the knots $K_1,K_2$ and $K_1\# K_2$.
\end{conv}

To prove Proposition \ref{prop: Kunneth diff}, we need the following lemma, which generalizes results in \cite[Section 3.2]{li2019tau}. 
\blem[]\label{lem: contact map}
% Suppose $K_1\subset Y_1$ is a rationally null-homologous knot and $K_2\subset Y_2$ is a core knot in a lens space. Suppose $n$ and $k$ are large integers. Then there exist maps$$C_{\pm,n+k}^{n,k}:(\sut{n}^1,?)\otimes (\sut{k}^2,?)\to (\sut{n+k}^{\#},?)$$$$C_{\pm,\mu}^{\mu,k}:(\sut{\mu}^1,?)\otimes (\sut{k}^2,?)\to (\sut{\mu}^{\#},?)$$so that we have the following commutative diagrams.

Suppose $K_i\subset Y_i$ for $i=1,2$ is a rationally null-homologous knot. Suppose $n$ and $k$ are large integers. Then there exist maps$$C_{\pm,n+k}^{n,k}:\sut{n}^1\otimes \sut{k}^2\to \sut{n+k}^{\#}\aand C_{\pm,\mu}^{\mu,k}:\sut{\mu}^1\otimes \sut{k}^2\to \sut{\mu}^{\#}$$such that we have the following commutative diagrams.

\begin{equation}\label{eq: commutative diagram a}
\xymatrix@R=3ex{
\sut{\mu}^1\otimes \sut{k}^2\ar[rr]^{\psi_{\pm,n}^{\mu,1}\ot \id}\ar[dd]_{C_{\pm,\mu}^{\mu,k}}&&\sut{n}^1\otimes \sut{k}^2\ar[rr]^{\psi_{\pm,\mu}^{n,1}\ot\id}\ar[dd]_{C_{\pm,n+k}^{n,k}}&&\sut{\mu}^1\otimes \sut{k}^2\ar[dd]^{C_{\pm,\mu}^{\mu,k}}\\
&&\\
\sut{\mu}^{\#}\ar[rr]^{\psi_{\pm,n+k}^{\mu,\#}}&&\sut{n+k}^{\#}\ar[rr]^{\psi_{\pm,\mu}^{n+k,\#}}&&\sut{\mu}^{\#}
}
\end{equation}
\begin{equation}\label{eq: commutative diagram b}
\xymatrix@R=3ex{
\sut{n}^1\otimes \sut{k}^2\ar[rr]^{\psi_{\pm,n+1}^{n,1}\ot\id}\ar[dd]_{C_{\pm,m+n}^{n,k}}&&\sut{n+1}^1\otimes \sut{k}^2\ar[dd]^{C_{\pm,n+k+1}^{n+1,k}}\ar[rr]^{F_{n+1}^1\ot F_{k}^2}&&I^\sharp(-Y_1)\ot I^\sharp(-Y_2)\ar[dd]^{=}\\
&&\\
\sut{n+k}^{\#}\ar[rr]^{\psi_{\pm,n+k+1}^{n+k,\#}}&&\sut{n+k+1}^{\#}\ar[rr]^{F_{n+k+1}^\#}&&I^\sharp(-Y_\#)
}
\end{equation}

\begin{equation}\label{eq: commutative diagram c}
\xymatrix@R=6ex{
\sut{\mu}^1\otimes \sut{k}^2\ar[rr]^{\id\ot\psi_{\pm,\mu}^{k,2}}\ar[dr]^{C^{\mu,k}_{\pm,\mu}}&&\sut{\mu}^1\otimes \sut{\mu}^2\\
&\sut{\mu}^{\#}\ar[ur]^{=}&
}
\end{equation}
where the identification in (\ref{eq: commutative diagram b}) comes from the connected sum formula for $I^\sharp$ (\textit{cf.} \cite[Section 1.8]{li2018contact}) and the identification in (\ref{eq: commutative diagram c}) comes from the sutured decomposition along the product annulus.
\elem
\bpf
The proof is similar to the arguments in \cite[Section 4]{li2019tau}, especially the proof of \cite[Lemma 4.3]{li2019tau} and the proof of \cite[Proposition 1.14]{li2019tau}. Although the proofs in \cite{li2019tau} were only carried out for knots inside $S^3$, the same argument essentially works for rationally null-homologous knots in general $3$-manifolds. Here we only sketch the proofs as follows. We only prove the case involving positive bypasses. The case for the negative bypasses is similar. 

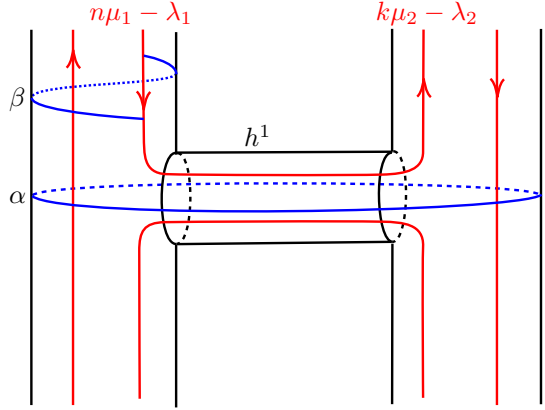
\begin{figure}[ht]
	\begin{overpic}[width=3in]{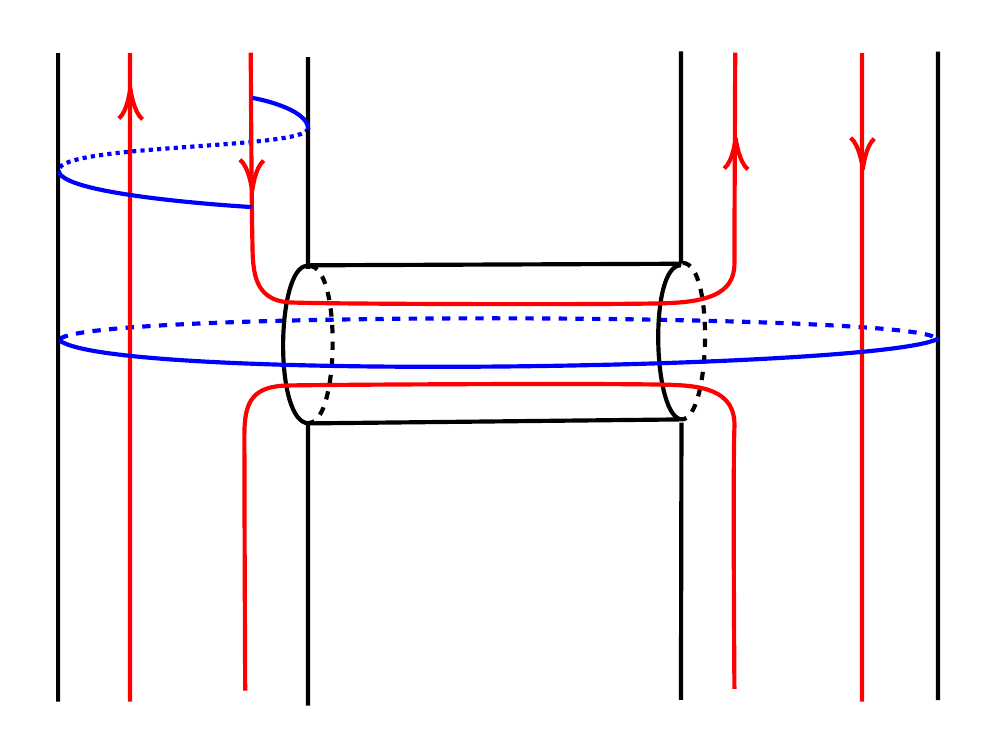}
		\put(16,70){\color{red}$n\mu_1-\lambda_1$}
		\put(66,70){\color{red}$k\mu_2-\lambda_2$}
		\put(2,55){$\beta$}
		\put(2,38){$\al$}
		\put(43,48){$h^1$}
	\end{overpic}
\caption{The $1$-handle $h^1$, the attaching curve $\al$ for the $2$-handle $h^2$, and the bypass arc $\beta$.}\label{fig: connected sum 1}
\end{figure}

We attach a $1$-handle $h^1$ to $(Y_1\backslash N(K_1),\Gamma^1_n)\sqcup (Y_2\backslash N(K_2),\Gamma_k^2)$ so that the two attaching points of the $1$-handles are on the curve $(n\mu_1-\lambda_1)\subset\Ga^1_n$ and $(k\mu_2-\lambda_2)\subset\Ga^2_k$ respectively. (For negative bypasses, we attach the $1$-handle to $(\lambda_1-n\mu_1)$ and $(\lambda_2-k\mu_2)$ accordingly. Note that the orientations of curves are different.) See Figure \ref{fig: connected sum 1}. Then we can attach a $2$-handle $h^2$ along the curve $\al$ which goes through the $1$-handle $h^1$ and intersects the suture obtained from attaching the $1$-handle twice, as shown in Figure \ref{fig: connected sum 1}. Define
$$C^{n,k}_{+,m+k}=C_{h^2}\circ C_{h^1}.$$
Here $C_{h^2}$ and $C_{h^1}$ are the contact handle attaching maps as introduced by Baldwin-Sivek in \cite{baldwin2016contact}. Pick the bypass arc $\beta$ such that it intersects the curve $(\lambda_2-k\mu_2)\subset\Ga^1_n$ once and its two endpoints are on the curve $(n\mu_1-\lambda_1)\subset\Ga^1_n$. See Figure \ref{fig: connected sum 1}. We know this is a positive bypass ({\it cf.} \cite{li2019direct}). Attaching a bypass along $\beta$ yields $(Y_1\backslash N(K_1),\Ga_{\mu})$. 

Let $h^{1\prime}$ and $h^{2\prime}$ be the corresponding $1$-handle and $2$-handle attached to $(Y_1\backslash N(K_1),\Ga_{\mu})\sqcup (Y_2\backslash N(K_2),\Gamma_k^2)$. Define
$$C^{\mu,k}_{+,\mu}=C_{h^{2\prime}}\circ C_{h^{1\prime}}.$$
The commutative diagrams (\ref{eq: commutative diagram a}) and (\ref{eq: commutative diagram b}) are straightforward since the bypass arc and the contact handles are disjoint from each other. See also \cite[Diagram (4.4)]{li2019tau} and proof of \cite[Proposition 1.14]{li2019tau} for more detailed discussions. The proof of (\ref{eq: commutative diagram c}) is similar to that of \cite[Equation (4.6)]{li2019tau}: let $\al'$ be the attaching curve of the $2$-handle $h^{2\prime}$. We can isotope $\al'$ into a suitable position so that this contact $2$-handle attachment corresponds to the one in the construction of the bypass map as in \cite{baldwin2016contact}.

\epf

% \blem\label{lem: compute grading of C}
% Suppose  $K_1,K_2$ and $C_{\pm,n+k}^{n,k}, C_{\pm,\mu}^{\mu,k}$ are knots and maps from Lemma \ref{lem: contact map}. Suppose further that the order of $K_1$ is $1$ (\textit{i.e.} $K_1$ is null-homologous) and the order of $K_2$ is $v$.  for $i=1,2$ and $\gcd(p_1,p_2)=1$. Then the order of $K_\#\deq K_1\# K_2$ is $p_1p_2$ and the grading shifts of $C_{\pm,n+k}^{n,k}$ and $C_{\pm,\mu}^{\mu,k}$ are as follows.

% \elem
% \bpf
% The knot $K_1$ is of order $1$ and $O_{v/r}$ is of order $v$. Then $$H_1(Y_1\backslash N(K_1))\cong H_1(Y_1)\oplus \mathbb{Z}\langle \mu_1\rangle\aand H_1(Y_2\backslash N(K_2))\cong\mathbb{Z}.$$We write $g_1=\mu_1$ and $g_2$ as the generator $H_1(Y_2\backslash N(K_2))$. Then $\mu_2=v\cdot g_2$ and $\lambda_2=r\cdot g_2$. A calculation on the homology shows \begin{equation}\label{eq: homology cal}
%     \begin{aligned}
%     H_1(Y_\#\backslash N(K_\#))&\cong \bigg(H_1(Y_1\backslash N(K_1))\oplus H_1(Y_2\backslash N(K_2))\bigg)/(\mu_1,\mu_2)\\&\cong \bigg(H_1(Y_1)\oplus\mathbb{Z}\langle g_1,g_2\rangle \bigg)/(g_1=v\cdot g_2)\\&\cong H_1(Y_1)\oplus\mathbb{Z}\langle g_\#\rangle\\&\cong H_1(Y_1\backslash N(K_1)),
% \end{aligned}
% \end{equation}where we write $g_\#$ as the generator. We also write ${\rm pr}$ as the projection to the summand generated by $g_\#$. Then $${\rm pr}(\mu_\#)=v\cdot g_\#\aand {\rm pr}(\lambda_\#)=r\cdot g_\#.$$Thus, the knot $K_\#$ is of order $v$ and $\partial S_\#\cdot \lambda_\#=-r$.
% \epf

\begin{proof}[Proof of Proposition \ref{prop: Kunneth diff}]
We only show the proof for $d_+$. The proof for $d_-$ is similar. To construct the differential $d_+$, we need to use triangles about positive bypass maps in (\ref{eq: useful triangles}). Suppose $m$ and $n$ are large integers. Consider the following diagram.

\begin{equation}\label{eq: useful triangles 2}
\xymatrix@R=6ex{
\cdots&\ar[l]\sut{n+1}^1\otimes \sut{k}^2\ar[dd]_{C_{+,n+k+1}^{n+1,k}}\ar[dr]_{\psi_{+,\mu}^{n+1,1}\ot\id}&&\sut{n}^1\otimes \sut{k}^2\ar[dd]_{C_{+,n+k}^{n,k}}\ar[ll]_{\psi_{+,n+1}^{n,1}\ot\id}\ar[dr]_{\psi_{+,\mu}^{n,1}\ot\id}&&\sut{n-1}^1\otimes \sut{k}^2\ar[dd]_{C_{+,n+k-1}^{n-1,k}}\ar[ll]_{\psi_{+,n}^{n-1,1}\ot\id}&\ar[l]\cdots
\\
\cdots&&\sut{\mu}^1\ot\sut{k}^2\ar[dd]\ar[ur]^{\psi_{+,n}^{\mu,1}\ot\id}&&\sut{\mu}^1\ot\sut{k}^2\ar[dd]\ar[ur]^{\psi_{+,n-1}^{\mu,1}\ot\id}&&\cdots
\\
\cdots&\ar[l]\sut{n+k+1}^{\#}\ar[dr]_{\psi_{+,\mu}^{n+k+1,\#}}&&\sut{n+k}^{\#}\ar[ll]_{\psi_{+,n+k+1}^{n+k,\#}}\ar[dr]_{\psi_{+,\mu}^{n+k,\#}}&&\sut{n+k-1}^{\#}\ar[ll]_{\psi_{+,n+k}^{n+k-1,\#}}&\ar[l]\cdots
\\
\cdots&&\sut{\mu}^{\#}\ar[ur]^{\psi_{+,n+k}^{\mu,\#}}&&\sut{\mu}^{\#}\ar[ur]^{\psi_{+,n+k-1}^{\mu,\#}}&&\cdots
}
\end{equation}
where the vertical maps from $\sut{\mu}^1\ot\sut{k}^2$ to $\sut{\mu}^\#$ is $C_{+,\mu}^{\mu,k}$ in Lemma \ref{lem: contact map}.

By (\ref{eq: commutative diagram a}) and (\ref{eq: commutative diagram b}) in Lemma \ref{lem: contact map}, the squares in (\ref{eq: useful triangles 2}) involving vertical maps commute. Hence the vertical maps induce a morphism $C_+$ between unrolled exact couples. This induces a filtered chain map from $B^+(K_1)\ot \sut{k}^2$ to $B^+(K_1\# K_2)$.

Since $k$ is large, from Corollary \ref{cor: split}, there is a canonical embedding of $\khii(-Y_2,K_2)=\sut{\mu}^2$ into $\sut{k}^2$ by the inverse of $\psi_{+,\mu}^{k,2}$. By precomposing this embedding, we obtain a filtered chain map from $B^+(K_1)\ot \sut{\mu}^2$ to $B^+(K_1\# K_2)$. Then by (\ref{eq: commutative diagram c}), we know this filtered chain map is an identification on the first page. This implies that it induces an identification on each page and then on the total filtered chain complex.
\end{proof}

\subsection{The formula for null-homologous knots}
In this subsection, we combine Proposition \ref{prop: Kunneth diff} and the integral surgery formula to obtain rational surgery formulae for null-homologous knots. First, we do similar calculations as in \cite[Lemma 7.2]{Ozsvath2011rational}.

\blem\label{lem: homology cal}
Suppose $K_1\subset Y_1$ is a null-homologous knot and $K_2=O_{v/r}\subset Y_2=L(v,-r)$. Let $(Y_\#,K_\#)$ be the connected sum of $K_1$ and $K_2$ and suppose $K_\#$ is framed by the longitude of $K_2$. Suppose $(\mu_\bullet,\lambda_\bullet)$ is the meridian and the longitude of $K_\bullet$ for $\bullet\in\{1,2,\#\}$. Then $$H_1(Y_\#\backslash N(K_\#))\cong H_1(Y_1\backslash N(K_1)).$$Moreover, the order of $K_\#$ is $v$ and the intersection number $\partial S_\#\cdot \lambda_\#$ is $-r$, where $S_\#$ is the Seifert surface of $K_\#$.
\elem
\bpf
The knot $K_1$ is of order $1$ and $O_{v/r}$ is of order $v$. Then $$H_1(Y_1\backslash N(K_1))\cong H_1(Y_1)\oplus \mathbb{Z}\langle \mu_1\rangle\aand H_1(Y_2\backslash N(K_2))\cong\mathbb{Z}.$$We write $g_1=\mu_1$ and $g_2$ as the generator of $H_1(Y_2\backslash N(K_2))$. Then $\mu_2=v\cdot g_2$ and $\lambda_2=r\cdot g_2$. A calculation on the homology shows \begin{equation}\label{eq: homology cal}
    \begin{aligned}
    H_1(Y_\#\backslash N(K_\#))&\cong \bigg(H_1(Y_1\backslash N(K_1))\oplus H_1(Y_2\backslash N(K_2))\bigg)/(\mu_1,\mu_2)\\&\cong \bigg(H_1(Y_1)\oplus\mathbb{Z}\langle g_1,g_2\rangle \bigg)/(g_1=v\cdot g_2)\\&\cong H_1(Y_1)\oplus\mathbb{Z}\langle g_\#\rangle\\&\cong H_1(Y_1\backslash N(K_1)),
\end{aligned}
\end{equation}where we write $g_\#$ as the generator. We also write ${\rm pr}$ as the projection to the summand generated by $g_\#$. Then $${\rm pr}(\mu_\#)=v\cdot g_\#\aand {\rm pr}(\lambda_\#)=r\cdot g_\#.$$Thus, the knot $K_\#$ is of order $v$ and $\partial S_\#\cdot \lambda_\#=-r$.
\epf

\bcor\label{cor: grading shift of C}
Let $K_\bullet\subset Y_\bullet$ for $\bullet\in\{1,2,\#\}$ be defined as in Lemma \ref{lem: homology cal}. Suppose $C^{n,k}_{\pm,m+k}$ and $C^{\mu,k}_{\pm,\mu}$ are defined as in Lemma \ref{lem: contact map}. Then we have explicit formulae of the grading shifts of the maps as follows.
\begin{equation}\label{eq: c1 graded}
    C^{n,k}_{\pm,m+k}\bigg((\sut{n}^1,i\pm \frac{n-1}{2})\ot(\sut{k}^2,j\pm \frac{(k-1)v+r}{2})\bigg)\subset (\sut{n+k+1}^\#,iv+j\pm \frac{(n+k-1)v+r}{2}).
\end{equation}
\begin{equation}\label{eq: c2 graded}
    C^{\mu,k}_{\pm,\mu}\bigg((\sut{\mu}^1,i)\ot(\sut{k}^2,j\pm \frac{(k-1)v+r}{2})\bigg)\subset (\sut{n+k+1}^\#,iv+j).
\end{equation}
\ecor
\bpf
First, we compute the grading shift of $C_{\pm,\mu}^{\mu,k}$. From the homology calculation in Lemma \ref{lem: homology cal} and the graded version of (\ref{eq: connected sum}) in \cite[Proposition 5.15]{LY2021large}, we have\begin{equation}\label{eq: graded}
    (\sut{\mu}^\#,s)\cong\bigoplus_{s_1v+s_2=s}(\sut{\mu}^1,s_1)\otimes (\sut{\mu}^2,s_2),
\end{equation}where we take the direct sum over $s_1v+s_2=s$ because $g_\#=v\cdot g_1=g_2$ under the third isomorphism in (\ref{eq: homology cal}). From Lemma \ref{lem: bypass n,n+1,mu}, we know the grading shift of the map $\psi^{k,2}_{\pm,\mu}$ is $\mp\frac{(k-1)v+r}{2}$. Then from (\ref{eq: commutative diagram c}), we know the grading shift of $C_{\pm,\mu}^{\mu,k}$ is described in (\ref{eq: c2 graded}).

\begin{figure}[ht]
	\begin{overpic}[width=3in]{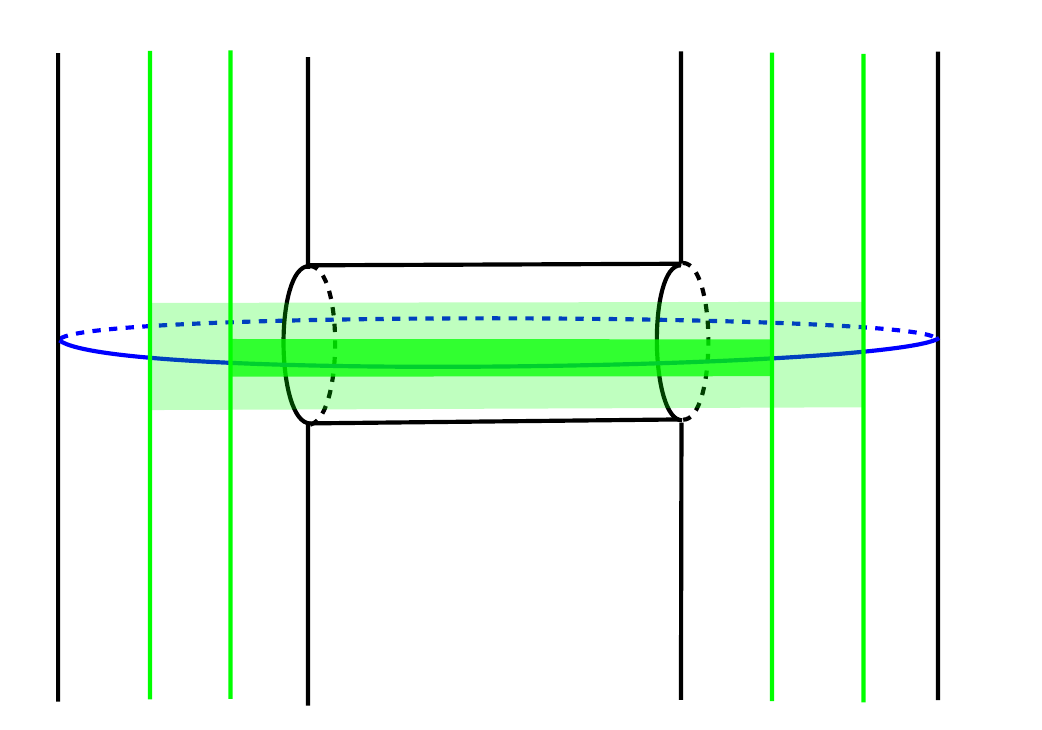}
		\put(20,66){$\partial S_1$}
		\put(69,66){$\partial S_2$}
		\put(11,66){$\partial S_1$}
		\put(80,66){$\partial S_2$}
		%\put(2,55){$\beta$}
		\put(2,38){$\al$}
		\put(43,48){$h^1$}
	\end{overpic}
\caption{The band sum for the case $v=2$. The two (green) shaded regions are the two bands.}\label{fig: connected sum 2}
\end{figure}
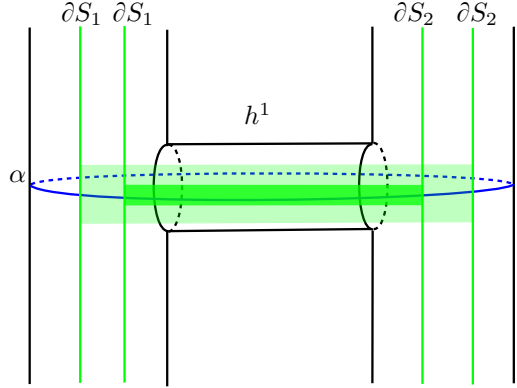

Also from Lemma \ref{lem: bypass n,n+1,mu} and Lemma \ref{lem: homology cal}, we know the grading shifts of $\psi_{\pm,\mu}^{n,1}$ and $\psi_{\pm,\mu}^{n+k,\#}$ are $\mp\frac{n-1}{2}$ and $\mp\frac{(n+k-1)v+r}{2}$, respectively. From (\ref{eq: commutative diagram a}) and (\ref{eq: c2 graded}), the expected grading shifts of $C^{n,k}_{\pm,m+k}$ are described in (\ref{eq: c1 graded}). Though in general $\psi_{\pm,\mu}^{n+k,\#}$ are not injective, we can still obtain (\ref{eq: c1 graded}) from the topological construction of $C^{n,k}_{\pm,m+k}$ in the proof of Lemma \ref{lem: contact map}. The proof is similar to the proof of \cite[Lemma 4.3]{li2019tau} and the only difference is that now the knot $K_2$ has order $v$ so that a (rational) Seifert surface of the connected sum knot $K_1\# K_2$ is obtained from one Seifert surface of $K_2$ and $p$ copies of Seifert surfaces of $K_1$ by $v$-many band sums. See Figure \ref{fig: connected sum 2}.
\epf

Then we provide an identification of bent complexes.

\bprop\label{prop: id of bent complex}
Let $K_\bullet\subset Y_\bullet$ for $\bullet\in\{1,2,\#\}$ be defined as in Lemma \ref{lem: homology cal}. Then for any grading $s$, there is an identification $$A(K_\#,s)=A(K_1,s^\p),$$where $s^\p$ is the unique grading satisfying $$|s-s^\p v|\le \frac{v-1}{2}.$$Moreover, we have the following commutative diagrams
\begin{equation*}
\xymatrix@R=3ex{
A(K_\#,s)\ar[rr]^{\pi^\pm(K_\#,s)}\ar[dd]_{=}&&B^\pm(K_\#,s)\ar[dd]^{=}\\
&&\\
A(K_1,s^\p)\ar[rr]^{\pi^\pm(K_1,s^\p)}&&B^\pm(K_1,s^\p)
}
\end{equation*}
\eprop
\bpf
From Lemma \ref{lem: dim}, we know $(\sut{\mu}^2,s_2)$ is nontrivial only for $|s_2|\le \frac{v-1}{2}$, for which the grading summand is $1$-dimensional. Due to the homology result in Lemma \ref{lem: homology cal}, we can apply the graded version of (\ref{eq: connected sum}) in (\ref{eq: graded}) to show that $$(\sut{\mu}^\#,s)\cong\bigoplus_{s_1v+s_2=s}(\sut{\mu}^1,s_1)\otimes (\sut{\mu}^2,s_2)\cong (\sut{\mu}^1,s^\p).$$Moreover, we have $$
\bigoplus_{k\ge 0}(\sut{\mu}^\#,s+kv)\cong \bigoplus_{k\ge 0}(\sut{\mu}^1,s^\p+k).$$Note that two sides of the isomorphism are underlying spaces of subcomplexes of $B^+(K_\#)$ and $B^+(K_1)$ since the orders of $K_\#$ and $K_1$ are $v$ and $1$, respectively. From Proposition \ref{prop: Kunneth diff}, the differentials $d_{+}$ on both sides are the same under the isomorphism. Similarly, we have$$\bigoplus_{k\le 0}(\sut{\mu}^\#,s+kv)\cong \bigoplus_{k\le 0}(\sut{\mu}^1,s^\p+k),$$and the differentials $d_-$ on both sides are the same. Hence we conclude the identification about the bent complex (\textit{cf.} Definition \ref{defn: complex A}). The commutative diagrams follow immediately.
\epf

\bthm\label{thm: rational surgery formula, bent complex}
Suppose $K\subset Y$ is a null-homologous knot. Suppose $u/v\in \mathbb{Q}\backslash\{0\}$. For any grading $s$, let $s^\p$ and $s^\pp$ be the unique gradings satisfying $$|s-s^\p v|\le \frac{v-1}{2}\aand |s+u-s^\pp v|\le \frac{v-1}{2}.$$ Then there exists an isomorphism 
$$\Xi_{u/v}:\bigoplus_{s\in\mathbb{Z}}H(B^+(s^\p))\xra{\cong}\bigoplus_{s\in\mathbb{Z}}H(B^-(s^\pp))$$as the direct sum of isomorphisms$$\Xi_{u/v,s}:H(B^+(s^\p))\xra{\cong}H(B^-(s^\pp))$$
so that 

$$I^\sharp(-Y_{-u/v}(K))\cong H\bigg(\cone(\pi^-+\Xi_{u/v}\circ \pi^+:\bigoplus_{s\in\mathbb{Z}}H(A(s^\p))\to  \bigoplus_{s\in\mathbb{Z}}H(B^-(s^\p)))\bigg).$$
% Let $$\mathbb{A}_{u/v}=\bigoplus_{i=0}^{u-1}\bigoplus_{s\in\mathbb{Z}}H(A(K,\lfloor\frac{i+us}{v}\rfloor))\aand \mathbb{B}_{u/v}=\bigoplus_{i=0}^{u-1}\bigoplus_{s\in\mathbb{Z}}H(B^-(s)).$$Let $D_{u/v}:\mathbb{A}_{u/v}\to \mathbb{B}_{u/v}$ be the map that sends $$a_w\in H(A(K,w))$$to $$\pi^-(a_w)+\Xi\circ \pi^+(a_w)\in H(B^-(K,w))\oplus H(B^-(K,w^\p))$$where $$w=\lfloor\frac{i+us}{v}\rfloor, w^\p=\lfloor\frac{i+u(s+1)}{v}\rfloor,$$and $$\Xi:H(B^+(K))\xra{\cong} H(B^-(K))$$ is an isomorphism. Then we have an isomorphism$$I^\sharp(-Y_{-u/v}(K))\cong H(\cone(D_{u/v})).$$
\ethm
\bpf
The statement is an analog of the rational surgery formula for $\widehat{HF}$ in \cite[Section 7.1]{Ozsvath2011rational}, where $\pi^-$ and $\Xi_{u/v}\circ \pi^+$ are analogs of $\hat{v}$ and $\hat{h}$. Let $m=\lfloor u/v\rfloor$ and $u/v=m+r/v$. Following the strategy at the start of this section, set $K_1=K$, and let $K_2=O_{v/r}$ to be a core knot in a lens space. Then $Y_{u/v}(K)$ is obtained by $m$-surgery on $K_\#=K_1\# K_2$. From Theorem \ref{thm: integral bent complex}, we know there exists an isomorphism$$\Xi_{m}^\#:\bigoplus_{s\in\mathbb{Z}}H(B^+(K_\#,s))\xra{\cong}\bigoplus_{s\in\mathbb{Z}}H(B^-(K_\#,s+mp-q))$$preserving the direct summand for each fixed $s\in\intg$ such that $$I^\sharp(-Y_{-u/v}(K))\cong H\bigg(\cone(\pi^{-,\#}+\Xi_m^\#\circ \pi^{+,\#}:\bigoplus_{s\in\mathbb{Z}}H(A(K_\#,s))\to  \bigoplus_{s\in\mathbb{Z}}H(B^-(K_\#,s)))\bigg).$$
Note that $(p,q)=(v,-r)$ so $mp-q=u$.

From Proposition \ref{prop: id of bent complex}, we replace complexes of $K_\#$ with complexes of $K_1$ to obtain the rational surgery formula, where $\Xi_{u/v}$ is induced by $\Xi_m^\#$ under the identification.
\epf

\brem
The rational surgery formula for a rationally null-homologous knot is more complicated but still doable. In such a case, the graded version of K\"{u}nneth formula is not enough and we need a torsion spin$^c$-like decomposition for sutured instanton homology (\textit{cf.} \cite{LY2021enhanced}, Remark \ref{rem: remove condition}, and Remark \ref{rem: remove condition2}).
\erem
% As mentioned in Remark \ref{rem: stronger}, the integral surgery formula via sutured instanton homology (Theorem \ref{thm: integral surgery formula}) is stronger than the one via bent complex (Theorem \ref{thm: integral bent complex}) in some cases. So we also state a rational surgery formula via sutured instanton homology.

% Suppose $K\subset Y$ is a rationally null-homologous knot with framing. Suppose $u/v$ is the surgery coefficient and $m=\lfloor u/v\rfloor$. Following the strategy at the start of this section, set $K_1=K$ and $K_2=O_{v/r}$ to be a core knot. Let $K^\p$ be the connected sum and we adapt notations in Section \ref{subsec: The connected sum with a core knot} for sutured instanton homology and bypass maps. From Theorem \ref{thm: integral surgery formula}, we have the following exact triangle

% \begin{equation*}\label{eq: connected sum triangle}
% 	\xymatrix{
% 	\mathbf{\Ga}_{\frac{2m+2k-1}{2}}^\#\ar[rr]^{\pi^\#}&&\mathbf{\Ga}_{m+2k-1}^\#\ar[dl]\\
% 	&I^\sharp(-Y_{-u/v}(K))\ar[ul]&
% 	}
% \end{equation*}
% where $$\pi^\#=\Psi^{n+k,\#}_{+,m+2k-1}\circ\psi^{\frac{2m+2k-1}{2},\#}_{-,m+k}+\Psi^{n+k,\#}_{-,m+2k-1}\circ\psi^{\frac{2m+2k-1}{2},\#}_{+,m+k}.$$
\section{The 0-surgery for knots in the 3-sphere}\label{sec: Knots inside $S^3$}

In this section, we deal with $0$-surgery for knots inside $S^3$. Recall that we have
$$\pi^+_{m,k}=\Psp{m+k}{m+2k-1}\circ\psm{\frac{2m+2k-1}{2}}{m+k}:\sut{\frac{2m+2k-1}{2}}\ra\sut{m+2k-1},$$
$$\pi^-_{m,k}=\Psm{m+k}{m+2k-1}\circ\psp{\frac{2m+2k-1}{2}}{m+k}:\sut{\frac{2m+2k-1}{2}}\ra\sut{m+2k-1},$$
and $\pi^{\pm,i}_{m,k}$ be the restriction of $\pi^{\pm}_{m,k}$ on $\sutg{\frac{2m+2k-1}{2}}{i}$. For knots inside $S^3$, we have a better description of the maps $\pi^{\pm,i}_{m,k}$ than in Lemma \ref{lem: the map pi^pm_m,k}.

\blem\label{lem: properties of pi^pm, sutured}
Suppose $K\subset S^3$ is a knot. Let $\tau=\tau_I(K)$ be defined as in Definition \ref{defn: tau}. For any fixed integer $m$ and large enough integer $k$, we have the following.
\begin{enumerate}
	\item When $i>\tau$, $\pi^{+,i}_{m,k}=0.$ When $i<-\tau$, $\pi^{-,i}_{m,k}=0$.
	\item When $i<\tau$, $\pi^{+,i}_{m,k}\neq0.$ When $i>-\tau$, $\pi^{-,i}_{m,k}\neq0$.
	\item When $i\leq-g(K)$, $\pi^{+,i}_{m,k}$ is an isomorphism. When $i\geq g(K)$, $\pi^{-,i}_{m,k}$ is an isomorphism. 
\end{enumerate}
\elem

% Due to the discussion in Step 1 of Section \ref{subsec: A strategy to prove the formula}, the maps $\pi^\pm_{m,k}$ are related to $\pi^\pm(s)$ constructed from the bent complex. Then we can restate the integral surgery formula for knots in $S^3$ as follows.

% \bthm[{Theorem \ref{thm: integral surgery formula for S^3}}]\label{thm: mapping cone restated}
% Suppose $K\subset S^3$ is a knot. For any fixed $m\in\intg\backslash\{0\}$ and any large enough integer $k$, we have
% $$I^{\sharp}(S^3_{-m}(K))\cong H({\rm Cone}(\pi^+_{m,k}+\pi_{m,k}^-)).$$
% \ethm

% We first establish some basic properties of the maps $\pi^{\pm}_{m,k}$.

\bpf
For part (1), we only prove the statement regarding $\pi^+_{m,k}$. The statement regarding $\pi^-_{m,k}$ follows from the symmetry between $K$ and $-K$, where $-K$ is the orientation reversal of $K$. Note that when we switch the orientation of the knot, the tau invariant remains the same, $\pi^{\pm}$ switches with each other, and the grading induced by the Seifert surface becomes the additive inverse. Let
$$\psm{\frac{2m+2k-1}{2},i}{m+k}=\psm{\frac{2m+2k-1}{2}}{m+k}|_{\sutg{\frac{2m+2k-1}{2}}{i}}.$$
We know that
$$\pi^{+,i}_{m,k}=\Psp{m+2k-1}{m+k}\circ\psm{\frac{2m+2k-1}{2},i}{m+k}.$$
From Lemma \ref{lem: bypass n+1,n,2n+1/2} we know
$$\im(\psm{\frac{2m+2k-1}{2},i}{m+k})\subset\sutg{m+k}{i+\frac{m+k-1}{2}}.$$
When $k$ is large enough so that $m+k$ is large, the map $\Psp{m+k}{m+2k-1}$  corresponds to the composition of $(k-1)$ many $U$-actions as in the construction of $\khi$ in Definition \ref{defn: khi minus}. By the definition of $\tau$ in Definition \ref{defn: tau}, we immediately conclude that 
$$\Psp{m+k}{m+2k-1}|_{\sutg{m+k}{j}}=0$$
whenever $j>\tau+\frac{m+k-1}{2}$. Hence, as a result, we have
$$\pi^{+,i}_{m,k}=0$$
when $i>\tau$.

For part (2), again we only prove the statement involving $\pi^+_{m,k}$. By the definition of $\tau$, and the correspondence between $\Psp{m+k}{m+2k-1}$ and $U^{k-1}$ on $\khi$, we know that when $k$ is large enough, there exists 
$$x\in\sutg{m+k-\tau+i}{\tau+\frac{m+k-\tau+i-1}{2}}$$ such that 
$$\Psp{m+k-\tau+i}{m+2k-1}(x)\neq 0.$$ 
Taking
$$y=\Psp{m+k-\tau+i}{m+k}(x)\in\sutg{m+k}{i},$$
we know that $$\Psp{m+k}{m+2k}(y)=\Psp{m+k-\tau+i}{m+2k-1}(x)\neq 0.$$
So it remains to show that $y\in\im (\psm{\frac{2m+2k-1}{2},i}{m+k})$. Indeed, from the construction of $y$ we know that 
$$\psp{m+k}{\mu}(y)=0.$$
Then from Lemma \ref{lem: comm diag for n+1 to n}, we know that 
$$\psm{m+k}{m+k-1}(y)= \psm{\mu}{m+k}\circ\psp{m+k}{\mu}(y)=0.$$
Hence by Lemma \ref{lem: bypass n+1,n,2n+1/2} we have
$$y\in \ke(\psm{m+k}{m+k-1})=\im (\psm{\frac{2m+2k-1}{2},i}{m+k}).$$

Part (3) is a restatement of Lemma \ref{lem: the map pi^pm_m,k}, part (3).
\epf

Next we study the $0$-surgery for knots inside $S^3$. The main obstruction to applying the proof of the integral surgery formula in \cite[Section 3.2]{LY2022integral1} to the $0$-surgery is that $\pi^{+}_{m,k}$ and $\pi^{-}_{m,k}$ have the same grading shift. Then $$H(\cone(c_{1}\pi^{+}_{m,k}+c_{2}\pi^{-}_{m,k}))$$may depend on the scalars. If either map vanishes, then the homology is still independent of the scalars. However, this is not true in general. Fortunately, we can make use of the $\mathbb{Z}$-grading on $I^\sharp(S_0^3(K))$ in (\ref{eq: z grading}). Note that one of the restrictions of $\pi_{m,k}^\pm$ on a single grading vanishes.

\bthm[$0$-surgery formula]\label{thm: 0surgery}
Suppose $K\subset S^3$ is a knot with $\tau_I(K)\le 0$. Suppose $A(s),B^\pm(s)$ and $\pi^\pm(s):A(s)\to B^\pm(s)$ are complexes and maps constructed in Definition \ref{defn: complex A}. For any $s\in\mathbb{Z}\backslash\{0\}$, there exists an isomorphism $$\Xi_{0,s}:H(B^+(s))\to H(B^-(s))$$such that $I^\sharp(-S_0^3(K),s)$ is isomorphic to  $$H\bigg(\cone(\pi^-(s)+\Xi_{0,s}\circ \pi^+(s):H(A(s))\to H(B^-(s)))\bigg).$$If $\tau_I(K)\neq 0$, then the same result also applies to $s=0$.
\ethm

\bpf
From Lemma \ref{lem: surgery triangles}, we have a long exact sequence $$\cdots\to\sut{\mu}\xra{A_{-1}}  \sut{-1}\to I^\sharp(-S_0^3(K))\to \sut{\mu}\to\cdots$$By the same reasoning in Lemma \ref{lem: -1 surgery and bypass}, we have \begin{equation}\label{eq: scalars}
    A_{-1}=c_1\psi_{+,-1}^\mu+c_2\psi_{-,-1}^\mu.
\end{equation}
Following the construction of the gradings induced by Seifert surfaces, the maps in the long exact sequence are all grading-preserving. We consider the following octahedral diagram that is used in \cite[Section 3.2]{LY2022integral1}.

\begin{equation}\label{eq: octahedral 3}
    \xymatrix@C=4ex{
    &&&\dehny{m}\ar@{..>}[dr]^{\psi}\ar[rr]^{}&&\sut{\mu}\ar[ddr]^{\begin{aligned}
        \psi_{+,m-1}^\mu\\+\psi_{-,m-1}^\mu
    \end{aligned}}&
    \\
    &&\sut{m-1}\ar[dr]^{\begin{aligned}
       (\Psi_{+,m-1+k}^{m-1},\\\Psi_{-,m-1+k}^{m-1}) 
    \end{aligned}}\ar[ur]^{}&&\sut{\frac{2m+2k-1}{2}}\ar@{..>}[ddrr]^{\phi}\ar[ur]^{l^\p}&&
    \\
 &&&\sut{m-1+k}\oplus\sut{m-1+k}\ar[drrr]_{\begin{aligned}
     \Psi_{-,m-1+2k}^{m-1+k}\\-\Psi_{+,m-1+2k}^{m-1+k}
 \end{aligned}}\ar[ur]_{h^\p}&&&\sut{m-1}
    \\
    \sut{\mu}\ar[uurr]^{\begin{aligned}
        \psi_{+,m-1}^\mu\\+\psi_{-,m-1}^\mu
    \end{aligned}}\ar[urrr]_{\begin{aligned}
        (\psi_{-,m-1+k}^\mu,\\\psi_{+,m-1+k}^\mu)
    \end{aligned}}&&&&&&\sut{m-1+2k}\ar[u]^{l}
    }
     \end{equation}
where $$h^\p=\psi_{-,\frac{2m+2k-1}{2}}^{m+k-1}-\psi_{+,\frac{2m+2k-1}{2}}^{m+k-1}.$$
When $m=0$, all maps are homogeneous, so we could consider the diagram grading-wise. Note that we may not know $c_1=c_2=1$ in (\ref{eq: scalars}), but we can add scalars to other maps to keep the diagram still being commutative. Following the same strategy in \cite[Section 3.2]{LY2022integral1}, we obtain for any $s\in\mathbb{Z}$,$$I^\sharp(-S_0^3(K),s)\cong H(\cone(c_{3}\pi_{0,k}^{+,i}+c_{4} \pi_{0,k}^{-,i}))$$for some scalars $c_3,c_4$.

When $\tau_I(K)\le 0$, from Lemma \ref{lem: properties of pi^pm, sutured}, we know for any $i\in\mathbb{Z}$, either $\pi_{0,k}^{+,i}$ or $\pi_{0,k}^{-,i}$ vanishes, and hence $H(\cone(c_{3}\pi_{0,k}^{+,i}+c_{4} \pi_{0,k}^{-,i}))$ is independent of the scalars. Then we have

$$\begin{aligned}
    I^\sharp(-S_0^3(K),s)&\cong H(\cone(c_{3}\pi_{0,k}^{+,i}+c_{4} \pi_{0,k}^{-,i}))\\&\cong H(\cone(\pi_{0,k}^{+,i}+\pi_{0,k}^{-,i}))\\&\cong H(\cone(\pi^-(s)+\Xi_{0,s}\circ \pi^+(s))),
\end{aligned}$$where $\Xi_{0,s}$ is constructed similarly to $\Xi_{m,s}$ in Theorem \ref{thm: integral bent complex} for $m\neq 0$.
\epf
\brem
From Lemma \ref{lem: basic properties of tau} part (2), we may pass to the mirror knot to satisfy the assumption $\tau_I\le 0$ in Theorem \ref{thm: 0surgery}.
\erem

Baldwin-Sivek also studied framed instanton homology with a twisted bundle for $0$-surgery, which is denoted by $I^\sharp(S_0^3(K),\mu)$, where $\mu$ is the meridian of the knot. There is also a $\mathbb{Z}$-grading on this homology induced by the Seifert surface and we also have a long exact sequence $$\cdots\to \sut{\mu}\xra{c_1^\p\psi_++c_2^\p\psi_-}  \sut{-1}\to I^\sharp(-S_0^3(K),\mu)\to \sut{\mu}\to\cdots$$such that all maps are grading-preserving (the coefficients $c_1^\p$ and $c_2^\p$ may be different from $c_1$ and $c_2$). Thus, we can use the similar octahedral diagram grading-wise to prove the result in Theorem \ref{thm: 0surgery} when replacing $I^\sharp(-S_0^3(K))$ with $I^\sharp(-S_0^3(K),\mu)$. As a result, we obtain the following corollary. We also write $I^\sharp(-S_0^3(K),\mu,i)$ for the summand of $I^\sharp(-S_0^3(K),\mu)$ having grading $i$.

\bcor\label{cor: same dimension for 0surgery}
Suppose $K\subset S^3$ is a knot. For any $s\in\mathbb{Z}\backslash\{0\}$, we have $$I^\sharp(-S_0^3(K),s)\cong I^\sharp(-S_0^3(K),\mu,s).$$
\ecor

\section{Surgeries on Borromean knots}\label{sec: Surgeries on Borromean knots}

In this section, we study surgeries on the connected sums of Borromean knots. 

\subsection{The Borromean knot}

First, we compute the $\khii$ for the Borromean knot. Let $T^3=S^1_1\times S^1_2\times S^1_3$. Let $Y$ be the result of a $0$-surgery along $S^1_1\subset T^3$ with respect to the surface framing induced by $T^2=S^1_1\times S^1_2$. Note $Y=\#^2(S^1\times S^2)$. Let $K$ be the core knot of the $0$-surgery, which is another description of the Borromean knot according to \cite[Section 9]{ozsvath2004holomorphicknot}. The knot $K$ bounds a genus-one Seifert surface $S=S^1_2\times S^1_3\backslash N(S^1_1)$. Let $\mu\subset \partial (Y\backslash N(K))$ be the meridian, and let $\lambda=\partial S\subset \partial (Y\backslash N(K))$ be the longitude.

\blem\label{lem: KHI of Borromean knot}
We have the following
\begin{equation*}
	\khii(Y,K,i)\cong \begin{cases}
		\mathbb{C}&|i|=1\\
		\mathbb{C}^2&i=0\\
		0&{\rm otherwise}
	\end{cases}
\end{equation*}
\elem

\bpf
We first figure out $\khii(Y,K)=\shi(Y\backslash N(K),\Gamma_{\mu})$. Using an annulus to form an auxiliary surface, we know from \cite[Lemma 5.2]{kronheimer2010knots} that a closure of $(Y\backslash N(K),\Gamma_{\mu})$ can be described as $S^1\times \Sigma_2$ where $\Sigma_2$ is a closed surface of genus $2$, obtained by gluing two once-punctured tori together. From the proof of \cite[Lemma 5.2]{kronheimer2010knots}, there is a pair of simple closed curves $\al,\be\subset \Sigma_2$ such that $\al\cdot \be=1$, the torus $S^1\times \alpha$ is the distinguishing surface of the closure, and $\beta$ serves as the $w_2$ that specifies the bundle over $S^1\times \Sigma_2$. By construction,
$$\shi(Y\backslash N(K),\Gamma_{\mu})\cong{\rm Eig}(I^{\be}(S^1\times \Sigma_2),\mu({\rm pt}),2),$$
where ${\rm Eig}(I^{\be}(S^1\times \Sigma_2),\mu({\rm pt}),2)$ means the generalized eigenspace of $\mu({\rm pt})$ on $I^{\be}(S^1\times \Sigma_2)$ with eigenvalue $2$.

On the other hand, take $T^2=\partial (Y\backslash N(K))$ and take the Seifert framing of $K$ on $T^2$. Let $M=[0,1]\times T^2$ and let $\Gamma_{\mu,\mu}$ be the suture on $\partial M$ which consists of two meridians on each boundary component of $M$. We can use an annulus to close up each boundary component of $(M,\Ga_{\mu,\mu})$ separately. A construction similar to that of \cite[Lemma 5.2]{kronheimer2010knots} implies that a closure of $(M,\Ga_{\mu,\mu})$ can be described as $S^1\times \Sigma_2$, and there are two pairs of curves $\al,\be,\al',\be'$ on $\Sigma_2$ such that $\al\cdot \be=1$, $\al'\cdot\be'=1$, the surface $S^1\times (\al\cup\al')$ is the distinguishing surface of the closure, and $\be\cup\beta'$ serves as the $w_2$. Furthermore, the two pairs $(\al,\be)$ and $(\al',\be')$ come from closing up two boundary components of $M$, so they are disjoint from each other. We know the following
$$\shi(M,\Ga_{\mu,\mu})\cong {\rm Eig}(I^{\be\cup\be'}(S^1\times \Sigma_2),\mu({\rm pt}),2).$$
Since $\be$ and $\be\cup\be'$ both represent primitive homology classes on $\Sigma_2$, there is an orientation-preserving diffeomorphism $h:\Sigma_2\to\Sigma_2$
such that $h([\be])=[\be]+[\be']$.
As a result, the map $h$ extends to a diffeomorphism between closures and we conclude
$$\shi(Y\backslash N(K),\Gamma_{\mu})\cong \shi(M,\Ga_{\mu,\mu}).$$

Observe that there is a sutured manifold decomposition
$$(M,\ga)\stackrel{A}{\leadsto}(V,\ga^6),$$
where $A=[0,1]\times\mu\subset M$ is a product annulus and $V\cong S^1\times D^2$ is a solid torus with $\ga^6$ consisting of six longitudes of $V$. From \cite[Proposition 1.4]{li2019direct} and \cite[Proposition 6.7]{kronheimer2010knots} we know that
$$\shi(Y\backslash N(K),\Gamma_{\mu})\cong \shi(M,\Ga_{\mu,\mu})\cong \mathbb{C}^4.$$

Now we compute the dimension of each graded part. Since $g(K)=1$, we know $\khii(Y,K,i)=0$ for $|i|>1$. For $|i|=1$, since $K\subset Y$ is fibered (the complement is $S^1\times (T^2\backslash D^2)$), we have
$$\khii(Y,K,1)\cong \khii(Y,K,-1)\cong \mathbb{C}.$$
As a result, we conclude that $\khii(Y,K,0)\cong \mathbb{C}^2.$
\epf

On connected sums of $S^1\times S^2$, the circles $S^1\times\{{\rm pt}\}$ induce nontrivial actions on the framed instanton homology. In particular, we have the following lemma.

\blem[{\cite[Section 7.8]{scaduto2015instanton} and \cite[Theorem 7.16]{donaldson2002floer}}]\label{lem: S^1 action}
Suppose $\widehat{Y}$ is the connected sum of copies of $S^1\times S^2$. There is a canonical action of $\Lambda^*H_1(\widehat{Y})$ on $I^{\sharp}(\widehat{Y})$, making $I^{\sharp}(\widehat{Y})$ a rank-one free module over $\Lambda^*H_1(\widehat{Y})$.
\elem

Since $Y=\# ^2 S^1\times S^2$, Lemma \ref{lem: S^1 action} implies
$$I^{\sharp}(Y)=\Lambda^*H_1(Y;\mathbb{C})=\mathbb{C}[x_1,x_2]\slash(x_1x_2+x_2x_1,x_1^2,x_2^2)=\mathbb{C}\lgl 1,x_1,x_2,x_1x_2\rgl.$$

Note on $Y$ we can pick two circles whose $\mu$-actions correspond to the multiplication of $x_1$ and $x_2$ on $\dehny{}=\Lambda^*H_1(-Y;\mathbb{C})$. We can pick these two circles to be disjoint from the Borromean knot $K$. Since the $\mu$-action of a circle commutes with all cobordism maps and all $\mu$-actions of surfaces and points, we know that there is an action of $\Lambda^*H_1(-Y;\mathbb{C})$ on $\sutg{*}{i}$ for any $*\in\mathbb{Q}\cup\{\mu\}$ and any grading $i$. This makes $\sutg{*}{i}$ a $\Lambda^*H_1(-Y;\mathbb{C})$-module and all bypass maps and surgery maps are module morphisms. We have the following structure lemma.
\begin{lem}\label{lem: gamma_n for Borromean knot}
	Suppose $K\subset Y$ is the Borromean knot. Then for any integer $n\geq 2$, we have an identification
	\begin{equation*}
		\sut{n}=\shi(-Y\backslash N(K),-\Gamma_n,i)=\begin{cases}
	\mathbb{C}\lgl x_1x_2\rgl&|i|=\frac{n+1}{2}\\
	\mathbb{C}\lgl x_1,x_2,x_1x_2\rgl&|i|=\frac{n-1}{2}\\
	\Lambda^*H_1(Y;\mathbb{C})&|i|<\frac{n-1}{2}\\
	0&{\rm otherwise}
\end{cases}
	\end{equation*}
\end{lem}

\bpf
The structure of $\sut{n}$ for large $n$ is understood by Lemma \ref{lem: structure of Gamma_n} so it suffices to work out the structures of $\sut{2}$ and $\sut{3}$. By Lemma \ref{lem: surgery triangles} and Lemma \ref{lem: bypass n,n+1,mu}, there are exact triangles
\begin{center}
\begin{minipage}{0.4\textwidth}
	\begin{equation*}
\xymatrix{
\sut{2}\ar[rr]^{H_2}&&\sut{3}\ar[dl]^{F_{3}}\\
&\dehny{}\ar[ul]^{G_2}&
}	
\end{equation*}
\end{minipage}
\begin{minipage}{0.4\textwidth}
		\begin{equation*}
\xymatrix{
\sut{2}\ar[rr]^{\psi^{2}_{\pm,3}}&&\sut{3}\ar[dl]^{\psi^{3}_{\pm,\mu}}\\
&\sut{\mu}\ar[ul]^{\psi^{\mu}_{\pm,2}}&
}	
\end{equation*}\end{minipage}
\end{center}

From Lemma \ref{lem: structure of Gamma_n} part (4), we know that $F_3$ is surjective so
$$\dim\sut{3}=\dim\sut{2}+\dim \dehny{}=\dim\sut{2}+4.$$
Since $\dim \sut{\mu}=4$, we know that the last two exact triangles split as well and in particular, the maps $\psi^2_{\pm,3}$ are both injective. From Lemma \ref{lem: structure of Gamma_n} part (4), we know that
\begin{equation}\label{eq: gamma_3,0}
	\sutg{3}{0}\cong \dehny{}=\Lambda^*H_1(-Y;\mathbb{C}).
\end{equation}
Hence from Lemma \ref{lem: KHI of Borromean knot}, Lemma \ref{lem: structure of Gamma_n} part (2), and Lemma \ref{lem: bypass n,n+1,mu} we know that
$$\dim \sutg{3}{\pm1}=\dim\sutg{2}{\pm\frac{1}{2}}=\dim\sutg{3}{0}-\dim\sutg{\mu}{\mp1}=3.$$
Similarly,
$$\dim \sutg{3}{\pm2}=\dim\sutg{2}{\pm\frac{3}{2}}=\dim\sutg{3}{\pm1}-\dim\sutg{\mu}{0}=1.$$

Since the isomorphism in (\ref{eq: gamma_3,0}) is induced by a cobordism map, it is an isomorphism between modules. We have an injective module morphism
$$\psp{2}{3}:\sutg{2}{\frac{1}{2}}\to\sutg{3}{0}.$$
We have the following claim.

{\bf Claim}. There is a unique $3$-dimensional submodule inside $\Lambda^*H_1(-Y;\mathbb{C})$.
\bpf[Proof of Claim.]
Indeed, suppose $\mathcal{M}\subset \Lambda^*H_1(-Y;\mathbb{C})$ is a $3$-dimensional submodule. Assume that $1+a\in \mathcal{M}$, where $a$ is spanned by $x_1$, $x_2$, and $x_1x_2$. Then, note that $x_1x_2=x_1x_2(1+a)\in\mathcal{M}$. Also $x_1(1+a)$ is of the form $x_1+c\cdot x_1x_2$ for some $c\in\mathbb{C}$ so we know $x_1\in\mathcal{M}$ and similarly $x_2\in\mathcal{M}$. As a result $1\in\mathcal{M}$ so $\mathcal{M}$ must be all of $\Lambda^*H_1(-Y;\mathbb{C})$. Hence we conclude that $\mathcal{M}$ does not have an element of the form $1+a$ so the only possibility is that $\mathcal{M}=\mathbb{C}\lgl x_1,x_2,x_1x_2\rgl$.
\epf

From the claim we know that 
$$\sutg{3}{1}\cong\sutg{2}{\frac{1}{2}}\cong \mathbb{C}\lgl x_1,x_2,x_1x_2\rgl.$$
From the injectivity of the map $\psp{2}{3}:\sutg{2}{\frac{3}{2}}\to\sutg{3}{1}$
we can conclude similarly that
$$\sutg{3}{2}\cong\sutg{2}{\frac{3}{2}}\cong \mathbb{C}\lgl x_1x_2\rgl.$$
\epf

\begin{cor}\label{cor: bypass maps for Borromean knot}
Under the description of Lemma \ref{lem: gamma_n for Borromean knot}, the bypass maps between $\sut{n}$ and $\sut{n+1}$ for $n\geq 2$ are described as follows.

\begin{itemize}
	\item If $i\geq0$, the map
	$$\psi^{n}_{\pm,n+1}:\sutg{n}{\pm i}\ra\sutg{n+1}{i\mp\frac{1}{2}}$$
	is an inclusion or the identity if the domain and the range are the same.
	\item If $i\leq0$, the map
	$$\psi^{n}_{\pm,n+1}:\sutg{n}{\pm i}\ra\sutg{n+1}{i\mp\frac{1}{2}}$$
	is the identity.
\end{itemize}
Moreover, the module structure on $\sut{\mu}$ is trivial, \textit{i.e.}, the module multiplications of $x_1$ and $x_2$ are both zeros.
% \begin{enumerate}
% 	\item When $i=\frac{n-1}{2},\frac{n+1}{2}$, the bypass map
% 	$$\psi^{n}_{\mp,n+1}:\sutg{n}{\pm i}\to \sutg{n+1}{\pm i \pm \frac{1}{2}}$$
% 	is the identity.
% 	\item When $i<0$, the bypass map
% 	$$\psi^{n}_{\pm,n+1}:\sutg{n}{ i\pm\frac{1}{2}}\to\sutg{n+1}{i}$$
% 	is the inclusion (or the identity if domain and codomain are the same).
% \end{enumerate}	

\end{cor}

\bpf
Note that all the bypass maps are module morphisms. The description of the bypass maps is straightforward from the proof of Lemma \ref{lem: gamma_n for Borromean knot}. For the module structure of $\sut{\mu}$, we know that $\sutg{\mu}{1}\cong \sutg{\mu}{-1}\cong \mathbb{C}$ so the structure must be zero. Also from Lemma \ref{lem: bypass n,n+1,mu}  we know
$$\sutg{\mu}{0}\cong H\bigg({\rm Cone}(\mathbb{C}\lgl x_1x_2\rgl\hookrightarrow\mathbb{C}\lgl x_1,x_2,x_1x_2\rgl)\bigg)$$
so the module structure on $\sutg{\mu}{0}$ is also trivial. 
\epf

Using the integral surgery formula Theorem \ref{thm: integral surgery formula} and the dual knot formula in \cite[Section 3.4]{LY2022integral1}, we can compute $I^{\sharp}(-Y_{-n}(K))$ and $\sut{n}$ for any $n\in\mathbb{Z}$ ($n\neq 0$ for $I^\sharp$). Since we will also deal with the connected sum of the Borromean knots, we omit the calculation here.

\subsection{The connected sums}\label{subsec: connected sum}

In this subsection, we compute the surgeries along $g$ copies of connected sums of $K\subset Y$. According to \cite[Section 5.2]{Ozsvath2008integral}, these surgeries give rise to nontrivial circle bundles over $\Sigma_g$. Write 
$$K^g=\#^gK\subset Y^g=\#^gY=\#^{2g}S^1\times S^2.$$
Note that the genus of $K^g$ is exactly $g$. For the rest of this subsection, for $*\in\intg\cup\{\mu\}$, write $\Ga_{*}$ for the suture on $Y^g\backslash N(K^g)$ and write $\sut{*}$ the corresponding sutured instanton homology. The connected sum formula (\ref{eq: connected sum}) for instanton knot homology gives rise to the following.
\bcor\label{cor: KHI of K^g}
We have the following
\begin{equation*}
	\dim \khii(Y^g,K^g,i)=\dim\sutg{\mu}{i}=\binom{2g}{g+i}.
\end{equation*}
Moreover, the module structure of $\sut{\mu}$ is trivial.
\ecor

Note that, from Lemma \ref{lem: S^1 action}, we know that
$$I^{\sharp}(-Y^g)=\Lambda^*H_1(-Y^g;\mathbb{C})=\mathbb{C}\lgl x_1,\dots,x_{2g}\rgl\slash(x_ix_j+x_jx_i).$$
For any $k\in[0,2g]\cap\intg$ write
\begin{equation}\label{eq: M_2g,k}
	\mathcal{M}_{2g,k}={\rm Span}\{\Pi_{j=1}^lx_{i_j}~|~l\geq k,~1\leq i_1<\dots<i_l\leq 2g\}.
\end{equation}
Note that $$\mathcal{M}_{2g,0}=\Lambda^*H_1(-Y^g;\mathbb{C})\aand \mathcal{M}_{2g,2g}\cong\mathbb{C}.$$It is straightforward to check that
$$\dim \mathcal{M}_{g,k}=\sum_{j=k}^{2g}\binom{2g}{k}.$$

\bdefn
Suppose $\mathcal{M}$ is a module over $\Lambda^*H_1(-Y^g;\mathbb{C})$. We say $\mathcal{M}$ is of {\bf degree $k>0$} if any monomial of degree at least $k+1$ annihilates $\mathcal{M}$ and there exists a monomial of degree $k$ acting nontrivially on $\mathcal{M}$. We say $\mathcal{M}$ is of {\bf degree $0$} if the module structure is trivial.
\edefn

\blem\label{lem: degree of modules}
Suppose $\mathcal{A}$, $\mathcal{B}$, and $\mathcal{C}$ are three modules over $\Lambda^*H_1(-Y^g;\mathbb{C})$ such that there is an exact triangle
\begin{equation*}
	\xymatrix{
	\mathcal{A}\ar[rr]^f&&\mathcal{B}\ar[dl]^g\\
	&\mathcal{C}\ar[lu]^h&\\
	}
\end{equation*}
where the three maps $f$, $g$, and $h$ are all module morphisms. Suppose further that $\mathcal{A}$ is of degree $k$ and $\mathcal{C}$ is of degree $0$. Then the degree of $\mathcal{B}$ is at most $k+1$.
\elem
\bpf
Suppose, on the contrary, that $\mathcal{B}$ is of degree $k+2$. Assume, without loss of generality, that there exists $b\in\mathcal{B}$ such that 
$$\bigg(\mathop{\prod}^{k+2}_{j=1}x_j\bigg)\cdot b\neq 0.$$

Suppose that $g(x_{k+2}\cdot b)\neq 0$. Then $x_{k+2}\cdot g(b)=g(x_{k+2}\cdot b)\neq 0$ and this contradicts the assumption that $\mathcal{C}$ is of degree $0$. As a result, there exists $a\in\mathcal{A}$ such that $f(a)=x_{k+2}\cdot b$. Then we have
$$f\bigg(\mathop{\prod}^{k+1}_{j=1}x_j\cdot a\bigg)=\bigg(\mathop{\prod}^{k+1}_{j=1}x_j\bigg)\cdot f(a)=\bigg(\mathop{\prod}^{k+2}_{j=1}x_j\bigg)\cdot b\neq 0.$$
As a result, we have
$$\bigg(\mathop{\prod}^{k+1}_{j=1}x_j\bigg)\cdot a\neq 0,$$
which contradicts the assumption that $\mathcal{A}$ has degree at most $k$.
\epf

\blem\label{lem: gamma_n for connected sums of Borromean knot}
For any $g\geq 1$, $n\geq 2g$, and any grading $i$, we have the following. 
\begin{equation*}
	\sutg{n}{i}\cong\begin{cases}
		\mathcal{M}_{2g,|i|+g-\frac{n-1}{2}}&|i|\ge \frac{n-1}{2}-g\\
		\Lambda^*H_1(-Y^g;\mathbb{C})&|i|\leq \frac{n-1}{2}-g\\
		0&{\rm otherwise}
	\end{cases}
\end{equation*}
\elem
\bpf
Again we only deal with $\sut{2g}$ and $\sut{2g+1}$. We prove this lemma by three claims. Note that, from Lemma, \ref{lem: structure of Gamma_n} part (4), we have 
$$\sutg{2g+1}{0}\cong I^{\sharp}(-Y^g)\cong \Lambda^*H_1(-Y^g;\mathbb{C}).$$ 
Also from Lemma \ref{lem: structure of Gamma_n} part (2) we know that $\sutg{2g}{i}\cong\sutg{2g+1}{i\pm\frac{1}{2}}$
for $\pm i\geq 0$.

{\bf Claim 1}. For $i>0$, the degree of $\sutg{2g}{\pm i}$ is at most $\frac{4g-1}{2}-i$.

\bpf[Proof of Claim 1.]
We only deal with $\sutg{2g}{i}$. The argument for $\sutg{2g}{-i}$ is similar. First from Lemma \ref{lem: structure of Gamma_n} part (2) we know that $\sutg{2g}{i}\cong\sutg{2g+1}{i\pm\frac{1}{2}}$ 
then we know from Lemma \ref{lem: bypass n,n+1,mu} that there exists an exact triangle
\begin{equation*}
	\xymatrix{
	\sutg{2g}{i}\ar[rr]&&\sutg{2g+1}{i-\frac{1}{2}}\cong\sutg{2g}{i-1}\ar[dl]\\
	&\sutg{\mu}{i-\frac{2g+1}{2}}\ar[ul]&
	}
\end{equation*}
Hence we can apply Lemma \ref{lem: degree of modules} to carry out an induction from the top grading of $\sut{2g}$ and the fact that $\sut{\mu}$ has degree $0$. The starting point is the top grading $\frac{4g-1}{2}$ for which we have $\sutg{2g}{\frac{4g-1}{2}}\cong\sutg{\mu}{g}\cong\mathbb{C}$. So clearly it has degree $0$.
\epf

{\bf Claim 2}. For $i>0$, we have the following.
$$\dim\sutg{n}{i}=\sum_{j=0}^{\frac{4g-1}{2}-i}\binom{2g}{j}.$$

\bpf[Proof of Claim 2.] From Lemma \ref{lem: structure of Gamma_n} part (4) we know that $F_{2g+1}$ is surjective. As a result, we have
$$\dim\sut{2g+1}-\dim\sut{2g}=\dim\dehny{}^g=2^{2g}=\dim\sut{\mu}.$$

% where $F_{2g+1}$ fits into the exact triangle
% \begin{equation*}
% 	\xymatrix{
% 	\sut{n}\ar[rr]&&\sut{n+1}\ar[ld]\\
% 	&\dehny{}^g\ar[lu]&
% 	}
% \end{equation*}

Hence the exact triangles
\begin{equation*}
	\xymatrix{
	\sut{n}\ar[rr]^{\psi^{2g}_{\pm,2g+1}}&&\sut{n+1}\ar[ld]\\
	&\sut{\mu}\ar[lu]&
	}
\end{equation*}
both split, which means the map $\psi^{2g}_{\pm,2g+1}$ is injective when restricted to the grading $\pm i$ for $i>0$. As a result, we can obtain the claim by induction and applying Corollary \ref{cor: KHI of K^g} and the fact
$$\dim\sutg{2g}{\pm\frac{4g-1}{2}}=\dim\sutg{\mu}{\pm g}=1=\binom{2g}{0}.$$
\epf

{\bf Claim 3}. The module $\mathcal{M}_{2g,k}$ is the only submodule of $\Lambda^*H_1(-Y^g;\mathbb{C})$ that has degree at most $2g-k$ and has dimension
$$\sum_{j=k}^{2g}\binom{2g}{k}.$$

The proof of the claim is straightforward. Note that there is a sequence of injective maps
$$\sutg{2g}{i}\hookrightarrow\sutg{2g+1}{i-\frac{1}{2}}\cong\sutg{2g}{i-1}\dots\hookrightarrow\sutg{2g+1}{0}\cong \Lambda^*H_1(-Y^g;\mathbb{C})$$
Hence the lemma follows from the above three claims.
\epf

From the proof of the above lemma, we also know the following
\bcor\label{cor: bypass maps for the connected sum of Borromean knots}
For any $g\geq 1$, $n\geq 2g$, and grading $i$, we have the following.
\begin{itemize}
	\item If $i\geq0$, the map
	$$\psi^{n}_{\pm,n+1}:\sutg{n}{\pm i}\ra\sutg{n+1}{i\mp\frac{1}{2}}$$
	is an inclusion or the identity if the domain and the range are the same.
	\item If $i\leq0$, the map
	$$\psi^{n}_{\pm,n+1}:\sutg{n}{\pm i}\ra\sutg{n+1}{i\mp\frac{1}{2}}$$
	is the identity.
\end{itemize}
\ecor

Again, based on Lemma \ref{lem: gamma_n for connected sums of Borromean knot} and Corollary \ref{cor: bypass maps for the connected sum of Borromean knots} and using the integral surgery formula and the dual knot formula in \cite[Section 3.4]{LY2022integral1}, we are able to compute $\sut{n}$ and $I^{\sharp}(-Y_{-n}(K))$ for any $n\in\intg$ ($n\neq 0$ for $I^\sharp$). Here we only present the computation for $I^{\sharp}(-Y_{-n}(K))$.

\bpf[Proof of Theorem \ref{thm: circle bundle computation}]
The manifold $Y^g_m$ is obtained from $Y^g=\# ^{2g}S^1\times S^2$ by $m$-surgery on the connected sum of the Borromean knot $K^g$. Note that as the Borromean knot $K\subset \#^2 S^1\times S^2$, we also know $Y_m^g(K^g)$ is diffeomorphic to $Y_{-m}^g(K^g)$. So it suffices to compute $Y_{-m}^g(K^g)$ for $m>0$. 

Since $\dim \sut{\mu}=\dim\dehny{}^g$, we know that all the differentials in the bent complexes are trivial. If $m\ge 2g-1$, then the argument follows directly from the large surgery formula in Proposition \ref{prop: truncated integral surgery formula}.

For smaller $m$, we can also use the truncation of the integral surgery formula in Proposition \ref{prop: truncated integral surgery formula} to make the computation easier. Suppose $k$ is large enough. From Lemma \ref{lem: vanishing grading}, Lemma \ref{lem: bypass n+1,n,2n+1/2}, and Lemma \ref{lem: gamma_n for connected sums of Borromean knot}, we know that the following two exact triangles both split.
\begin{equation*}
\xymatrix{
\sut{m+k-1}\ar[rr]&&\sut{\frac{2m+2k-1}{2}}\ar[dl]^{\psi^{\frac{2m+2k-1}{2}}_{\pm, m+k}}\\
&\sut{m+k}\ar[lu]&
}	
\end{equation*}
This implies that $\psi^{\frac{2m+2k-1}{2}}_{\pm, m+k}$ are both surjective. Since
$$\pi_{m,k}^{\pm}=\Psi^{m+k}_{\pm,m+2k-1}\circ\psi^{\frac{2m+2k-1}{2}}_{\mp,m+k},$$
we have
\begin{equation}\label{eq: image of pi}
	\im\pi_{m,k}^{\pm}=\im \Psi^{m+k}_{\pm,m+2k-1}.
\end{equation}

% From Corollary \ref{cor: truncated integral surgery formula}, we have the following
% $$I^{\sharp}(-Y_{-m}(K))\cong H\bigg({\rm Cone}(\sum_{j=-m-g}^{m+g}\pi_{m,k}^{+,j}+\pi_{m,k}^{-,j}:\sum_{j=-m-g}^{m+g}\sutg{\frac{2m+2k-1}{2}}{j}\to\sum_{j=-\frac{m}{2}-g}^{\frac{m}{2}+g}\sutg{m+2k-1}{j})\bigg)$$
Note that when $|j|\leq m+g$ we know that
$$\sutg{\frac{2m+2k-1}{2}}{j}\cong\mathbb{C}^{2^{2g}}\cong \sutg{m+2k-1}{j}.$$
The truncation of the integral surgery formula implies the following.
% A special case is when $m>2g-2$. In this case $$1-g+\frac{m}{2}>g-1-\frac{m}{2}$$
% so from Lemma \ref{lem: the map pi^pm_m,k} the map 
% $$\sum_{j=-m-g}^{m+g}\pi_{m,k}^{+,j}+\pi_{m,k}^{-,j}$$
% is surjective. As a result,
% $$\dim I^{\sharp}(-Y_{-m}(K))=\dim\bigg(\sum_{i=-m-g}^{m+g}\sutg{\frac{2m+2k-1}{2}}{i}\bigg)-\dim\bigg(\sum_{i=-\frac{m}{2}-g}^{\frac{m}{2}+g}\sutg{m+2k-1}{i}\bigg)=2^{2g}\cdot m.$$

% For smaller $m$, by Lemma \ref{lem: the map pi^pm_m,k}, we know that $\pi_{m,k}^{-,-m-g}=0$ and
% $$\pi_{m,k}^{+,i}:\sutg{\frac{2m+2k-1}{2}}{-m-g}\xra{\cong}\sutg{m+2k-1}{-\frac{m}{2}-g}.$$
% So we can remove $\sutg{\frac{2m+2k-1}{2}}{-m-g}$ and $\sutg{m+2k-1}{-\frac{m}{2}-g}$ from the integral surgery formula by Lemma \ref{lem: reduce the chain complex}. We can repeat this argument using Lemma \ref{lem: the map pi^pm_m,k} and Lemma \ref{lem: reduce the chain complex} to obtain the following.
\begin{equation}\label{eq: truncated mapping cone}
	I^{\sharp}(-Y_{-m}(K))\cong H\bigg({\rm Cone}(\pi^T_{m,k}:\sum_{j=1-g}^{g-1}\sutg{\frac{2m+2k-1}{2}}{j}\to\sum_{j=\frac{m}{2}+1-g}^{g-1-\frac{m}{2}}\sutg{m+2k-1}{j})\bigg)
\end{equation}
where the map
$$\pi^T_{m,k}=\sum_{i=1-g}^{g-1-m}\pi_{m,k}^{+,i}+\sum_{i=1-g+m}^{g-1}\pi_{m,k}^{-,i}.$$
Now we compute the image of the map $\pi^T$. We discuss two different cases.

{\bf Claim 1}. If $m=2l-1$ where $1\leq l\leq g-1$, we have
$$\im\pi_{m,k}^T=\bigoplus_{j=1}^{g-l}(\mathcal{M}^+_{2g,j}\oplus \mathcal{M}^-_{2g,j}),$$
where $\mathcal{M}^{\pm}_{2g,j}\cong\mathcal{M}_{2g,j}$ is defined as in (\ref{eq: M_2g,k}) and 
$$\mathcal{M}^{\pm}_{2g,j}\subset\sutg{m+2k-1}{\pm(g-l+\frac{1}{2}-j)}.$$ 

\bpf[Proof of Claim 1.] For any $\frac{m}{2}+1-g\leq j\leq g-1-\frac{m}{2}$, we have
$$\im\pi_{m,k}^T\cap\sutg{m+2k-1}{j}=\im\pi_{m,k}^{+,j-\frac{m}{2}}\cup\im\pi_{m,k}^{-,j+\frac{m}{2}}.$$
When $\pm j\leq 0$, we have 
\beq
\im\pi_{m,k}^{\mp,j+\frac{m}{2}}&\subset \im\pi_{m,k}^{\pm,j-\frac{m}{2}}\\
(\ref{eq: image of pi})&=\Psi^{n+k}_{\pm,m+2k-1}\bigg(\sutg{n+k}{j\pm\frac{k-1}{2}}\bigg)\\
&=\mathcal{M}_{-|j|+g-\frac{m}{2}}.
\eeq
As a result, we conclude that 
$$\im\pi_{m,k}^T\subset\bigoplus_{j=1}^{g-l}(\mathcal{M}^+_{2g,j}\oplus \mathcal{M}^-_{2g,j}),$$

To show that this inclusion is an equality, assume that 
$$b\in \bigoplus_{j=1}^{g-l}(\mathcal{M}^+_{2g,j}\oplus \mathcal{M}^-_{2g,j}).$$
Without loss of generality, we can assume that $b\in \mathcal{M}^+_{g+\frac{m}{2}-j}\subset \sutg{m+2k-1}{j}$ for some $j\geq 0$. We will prove that $b\in\im\pi_{m,k}^T$. By the argument above, there exists
$$a_{j+\frac{m}{2}}\in \sutg{\frac{2m+2k-1}{2}}{j+\frac{m}{2}}$$
such that
$$\pi_{m,k}^{-,j+\frac{m}{2}}(a_{j+\frac{m}{2}})=b.$$ 
Note
$$\im\pi_{m,k}^{+,j+\frac{m}{2}}\subset\im\pi_{m,k}^{-,j+\frac{3m}{2}},$$
so we can pick
$$a_{j+\frac{3m}{2}}\in \sutg{\frac{2m+2k-1}{2}}{j+\frac{3m}{2}}$$
such that
$$\pi_{m,k}^{-,j+\frac{3m}{2}}(a_{j+\frac{3m}{2}})=-\pi_{m,k}^{+,j+\frac{m}{2}}(a_{j+\frac{m}{2}}).$$
We can repeat this argument inductively to obtain an element
$$a=a_{j+\frac{m}{2}}+a_{j+\frac{3m}{2}}+a_{j+\frac{5m}{2}}+\dots$$
such that
$$\pi_{m,k}^T(a)=b.$$
\epf

From Claim 1 and (\ref{eq: truncated mapping cone}), we can compute the dimension of $I^{\sharp}(-Y_{-m}(K))$ as:
\beq
\dim I^{\sharp}(-Y_{-m}(K))=&\sum_{j=1-g}^{g-1}\dim\sutg{\frac{2m+2k-1}{2}}{j}+\sum_{j=\frac{m}{2}+1-g}^{g-1-\frac{m}{2}}\dim\sutg{m+2k-1}{j})\\
&-2\cdot \dim \im\pi_{m,k}^T.\\
=&2^{2g}\cdot m+4\cdot\sum_{j=1}^{g-l}\sum_{i=0}^{j-1}\binom{2g}{i}.
\eeq

{\bf Claim 2}. When $m=2l$ for $1\leq l\leq g-1$, we have
$$\im\pi_{m,k}^T=\mathcal{M}_{2g,g-l}\oplus\bigoplus_{j=1}^{g-l-1}(\mathcal{M}^+_{2g,j}\oplus \mathcal{M}^-_{2g,j}),$$
where $\mathcal{M}^{\pm}_{2g,j}\cong\mathcal{M}_{2g,j}$ is defined as in (\ref{eq: M_2g,k}), 
$$\mathcal{M}^{\pm}_{2g,j}\subset\sutg{m+2k-1}{\pm(g-l+\frac{1}{2}-j)},~{\rm and}~\mathcal{M}_{2g,g-l}\subset\sutg{m+2k-1}{0}.$$

The proof of Claim 2 is similar to that of Claim 1. As a result we can compute
$$\dim I^{\sharp}(-Y_{-m}(K))=2^{2g}\cdot m+4\cdot\sum_{j=1}^{g-l-1}\sum_{i=0}^{j-1}\binom{2g}{i}+2\cdot \sum_{i=0}^{g-l-1}\binom{2g}{i}.$$
\epf
\subsection{Seifert fibered manifolds}

In this subsection, we use the generalized rational surgery in \cite[Section 10.2]{Ozsvath2011rational} to obtain the Seifert fibered manifolds by surgeries and then compute the framed instanton homology.

Following the notation in Section \ref{A rational surgery formula} and Section \ref{subsec: connected sum}, we denote the connected sum of $g$ copies of the Borromean knot by $K^g\subset Y^g=\#^{2g}S^1\times S_2$ and denote the core knot in $L(v,-r)$ by $O_{v/r}$. Let $K_\#\subset Y_\#$ be the connected sum of $K_0\deq K^g$ and $K_1\deq O_{v_1/r_1},\dots,K_n\deq O_{v_n/r_n}$. Then from \cite[Section 10.2]{Ozsvath2011rational}, the $m$-surgery on $K_\#$ gives the Seifert fibered space over a genus $g$ base orbifold with Seifert invariants $(m,r_1/v_1,\dots,r_n/v_n)$.

Similar to the calculation in Lemma \ref{lem: homology cal}, we have \begin{equation*}
    H_1(Y_\#\backslash N(K_\#))\cong \bigg(H_1(Y^g)\oplus \langle g_0,g_1,\dots,g_n\rangle\bigg)/(g_0= v_i\cdot g_i \text{ for }i=1,\dots,n)
\end{equation*}where $g_0$ is the meridian of $K^g$, $g_i$ is the generator of $$H_1(L(v_i,-r_i)\backslash N(O_{v_i/r_i}))\cong\mathbb{Z},$$and the meridian of $O_{v_i/r_i}$ is $v_i\cdot g_i$. Suppose \begin{equation}\label{eq: gcd}
    \gcd(v_i,v_j)=1 \text{ for }i\neq j\in \{1,\dots,n\}.
\end{equation}Let $$v=\prod_{i=1}^nv_i\aand v_j^\p=\frac{v}{v_j}=\prod_{i\neq j}v_i.$$ Suppose $g_i=v_i^\p\cdot g_i^\p$ for $i=1,\dots,n$ and $g_0=v\cdot g_0^\p$. Then we have\begin{equation}\label{eq: homology cal2}
    \begin{aligned}
    H_1(Y^\p\backslash N(K^\p))&\cong \bigg(H_1(Y^g)\bigoplus \mathbb{Z}\langle g_0^\p,g_1^\p,\cdots,g_n^\p\rangle\bigg)/(g_0^\p=g_i^\p\text{ for }i=1,\dots,n)\\&\cong H_1(Y^g)\oplus\mathbb{Z}\langle g_0^\p\rangle.
    \end{aligned}
\end{equation}

For $\bullet\in\{0,\dots,n,\#\}$, let $\sut{\mu}^\bullet,\sut{n}^\bullet,\psi_{\pm,n+1}^{n,\bullet},\psi_{\pm,\mu}^{n,\bullet},\psi_{\pm,n}^{\mu,\bullet}$, and $F_n^\bullet$ denote the sutured instanton homologies, the bypass maps, and the cobordism maps in Lemma \ref{lem: surgery triangles} for $K_\bullet$. Similar to Proposition \ref{prop: id of bent complex}, we have the following identifications of bent complexes.

\bprop\label{prop: id of bent complex 2}
Suppose Equation (\ref{eq: gcd}) holds. For any grading $s$, there is an identification $$A(K_\#,s)=A(K_0,s_0)\aand B^\pm(K_\#,s)=B^\pm(K_0,s_0),$$where $s_0$ is the unique grading satisfying \begin{equation}\label{eq: s condition}
    s=s_0 v+\sum_{i=1}^n s_iv_i^\p\aand |s_i|\le \frac{v_i-1}{2} \text{ for }i=1,\dots,n.
\end{equation}Moreover, we have the following commutative diagrams
\begin{equation*}
\xymatrix@R=3ex{
A(K_\#,s)\ar[rr]^{\pi^\pm(K_\#,s)}\ar[dd]_{=}&&B^\pm(K_\#,s)\ar[dd]^{=}\\
&&\\
A(K_0,s_0)\ar[rr]^{\pi^\pm(K_0,s_0)}&&B^\pm(K_0,s_0)
}
\end{equation*}
\eprop
\bpf
From Lemma \ref{lem: dim}, since $K_i$ for $i=1,\dots,n$ are core knots in lens spaces, we know $(\sut{\mu}^i,s_i)$ is nontrivial only for $|s_i|\le \frac{v_i-1}{2}$, for which the grading summand is $1$-dimensional. We can apply the graded version of (\ref{eq: connected sum}) in \cite[Proposition 5.15]{LY2021large} to show that \begin{equation}\label{eq: graded 2}
    (\sut{\mu}^\#,s)\cong\bigoplus_{s_0v+\sum_{i=1}^n s_iv_i^\p=s}\bigotimes _{i=0}^n(\sut{\mu}^i,s_i),
\end{equation}where the direct sum is again from the homology calculation in (\ref{eq: homology cal2}). 

For any fixed $s$, suppose there are integers $(s_0,\dots,s_n)$ and $(s_0^\p,\dots,s_n^\p)$ satisfying $$s=s_0v+\sum_{i=1}^n s_iv_i^\p=s_0^\p v+\sum_{i=1}^n s_i^\p v_i^\p.$$
Then we have$$(s_0-s_0^\p)v+\sum_{i=1}^n (s_i-s_i^\p)v_i^\p=0.$$For any $j\in\{1,\dots,n\}$, we have $v$ and $v_i^\p$ are divisible by $v_j$ for $i\neq j$ and $v'_j$ is not divisible by $v_j$. Hence we must have $s_j-s_j^\p$ to be divisible by $v_j$. If $$|s_j|,|s_j^\p|\le \frac{v_i-1}{2},$$then we must have $s_j=s_j^\p$. As a conclusion, for fixed $s$, there is a unique $s_0$ satisfying (\ref{eq: s condition}). From (\ref{eq: graded 2}), we have $$(\sut{\mu}^\#,s)\cong (\sut{\mu}^0,s_0).$$The remainder of the proof is similar to that of Proposition \ref{prop: id of bent complex}. Indeed, there are no differentials for $K_\#$ and $K_0$, so we do not need to identify differentials in bent complexes.

% Moreover, we have $$
% \bigoplus_{k\ge 0}(\sut{\mu}^\#,s+kv)\cong \bigoplus_{k\ge 0}(\sut{\mu}^0,s^\p+k).$$Note that two sides of the isomorphism are underlying spaces of subcomplexes of $B^+(K_\#)$ and $B^+(K_0)$ since the orders of $K_\#$ and $K_0$ are $v$ and $1$, respectively. From Proposition \ref{prop: Kunneth diff}, the differentials $d_{+}$ on both sides are the same under the isomorphism. Similarly, we have$$\bigoplus_{k\le 0}(\sut{\mu}^\#,s+kv)\cong \bigoplus_{k\le 0}(\sut{\mu}^1,s^\p+k),$$and the differentials $d_-$ on both sides are the same. Hence we conclude the identification about the bent complex (\textit{cf.} Definition \ref{defn: complex A}). The commutative diagrams follow immediately.
\epf

Even though we obtain the identification of the bent complexes as in Proposition \ref{prop: id of bent complex 2}, we still need to use the $\Lambda^*H_1(-Y^g;\mathbb{C})$-action studied in the previous two subsections to identify $B^\pm(K_\#)$. Iterating the construction of $C_{\pm,n+k}^{n,k}$ in the proof of Lemma \ref{lem: contact map}, we can construct maps \begin{equation}\label{eq: iterated C}
    C_{\pm,k_0+k_1+\cdots+k_n}^{k_0,k_1,\dots,k_n}:\sut{k_0}^0\ot \sut{k_1}^1\ot\cdots\ot \sut{k_n}^n\to \sut{k_0+k_1+\cdots+k_n}^\#.
\end{equation}
Moreover, we have the commutative diagram
\begin{equation}\label{eq: commutative diagram b iterated}
\xymatrix@R=3ex{
\sut{k_0}^0\ot \sut{k_1}^1\ot\cdots\ot \sut{k_n}^n\ar[rr]^{\psi_{\pm,k_0+1}^{k_0,0}\ot\id}\ar[dd]_{C_{\pm,k_0+k_1+\cdots+k_n}^{k_0,k_1,\dots,k_n}}&&\sut{k_0+1}^0\ot \sut{k_1}^1\ot\cdots\ot \sut{k_n}^n\ar[dd]_{C_{\pm,k_0+k_1+\cdots+k_n}^{k_0+1,k_1,\dots,k_n}}\ar[rr]^{F_{k_0+1}^0\ot \cdots\ot F_{k_n}^n}&&\bigotimes_{i=0}^n I^\sharp(-Y_i)\ar[dd]^{=}\\
&&\\
\sut{k_0+k_1+\cdots+k_n}^\#\ar[rr]^{\psi_{\pm,k_0+1+k_1+\cdots+k_n}^{k_0+k_1+\cdots+k_n,\#}}&&\sut{k_0+1+k_1+\cdots+k_n}^\#\ar[rr]^{F_{k_0+1+k_1+\cdots+k_n}^\#}&&I^\sharp(-Y_\#)
}
\end{equation}
Since the construction of the map in (\ref{eq: iterated C}) only involves the neighborhoods of the knots, the map commutes with the $\Lambda^*H_1(-Y^g;\mathbb{C})$-action and we regard it as a map between the $\Lambda^*H_1(-Y^g;\mathbb{C})$ modules.

Similar to the computation in Corollary \ref{cor: grading shift of C} (\textit{cf.} (\ref{eq: graded 2})), we have\begin{equation}\label{eq: C grading shift}
    \begin{aligned}
    C_{\pm,k_0+k_1+\cdots+k_n}^{k_0,k_1,\dots,k_n}&\bigg((\sut{k_0}^0,s_0\pm \frac{k_0-1}{2})\ot\bigotimes_{i=1}^n(\sut{k_i}^i,s_i\pm \frac{(k_i-1)v_i+r_i}{2})\bigg)\\&\subset (\sut{k_0+k_1+\cdots+k_n}^\#,s_0v+\sum_{i=1}^ns_iv_i^\p\pm \frac{(\sum_{i=0}^nk_i-1)v+\sum_{i=1}^nv_i^\p r_i}{2}).
    \end{aligned}
\end{equation}
Note that the last sum $\sum_{i=1}^nv_i^\p r_i$ comes from the fact that the homology class of the longitude of $K_\#$ is the sum of the homology classes $r_i\cdot g_i$ of the longitudes of $K_i$ for $i=1,\dots,n$ under the isomorphism (\ref{eq: homology cal2}).
\bprop[]\label{prop: module inclusion}
Suppose Equation (\ref{eq: gcd}) holds. For any large enough integer $l$ and any grading $i$, the summand $(\sut{k}^\#,i)$ is a module $\mathcal{M}_{g,l}$ over $\Lambda^*H_1(-Y^g;\mathbb{C})$ for some $l$, as constructed in (\ref{eq: M_2g,k}). Moreover, the bypass maps $$\psi_{\pm,k+1}^{k,\#}:(\sut{k}^\#,i)\to (\sut{k+1}^\#,i\mp \frac{v}{2})$$ are either an inclusion of modules or the identity.
\eprop
\bpf
From \cite[Proposition 3.15]{LY2022integral1}, if $s>g+\frac{p-1}{2}-kp$, then we have an isomorphism $$H(B^+ (\ge s))\cong (\sut{k},s+\frac{(k-1)p-q}{2})$$and if $s>-(g+\frac{p-1}{2}-kp)$, then we have an isomorphism $$H(B^- (\le s))\cong (\sut{k},s-\frac{(k-1)p-q}{2}),$$where $g=g(K_\#)=g(K_0)$, $(p,q)$ are defined as in Subsection \ref{sec: SHI, 1}, and $B^+(\ge s),B^-(\le s)$ are subcomplexes of $B^\pm(s)$. Hence for $k$ large enough, we can compute $(\sut{k}^\#,i)$ by $B^\pm(K_\#)$ and Lemma \ref{lem: structure of Gamma_n}. Note that $(p,q)=(v,-\sum_{i=1}^nv_i^\p r_i)$ for $K_\#$ and $(p,q)=(1,0)$ for $K_0$. We only show the computation for the large grading $i$ and the positive bypass map as follows.

From the identification of $B^\pm$ in Proposition \ref{prop: id of bent complex 2}, we know $$(\sut{k}^\#,s+\frac{(k-1)v+\sum_{i=1}^nv_i^\p r_i}{2})\cong H(B^+ (K_\#,\ge s))\cong H(B^+(K_0,\ge s_0))\cong (\sut{k_0}^0,s_0+\frac{k_0-1}{2}),$$where $s=s_0 v+\sum_{i=1}^ns_iv_i^\p$ with $|s_i|\le \frac{v_i-1}{2}$ and $s,s_0$ satisfy the inequality $s,s_0>g+\frac{p-1}{2}-kp$ for their corresponding $(p,q)$. Let $k=\sum_{i=0}^n k_i$. From Lemma \ref{lem: dim}, we know $$(\sut{k_i}^i,s_i+\frac{(k_i-1)v_i+r_i}{2})\cong\mathbb{C}.$$From (\ref{eq: C grading shift}), the map $C_{+,k_0+k_1+\cdots+k_n}^{k_0,k_1,\dots,k_n}$ induces a map from $(\sut{k_0}^0,s_0+\frac{k_0-1}{2})$ to $(\sut{k}^\#,s+\frac{(k-1)v+\sum_{i=1}^nv_i^\p r_i}{2})$. Since there are no differentials for $K_0$ and $K_\#$, the triangles involving $\psi_{+,k_0+1}^{k_0,0}$ and $\psi_{+,k_0+1+k_1+\cdots+k_n}^{k_0+k_1+\cdots+k_n,\#}$ split and these two maps are injective. From Lemma \ref{lem: structure of Gamma_n}, after applying the bypass maps for sufficiently many times $k^\p$, the restrictions of maps $F_{k_0+k^\p}^0$ and $F_{k_0+k^\p+k_1+\cdots+k_n}^\#$ are isomorphisms. Then the commutativity in (\ref{eq: commutative diagram b iterated}) implies that $$C_{+,k_0+k_1+\cdots+k_n}^{k_0,k_1,\dots,k_n}:(\sut{k_0}^0,s_0+\frac{k_0-1}{2})\to (\sut{k}^\#,s+\frac{(k-1)v+\sum_{i=1}^nv_i^\p r_i}{2})$$ is an isomorphism. Thus, the proposition follows from the computation in Lemma \ref{lem: gamma_n for connected sums of Borromean knot}, Corollary \ref{cor: bypass maps for the connected sum of Borromean knots}, and the commutativity in (\ref{eq: commutative diagram b iterated}).
\epf

\bpf[Proof of Theorem \ref{thm: SF space}]
Similar to the proof of Theorem \ref{thm: circle bundle computation} in the last subsection, we apply the integral surgery formula Theorem \ref{thm: integral surgery formula} and its truncation Proposition \ref{prop: truncated integral surgery formula} to the connected sum $K_\#$ of Borromean knots and core knots in lens spaces. Since there are no differentials for $K_\#$, again $\psi_{\pm,m+k}^{\frac{2m+2k+1}{2}}$ is surjective and \begin{equation*}
	\im\pi_{m,k}^{\pm}=\im \Psi^{m+k}_{\pm,m+2k-1}.
\end{equation*}as in (\ref{eq: image of pi}). Then the dimension of $\im \pi_{m,k}^T$ in the truncation can be computed from Proposition \ref{prop: module inclusion} and the two claims in the proof of Theorem \ref{thm: circle bundle computation}.

In Heegaard Floer theory, we apply Ozsv\'ath-Szab\'o's integral surgery formula for $\widehat{HF}$. There is an explicit identification between the Heegaard Floer version of $B^\pm(K_0)$ in \cite[Lemma 10.4]{Ozsvath2011rational} (the lemma is for the plus version, but setting $U=0$ gives the identification for the hat version), which coincides with the identification from the $\Lambda^*H_1(-Y_g;\mathbb{C})$-module structure.

Since the integral surgery formulae in instanton and Heegaard Floer theories have a similar form and we have already shown that the complexes and the maps in the formulae coincide, the dimensions of $I^\sharp$ and $\widehat{HF}$ for the surgery manifold are the same. Note that we have to use the dimension over $\ft$ for $\widehat{HF}$ since \cite[Lemma 10.4]{Ozsvath2011rational} works over $\ft$. The computations by two claims in the proof of Theorem \ref{thm: circle bundle computation} are independent of the underlying field.
\epf

\section{Surgeries on some alternating knots}\label{sec: Surgeries on some alternating knots}
In this section, we use oriented skein relation and an inductive argument to study differentials for a special family of alternating knots in $S^3$.

\subsection{Knots with torsion order one}
In this subsection, we introduce a condition on the differentials that is closely related to the thin complex in \cite[Definition 6]{Petkova2009thin}. Inspired by the $U$ map in Definition \ref{defn: khi minus}, we have the following definition.
\bdefn\label{defn: the U map}
Suppose $K\subset Y$ is a rationally null-homologous knot of order $p$. For a large enough integer $n$, and a grading $i$, we define the map
\begin{equation}\label{eq: U map}
    U=(\psm{n}{n+1})^{-1}\circ\psp{n}{n+1}:\sutg{n}{i}\to\sutg{n}{i-p}.
\end{equation}
\edefn

\blem\label{lem: the U map}The following are some basic properties of the map $U$.
\begin{enumerate}
	\item The map $U$ is well-defined for any $i\geq g- \frac{(n-1)p-q-1}{2}$.
	\item For any $i$ such that $U$ is defined, there exists an exact triangle
	\begin{equation*}
	\xymatrix{
	\sutg{n}{i}\ar[rr]^{U}&&\sutg{n}{i-p}\ar[dl]^{\psp{n}{\mu}}\\
	&\sutg{\mu}{i-\frac{(n-1)p-q}{2}}\ar[ul]^{\psp{\mu}{n}}&
	}
\end{equation*}
	\item We have $F_n\circ U=F_n$.
\end{enumerate}
\elem
\bpf
Part (1) follows from Lemma \ref{lem: structure of Gamma_n} part (2). Part (2) follows from Lemma \ref{lem: bypass n,n+1,mu} and Lemma \ref{lem: comm diag for n,n+1,mu}. Part (3) follows from Lemma \ref{lem: comm diag for n,n+1,dehn} and Lemma \ref{lem: F_n and G_n are iso when n large} Part (1). 
\epf

From Lemma \ref{lem: the U map} part (1), the map $U$ is well-defined on most of the gradings of $\sut{n}$. Since $n$ is large, it is enough to focus on $i>0$. 

From diagram (\ref{eq: useful triangles}), the differentials $d_\pm$ induce \begin{equation}\label{eq: d1 diff}
    d_{1,\pm}\deq \psi_{\pm,\mu}^n\circ \psi_{\pm,n}^{\mu}
\end{equation}on the first pages of corresponding spectral sequences. The definition is independent of the choice of $n$ due to Lemma \ref{lem: comm diag for n,n+1,mu}. 

\begin{lem}\label{lem: torsion order 1}
	Suppose $K\subset Y$ is a rationally null-homologous knot. The following are equivalent.
	\begin{itemize}
		\item [(i)] $\dim H(\sut{\mu},d_{1,+})=\dim I^{\sharp}(-Y)$.
		\item [(ii)] $\dim H(\sut{\mu},d_{1,-})=\dim I^{\sharp}(-Y)$.
		\item [(iii)] For large enough $n$ and any element $x\in\sutg{n}{i}$ with $i>0$, if there exists $k\in\mathbb{N}_+$ such that $U^k(x)=0$, then $U(x)=0$.
		\item [(iv)] For large enough $n$ and any grading $i>0$, we have
		$$U\bigg(\sutg{n}{i}\cap\ker F_n\bigg)=0.$$
	\end{itemize}
\end{lem}

\bpf
If we reverse the orientation of the knot, then positive and negative bypasses in defining the differentials $d_{1,\pm}$ in Equation (\ref{eq: d1 diff}) exchange with each other. As a result, the two differentials $d_{1,+}$ and $d_{1,-}$ also switch with each other. Hence we conclude that (i) and (ii) are equivalent. The equivalence between (iii) and (iv) follows easily from Lemma \ref{lem: the U map} Part (3). To show that (i) and (iii) are equivalent, recall that in \cite{LY2021large} the construction of the differentials $d_{+}$ involves a series of differentials $d_{k,+}$ defined as
$$d_{k,+}=\psp{n}{\mu}\circ(\Psp{n}{n+k})^{-1}\circ\psp{\mu}{n}.$$In \cite{LY2021large}, we proved that $d_{k,+}$ is well-defined on $\ker d_{k-1,+}\slash\im d_{k-1,+}$ and
$$\ker d_{k,+}\slash\im d_{k,+}\cong I^{\sharp}(-Y)$$
for any large enough $k$. For simplicity, we suppose $n$ is large enough. Since $\im\psp{\mu}{n}$ lies in the top few gradings of $\sut{n}$, by Definition \ref{defn: the U map} the map $U$ is well-defined on related gradings. Also from Lemma \ref{lem: structure of Gamma_n} we know that $\psm{n}{n+1}$ is an isomorphism on such gradings, so we can rewrite $d_{k,+}$ as
$$d_{k,+}=\psp{n}{\mu}\circ U^{-k}\circ \psp{\mu}{n}.$$
Now statement (i) is equivalent to the fact that $d_{k,+}=0$ for all $k\geq 2$ and it remains to show that this is equivalent to (iii).

If there exists $u\in\sut{\mu}$ and $k\geq 2$ such that $d_{i,+}(u)=0$ for all $i<k$ and $d_{k,+}(u)\neq0$. Then by definition there exists $y\in\sut{n}$ such that $\psp{n}{\mu}(y)\neq0$ and $U^{k-1}(y)=\psp{\mu}{n}(x)\neq 0$. Since by Lemma \ref{lem: the U map}
$$U^{k}(y)=U\circ\psp{\mu}{n}(x)=0$$
we know that (iii) does not hold. Conversely, if there exists $x\in\sut{n}$ such that $x\notin\im U$, $U^k(x)\neq0$ and $U^{k+1}(x)=0$ for some $k\geq 1$. We know that there exists $u\in\sut{\mu}$ such that $\psp{\mu}{n}(u)=U^k(x)$ and hence
$$d_{k+1,+}(u)=\psp{n}{\mu}(x)\neq0.$$
Hence we conclude that (i) and (iii) are equivalent.
%check
\epf

\bdefn\label{defn: torsion order one}
A knot $K\subset Y$ has {\bf torsion order one} if it satisfies any equivalent statement in Lemma \ref{lem: torsion order 1}.
\edefn

\subsection{Commutativity of the first differentials}\label{subsec: commutativity for d_1}

In this subsection, we prove the commutativity of two first differentials, which will provide a strong restriction for knots with torsion order one.
\begin{thm}\label{thm: main equality}
	For any rationally null-homologous $K\subset Y$, we have
	\[d_{1,-}\circ d_{1,+}\doteq d_{1,+}\circ d_{1,-},\]
	where $\doteq$ means the equation holds up to a scalar.
\end{thm}

From \cite{li2018gluing}, there is a gluing map
\[G:\sut{\mu}\otimes \shi(-[0,1]\times T^2,-\Ga_{\mu}\cup-\Ga_{\mu})\ra \sut{\mu}.\]
Here we can identify $\{0\}\times T^2$ with $\partial (S^3\backslash N(K))$ and then use the Seifert framing on $\partial (S^3\backslash N(K))$ to be the framing on $T^2$ as well. Let $\xi_{st}$ be the product contact structure on $[0,1]\times T^2$, and \[\theta(\xi_{st})\in \shi(-[0,1]\times T^2,-\Ga_{\mu}\cup-\Ga_{\mu})\] be its contact element \cite{baldwin2016instanton}. Then we know from \cite[Theorem 1.1]{li2018gluing} that
\begin{equation}\label{eq: G leads to the identity}
G(-\otimes \theta(\xi_{st}))\doteq \id:\sut{\mu}\ra \sut{\mu}.
\end{equation}
We write $Y_{T^2}=[0,1]\times T^2$ and take $n=0$ in the definition of $d_{1,\pm}$ in (\ref{eq: d1 diff}) for simplicity. We can view the bypasses, which are originally attached to $(S^3\backslash N(K),\Ga_{\mu})$, to be attached to $(Y_{T^2},\Ga_{\mu}\cup\Ga_{\mu})$ on the $\{1\}\times T^2$ side instead, and they lead to new exact triangles
\begin{equation}\label{eq: bypasses on product}
\xymatrix{
\shi(-Y_{T^2},-\Ga_{\mu}\cup-\Ga_{0})\ar[rr]^{\hat{\psi}_{\pm,1}^{0}}&&\shi(-Y_{T^2},-\Ga_{\mu}\cup-\Ga_{1})\ar[dl]^{\hat{\psi}_{\pm,\mu}^{1}}\\
&\shi(-Y_{T^2},-\Ga_{\mu}\cup-\Ga_{\mu})\ar[ul]^{\hat{\psi}_{\pm,0}^{\mu}}&
}	
\end{equation}

Using these bypass maps, we could construct the map $\hat{d}_{\pm}=\hat{\psi}_{\pm,\mu}^0\circ\hat{\psi}_{\pm,0}^{\mu}$ in the same way as the construction of the maps $d_{1,\pm}$. We have the following key proposition:
\begin{prop}\label{prop: key identity on the product}
	We have
	\[\hat{d}_+\circ\hat{d}_-(\theta(\xi_{st}))\doteq \hat{d}_-\circ\hat{d}_+(\theta(\xi_{st}))\neq 0\in \shi(-[0,1]\times T^2,-\Ga_{\mu}\cup-\Ga_{\mu})\]
\end{prop}

\bpf[Proof of Theorem \ref{thm: main equality} using Proposition \ref{prop: key identity on the product}.]
From the functoriality of gluing maps in \cite[Theorem 1.1]{li2018gluing}, we have two commutative diagrams
\begin{equation*}
\xymatrix{
\sut{\mu}\otimes \shi(-Y_{T^2},-\Ga_{\mu}\cup-\Ga_{\mu})\ar[rr]^{\id\otimes(\hat{d}_+\circ\hat{d}_-)}\ar[d]^{G}&&\sut{\mu}\otimes \shi(-Y_{T^2},-\Ga_{\mu}\cup-\Ga_{\mu})\ar[d]^{G}\\
\sut{\mu}\ar[rr]^{d_{1,+}\circ d_{1,-}}&&\sut{\mu}
}	
\end{equation*}
and
\begin{equation*}
\xymatrix{
\sut{\mu}\otimes \shi(-Y_{T^2},-\Ga_{\mu}\cup-\Ga_{\mu})\ar[rr]^{\id\otimes(\hat{d}_-\circ\hat{d}_+)}\ar[d]^{G}&&\sut{\mu}\otimes \shi(-Y_{T^2},-\Ga_{\mu}\cup-\Ga_{\mu})\ar[d]^{G}\\
\sut{\mu}\ar[rr]^{d_{1,-}\circ d_{1,+}}&&\sut{\mu}
}	
\end{equation*}
From (\ref{eq: G leads to the identity}), we have  $G(x\otimes \theta(\xi_{st}))\doteq x$. Hence, from the commutative diagrams and Proposition \ref{prop: key identity on the product}, we have\begin{equation}\label{eq: composition of differentials}
d_{1,+}\circ d_{1,-}(x)\doteq G(x\otimes \hat{d}_+\circ\hat{d}_-(\theta(\xi_{st})))\doteq G(x\otimes \hat{d}_-\circ\hat{d}_+(\theta(\xi_{st})))=d_{1,-}\circ d_{1,+}(x).
\end{equation}
\epf
Then we prove Proposition \ref{prop: key identity on the product}. First note that $\hat{d}_+$ and $\hat{d}_-$ are both constructed by bypasses, and contact elements are preserved by the gluing maps as in \cite[Theorem 1.1]{li2018gluing}. As a result, there are two contact structures $\xi_{+-}$ and $\xi_{-+}$ on $(Y_{T^2},\Ga_{\mu}\cup\Ga_{\mu})$, which are both obtained from $\xi_{st}$ by attaching four bypasses according to $d_{1,+}\circ d_{1,-}$ and $d_{1,-}\circ d_{1,+}$, respectively, so that
\[\theta(\xi_{+-})\doteq\hat{d}_+\circ\hat{d}_-(\theta(\xi_{st}))~{\rm and}~\theta(\xi_{-+})\doteq\hat{d}_-\circ\hat{d}_+(\theta(\xi_{st})).\]

\blem\label{lem: nonvanishing}
The contact elements $\theta(\xi_{+-})$ and $\theta(\xi_{-+})$ are both nonzero.
\elem
\bpf
From (\ref{eq: composition of differentials}) we know that for any knot $K\subset S^3$, we have
\[d_{1,+}\circ d_{1,-}\doteq G(-\otimes \theta(\xi_{+-}))~{\rm and~}d_{1,-}\circ d_{1,+}\doteq G(-\otimes \theta(\xi_{+-})).\]
We computed the differentials for the figure-eight knot in \cite[Section 6]{LY2021large}, for which we have $d_{1,+}\circ d_{1,-}\neq0$ and $d_{1,-}\circ d_{1,+}\neq 0$. Thus the lemma follows.
\epf
Next, to better study the two contact elements, we construct a $\intg^2$-grading on $\shi(-Y_{T^2},-\Ga_{\mu}\cup-\Ga_{\mu})$ as follows. View $T^2=S^1\times S^1$. We call curves that are isotopic to $S^1\times\{\operatorname{pt}\}$ and $\{\operatorname{pt}\}\times S^1$ longitudes and meridians, respectively. Taking a meridian $m$ on $T^2$, we have an annulus $A_m=[0,1]\times m\subset Y_{T^2}$. We can arrange $A_m$ as a product annulus inside $(Y_{T^2},\Ga_{\mu}\cup\Ga_{\mu})$. The decomposition along $A_m$ yields a solid torus with the suture being six copies of the longitude. According to \cite[Lemma 2.29]{li2019decomposition}, we have
\[\shi(-Y_{T^2},-\Ga_{\mu}\cup-\Ga_{\mu})\cong\mathbb{C}^4.\]
As in \cite[Theorem 2.28]{li2019decomposition}, the surface $A_m$ induces a $\intg$-grading on $\shi(-Y_{T^2},-\Ga_{\mu}\cup-\Ga_{\mu})$ such that all six dimensions are supported at grading $0$.

For a second surface, we pick a longitude of $l$ of $T^2$ and obtain a second annulus $A_l=[0,1]\times l$. Note that each component of $\partial A_l$ intersects the suture $\Ga_{\mu}$ twice, so as in \cite[Theorem 2.28]{li2019decomposition}, $A_l$ induces a $\intg$-grading on $\shi(-Y_{T^2},-\Ga_{\mu}\cup-\Ga_{\mu})$ which is supported at three gradings $-1$, $0$, $1$. The decompositions along $A_l$ and $-A_l$ both yield a solid torus with sutures being two copies of the longitude. According to \cite[Lemma 2.29]{li2019decomposition}, we have
\[\shi(-Y_{T^2},-\Ga_{\mu}\cup-\Ga_{\mu},A_l,1)\cong \shi(-Y_{T^2},-\Ga_{\mu}\cup-\Ga_{\mu},A_l,-1)\cong\mathbb{C}.\]
As a result,
\[\shi(-Y_{T^2},-\Ga_{\mu}\cup-\Ga_{\mu},A_l,0)\cong\mathbb{C}^2.\]
From \cite[Section 5.1]{li2019decomposition}, the surfaces $A_m$ and $A_l$ together induce a $\intg^2$-grading.
\blem\label{lem: bypass maps are homogeneous}
Suppose $(M,\ga)$ is a balanced sutured manifold and $S\subset M$ is a properly embedded surface. Let $\be\subset \partial M$ be a bypass arc, and let the bypass attachment along $\be$ changes the suture $\ga$ to $\ga'$. Let
\[\psi:\shi(-M,-\ga)\ra \shi(-M,-\ga')\]
be the corresponding bypass map. Then $\psi$ is homogeneous with respect to the grading induced by $S$ on $\shi(-M,-\ga)$ and $\shi(-M,-\ga')$.
\elem
\bpf
Since $\be$ is an arc, we can always perform stabilizations on $S$ in the sense of \cite[Definition 3.1]{li2019direct} to make $S$ disjoint from $\be$. Then as in the proof of \cite[Proposition 5.5]{li2019direct}, $\psi$ is clearly homogeneous.
\epf
\blem\label{lem: extremal bi-grading}
Suppose $(M,\ga)$ is a balanced sutured manifold and $S_1$ and $S_2$ are two admissible surfaces in the sense of \cite[Definition 2.26]{li2019decomposition} in $(M,\ga)$. Let $(i,j)$ denote the $\intg^2$-grading on $\shi(M,\ga)$ induced by the pair of surfaces $(S_1,S_2)$. Let
\[i_0=\frac{1}{4}|S_1\cap\ga|-\frac{1}{2}\chi(S_1),~{\rm and~}j_0=\frac{1}{4}|S_2\cap\ga|-\frac{1}{2}\chi(S_2).\]
Suppose $(M_1,\ga_1)$ is obtained from $(M,\ga)$ by decomposing along $S_1$, and $S_2'\subset(M_1,\ga_1)$ is obtained from $S_2$ by cutting along $S_1$. Suppose $(M_2,\ga_2)$ is obtained from $(M_1,\ga_1)$ by decomposing along $S_2'$. Then we have an isomorphism
\[\shi(M,\ga,(i_0,j_0))\cong \shi(M_2,\ga_2).\]
\elem
\bpf
By \cite[Lemma 2.29]{li2019decomposition}, we have
\[\shi(M,\ga,S_1,i_0)\cong \shi(M_1,\ga_1).\]
Applying this fact again, we conclude that
\[\shi(M,\ga,(S_1,S_2),(i_0,j_0))\cong \shi(M_2,\ga_2).\]
\epf

Next, we want to study a graded version of exact triangle (\ref{eq: bypasses on product}). First, we want to figure out the double grading on $\shi(-Y_{T^2},-\Ga_{\mu}\cup-\Ga_{0})$ and $\shi(-Y_{T^2},-\Ga_{\mu}\cup-\Ga_{1})$ induced by the pair of annuli $(A_l,A_m)$. For the sutured manifold $(-Y_{T^2},-\Ga_{\mu}\cup-\Ga_{0})$, $A_l$ and $A_m$ each intersect the suture at two points, so we perform a negative stabilization in the sense of \cite[Definition 3.1]{li2019direct} on each of them to obtain two surfaces $A_l^-$ and $A_m^-$. Then the $\intg$-gradings associated to $A_l^-$ and $A_m^-$ are both supported at grading $0$ and $1$.
\begin{lem}\label{lem: Gamma_0 on product}
We have
\[\shi(-Y_{T^2},-\Ga_{\mu}\cup-\Ga_{0})\cong \mathbb{C}^4\]
and the four generators are supported at bi-gradings $(0,0)$, $(0,1)$, $(1,0)$, and $(1,1)$.
\end{lem}
\bpf
This is a direct application of Lemma \ref{lem: extremal bi-grading} by looking at the four pairs of surfaces $(-A_l,-A_m)$, $-A_l,A_m$, $(A_l,-A_m)$, and $(A_l,A_m)$. Note that when dealing with $-A_l$ and $-A_m$, we need to use a positive stabilization instead, and use the grading shifting property in \cite[Theorem 1.12]{li2019decomposition} to relate the grading induced by $A_l^{\pm}$ and $A_m^{\pm}$.
\epf

For the sutured manifold $(-Y_{T^2},-\Ga_{\mu}\cup-\Ga_{1})$, the annulus $A_l$ intersects the suture four times, so it induces a $\intg$-grading where all non-vanishing gradings are $-1,0,1$. The annulus $A_m$ intersects the suture twice, so we perform a negative stabilization as above and use $A_m^-$ to construct a $\intg$-grading on $\shi(-Y_{T^2},-\Ga_{\mu}\cup-\Ga_{1})$. The non-vanishing gradings are $0$ and $1$.
\begin{lem}\label{lem: Gamma_1 on product}
We have
\[\shi(-Y_{T^2},-\Ga_{\mu}\cup-\Ga_{1})\cong \mathbb{C}^4\]
and the four generators are supported at bi-gradings $(0,-1)$, $(0,0)$, $(1,0)$, and $(1,1)$.
\end{lem}
\bpf
This is a direct application of Lemma \ref{lem: extremal bi-grading} by looking at the four pairs of surfaces $(-A_l,-A_m)$, $(-A_l,A_m)$, $(A_l,-A_m)$, and $(A_l,A_m)$. Note that when dealing with $-A_m$, we need to use a positive stabilization instead, and use the grading shifting property in \cite[Theorem 1.12]{li2019decomposition} to relate the grading induced by $A_m^{+}$ and $A_m^{-}$.
\epf

\bpf[Proof of Proposition \ref{prop: key identity on the product}.] 

By Lemma \ref{lem: bypass maps are homogeneous}, there is a graded version of (\ref{eq: bypasses on product}) as follows.
\begin{equation}\label{eq: bypasses on product, graded}
\xymatrix{
\shi(-Y_{T^2},-\Ga_{\mu}\cup-\Ga_{0},(i_0',j_0'))\ar[r]^{\hat{\psi}_{+,1}^{0,(i_0',j_0')}}&\shi(-Y_{T^2},-\Ga_{\mu}\cup-\Ga_{1},(i_1',j_1'))\ar[dl]^{\hat{\psi}_{+,\mu}^{1,(i_1',j_1')}}\\
\shi(-Y_{T^2},-\Ga_{\mu}\cup-\Ga_{\mu},(0,0))\ar[u]^{\hat{\psi}_{+,0}^{\mu,(0,0)}}&
}	
\end{equation}
for some indices $(i_0',j_0')$ and $(i_1',j_1')$. From the above argument, we know that
\[\shi(-Y_{T^2},-\Ga_{\mu}\cup-\Ga_{\mu},(0,0))\cong \mathbb{C}^2.\]
From Lemma \ref{lem: Gamma_0 on product}, we know that
\[\dim\shi(-Y_{T^2},-\Ga_{\mu}\cup-\Ga_{0},(i_0',j_0'))\leq 1.\]
From Lemma \ref{lem: Gamma_1 on product}, we know that
\[\dim\shi(-Y_{T^2},-\Ga_{\mu}\cup-\Ga_{1},(i_1',j_1'))\leq 1.\]
Since the three terms fit into an exact triangle as in (\ref{eq: bypasses on product, graded}), we must have
\[\dim\shi(-Y_{T^2},-\Ga_{\mu}\cup-\Ga_{0},(i_0',j_0'))=1\]
and the map $\hat{\psi}_{\pm,0}^{\mu,(0,0)}$ is surjective. Since we already have the nonvanishing result in Lemma \ref{lem: nonvanishing}, to show that
\[\theta(\xi_{+-})\doteq \theta(\xi_{-+}),\]
it suffices to prove that
$$\hat{\psi}_{+,0}^{\mu}(\theta(\xi_{+-}))=\hat{\psi}_{+,0}^{\mu}(\theta(\xi_{-+}))=0.$$

For the contact structure $\xi_{-+}$, the image $\hat{\psi}_{+,0}^{\mu}(\theta(\xi_{-+}))=0$ because after attaching the last bypass to $\xi_{-+}$, the resulting contact structure admits a Giroux torsion so it has vanishing contact element by \cite[Section 4]{LY2021large}. For the contact structure $\xi_{+-}$, note that we have
\beq
\hat{\psi}_{+,0}^{\mu}(\theta(\xi_{-+}))&=\hat{\psi}_{+,0}^{\mu}\circ \hat{d}_+\circ\hat{d}_-(\theta(\xi_{st}))\\
&=\hat{\psi}_{+,0}^{\mu}\circ (\hat{\psi}_{+,\mu}^{1}\circ\hat{\psi}_{+,1}^{\mu})\circ (\hat{\psi}_{-,\mu}^1\circ\hat{\psi}_{-,1}^{\mu})(\theta(\xi_{st}))\\
&=0\\
\eeq
It is finally zero since $\hat{\psi}_{+,0}^{\mu}$ and $\hat{\psi}_{+,\mu}^{1}$ fit into an exact triangle.
\epf

\subsection{Classification of complexes}\label{subsec: Bent complexes for knots of torsion order one}
Suppose $K\subset S^3$ is a knot of torsion order one ({\it cf.} Definition \ref{defn: torsion order one}). From Lemma \ref{lem: torsion order 1}, we know $d_{\pm}=d_{1,\pm}$, \textit{i.e.}, differentials in higher pages vanish. Then Theorem \ref{thm: main equality} imposes strong restrictions on the differentials. In this subsection, we prove a classification theorem for complexes of knots of torsion order one.

\blem
Suppose $K\subset S^3$ is a knot with torsion order one and
$$\dim \sut{\mu}=||\Delta_{K}(t)||,$$
where $||\cdot||$ is the sum of absolute values of coefficients. Write $d_{\pm}=d_{1,\pm}$ for simplicity. Then, up to changing a basis, the pair $(\sut{\mu},d_++d_-)$ is the direct sum of the following three basic types of complexes, which are called squares for $C$ and staircases for $C_l$.
\[
\xymatrix{\\\\
	c\ar[d]_{\lambda\cdot d_-}&a\ar[l]_{d_+}\ar[d]^{d_-}\\
	d&b\ar[l]^{d_+}
	}
	\xymatrix{
	a_1\ar[d]_{d_-}&\\
	a_2&a_3\ar[l]^{d_+}\ar[d]_{d_-}\\
	&\cdots&a_{2|l|-1}\ar[l]^{d_+}\ar[d]_{d_-}\\
	&&a_{2|l|}&a_{2|l|+1}\ar[l]^{d_+}
	}
	\xymatrix{
	a_{2l+1}&a_{2l}\ar[l]_{d_+}\ar[d]^{d_-}\\
	&a_{2l-1}&\cdots\ar[l]_{d_+}\ar[d]^{d_-}\\
	&&a_3&a_2\ar[l]_{d_+}\ar[d]^{d_-}\\
	&&&a_1
	}\]
	\[\quad \quad\quad C\quad\quad\quad\quad\quad\quad\quad\quad\quad\quad\quad C_{l} \text{ for }l\le 0\quad\quad\quad\quad\quad\quad\quad\quad\quad\quad\quad\quad\quad\quad C_{l} \text{ for }l>0\quad\quad\quad\quad\]
	where $\lambda$ is the scalar from Theorem \ref{thm: main equality} that makes the diagram in $C$ commute.
\elem

\bpf
The proof is an adaption of the proof of \cite[Lemma 7]{Petkova2009thin} to our setup. Note that the proof in the reference studied spaces with coefficients $\mathbb{F}_2$, while we deal with coefficients $\mathbb{C}$ here. Theorem \ref{thm: main equality} shows that $d_+$ and $d_-$ commute up to a scalar, so $(d_++d_-)^2$ is not necessarily zero if the scalar is not $-1$. But we can still carry out the proof similarly.

We now treat $(\sut{\mu},d_++d_-)$ as a purely algebraic object and prove by induction on the dimension of $\sut{\mu}$. Fix a basis of $\sutg{\mu}{i}$ for each grading $i$ that is homogeneous with respect to the $\mathbb{Z}_2$ homological grading. For a basis element $b$, we say that there is an upward arrow from an element $a$ to $b$ if
$$d_+(a)=\lambda\cdot b+{\rm (linear~combination~of~other~basis~elements)}$$
for some $\lambda\neq 0$. To be consistent with the complex in \cite[Lemma 7]{Petkova2009thin}, we use leftward arrows to represent upward arrows. In particular, $w$ and $z$ arrows correspond to $d_+$ and $d_-$ arrows, respectively. Note that if $b\in\sutg{\mu}{i}$, then $a\in\sutg{\mu}{i-1}$ since $d_+=d_{1,+}$ shifts the Seifert grading by $+1$. We define downward arrows using $d_-$ similarly.

We start with a basis element $b_1\in \sutg{\mu}{-g}$, where $g=g(K)$. Note that $-g$ is the minimal nontrivial grading of $\sut{\mu}$. 

{\bf Case 1}. There is a downward arrow from $a$ to $b$ for some $a\in\sutg{\mu}{-g+1}$, \textit{i.e.}, we have $$d_-(a)=\sum_{i=1}^n\lambda_i\cdot b_i$$for some basis elements $b_2,\cdots b_n$ (possibly $n=1$). We change the basis by replacing $\{b_1,\cdots,b_n\}$ with $\{b=\sum_{i=1}^n\lambda_i\cdot b_i,b_2,\cdots,b_n\}$. Then $d_-(a)=b$. Since $b$ lives in the minimal nontrivial grading, there is no downward arrow originating at $b$ and no upward arrow pointing to $b$. If there are other basis elements with downward arrows to $b$, then we can add $a$ to each of them with proper coefficients, so that in the new basis only $a$ has a downward arrow to $b$. If there is an upward arrow from $b^\p$ to $a$ for some $b^\p\in \sutg{\mu}{-g}$, then $d_-\circ d_+ (b^\p)$ must have a nonzero coefficient on $b$, which contradicts the fact that $d_+\circ d_+(b^\p)=0$ and the commutativity from Theorem \ref{thm: main equality}. Hence there is no upward arrow pointing to $a$. 

% We can change the basis as in the proof of \cite[Lemma 7]{Petkova2009thin} so that $d_{-}(a)=b$ and $a$ is the only basis element that admits a downward arrow to $b$. If there exists $b'$ that admits an upward arrow to $a$, then we know from Theorem \ref{thm: main equality} that $d_-(b')\neq 0\in \sutg{\mu}{-g-1}$ which contradicts the minimal nontrivial grading fact. As a result, there is no upward arrow to $a$. The grading of $b$ guarantees that there is no upward arrow to $b$ either.

{\bf Case 1.1}. We have $d_+(b)\neq 0$. We will split off a $C$ summand and hence the induction applies. Indeed, let $d=d_+(b)\neq0$ and $c=d_+(a)$. From Theorem \ref{thm: main equality}, we have $$d_-(c)=d_-\circ d_+(a)\doteq d_+\circ d_-(a)=d_+(b)=d\neq 0.$$Then Case 1.1 in the proof of \cite[Lemma 7]{Petkova2009thin} applies verbatim and we can change the basis to make the following conditions hold:
\benu
\item $c$ and $d$ are basis elements;
\item $b$ is the only basis element with an upward arrow to $d$;
\item $c$ is the only basis element with a downward arrow to $d$;
\item $a$ is the only basis element with an upward arrow $c$;
\item $a$ is the only basis element with a downward arrow to $b$.
\eenu
Hence the span of $a,b,c,d$ is a $C$ summand.

{\bf Case 1.2}. We have $d_+(b)=0$. The grading of $b$ guarantees that $b\notin \im d_+$ so $$[b]\neq0\in H(\sut{\mu},d_+)\cong I^\sharp(-S^3)\cong \mathbb{C}.$$As a result, there are no other generators of $H(\sut{\mu},d_+)$. In particular, $c=d_+(a)\neq 0$ since we have already argued that there is no upward arrow to $a$. Now $d_-(c)=d_+\circ d_-(b)=0$ by a grading argument and $d_+(c)=0$ since $d^2_+=0$. As in the proof of \cite[Lemma 7]{Petkova2009thin}, we can change the basis so that $c$ is a basis element and $a$ is the only basis element with an upward arrow to $c$. 

{\bf Case 1.2.1}. There is no downward arrow to $c$. In this case we can split off the staircase spanned by $a$, $b$, and $c$.

{\bf Case 1.2.2}. There is a downward arrow to $c$. As in the proof of \cite[Lemma 7]{Petkova2009thin}, after a suitable change of basis, we either eliminate the arrow to $c$ so that we can split off a staircase spanned by $a$, $b$, and $c$, or we can find $d$ such that $d_-(d)=c$ and $d$ is the only basis element with a downward arrow to $c$, and we can repeat the argument in Case 1.2 to further trace along the staircase.

{\bf Case 2}. There is no downward arrow to $b$. We will split off from a staircase. If $d_+(b)=0$ we split off the single $b$. If $c=d_+(b)\neq 0$ we can change the basis to make $c$ a basis element and $b$ is the only basis element with an upward arrow to $c$. As above also know that $d_+(c)=0$ and $d_-(c)=0$. Note that now $[b]$ is the unique generator of $H(\sut{\mu},d_-)$ so we know that there exists $d$ with $d_-(d)=c$. As in the proof of \cite[Lemma 7]{Petkova2009thin} we can keep this argument to split off a staircase starting from $b$.

In any case we can split off either a square or a staircase hence the induction applies.
\epf

\bcor\label{cor: bent complex for torsion order one}
Suppose $K\subset S^3$ is a knot with torsion order one and
$$\dim \sut{\mu}=||\Delta_{K}(t)||.$$
Then the structure of $(\sut{\mu},d_++d_-)$ is determined by $\Delta_{K}(t)$ and $\tau_I(K)$.
\ecor
\bpf
The proof of \cite[Theorem 4]{Petkova2009thin} applies verbatim. In particular, there is a unique staircase $C_l$ with $l=\tau_I(K)$, \textit{i.e.}, $$a_1\in \sutg{\mu}{-\tau_I(K)}\aand a_{2|l|+1}\in \sutg{\mu}{\tau_I(K)}.$$The remaining squares can then be fixed by $\Delta_{K}(t)$.
\epf

\bcor\label{cor: Dehn surgery for torsion order one}
Suppose $K\subset S^3$ is a knot such that $K$ has torsion order one and
$$\dim \sut{\mu}=||\Delta_{K}(t)||.$$
Suppose further that $\tau_I(K)=\tau(K)=\tau$. Then for any $r=p/q\in\mathbb{Q}\backslash\{0\}$ with $q\ge 1$, we have
\[\dim I^{\sharp}(S^3_r(K))=\dim_{\ft}\widehat{HF}(S^3_r(K))=\begin{cases}
(||\Delta_K(t)||+2|\tau|-3)\cdot q/2+ |p-q\cdot (2|\tau|-1)| & \tau >0\\(||\Delta_K(t)||+2|\tau|-3)\cdot q/2+ |-p-q\cdot (2|\tau|-1)| & \tau <0\\(||\Delta_K(t)||-1)\cdot q/2+ |p| & \tau =0.\end{cases}\]
\ecor
\bpf
Since $\tau_I(K)=\tau(K)$, we know from Lemma \ref{cor: bent complex for torsion order one} and \cite[Theorem 1.4]{Petkova2009thin} that the differentials in instanton and Heegaard Floer theory have exactly the same structure. Explicitly, there is one staircase $C_\tau$ and $k$ squares for $$k=\frac{||\Delta_K(t)||-2|\tau|-1}{4}.$$

Despite the difference in coefficients, we can apply the large surgery formulae in \cite{ozsvath2004holomorphicknot} and \cite{LY2021large} to obtain 
$$\dim I^{\sharp}(S^3_{\pm n}(K))=\dim_{\ft}\widehat{HF}(S^3_{\pm n}(K))$$
for any large enough $n$. Explicitly, we have the following.
\benu
\item A square $C$ contributes two-dimensional subspaces for both $(\pm n)$-surgeries
\item A staircase $C_l$ with $l< 0$ contributes an $n$-dimensional subspace for $(-n)$-surgery and a $(n+4|l|-2)$-dimensional subspace for $+n$-surgery.
\item A staircase $C_0$ contributes an $n$-dimensional subspace for both $(\pm n)$-surgeries.
\item A staircase $C_l$ with $l> 0$ contributes an $(n+4l-2)$-dimensional subspace for $(-n)$-surgery and an $n$-dimensional subspace for $+n$-surgery.
\eenu
Note that a figure-eight has one square and a staircase $C_0$ and the torus knot $T(2,2l+1)$ has a staircase $C_l$ and no square. Then the corollary follows from the dimension formulae in \cite[Theorem 1.1]{baldwin2020concordance} and \cite[Proposition 15]{Hanselman20surgery} for both instanton and Heegaard Floer theory.
\epf

\subsection{Induction using oriented skein relation}\label{subsec: alternating knots}
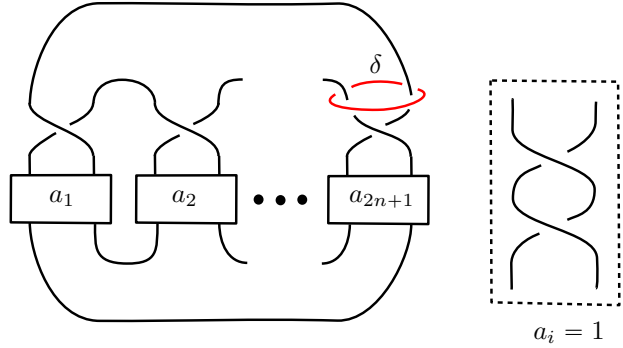
\begin{figure}[ht]
	\begin{overpic}[width=0.6\textwidth]{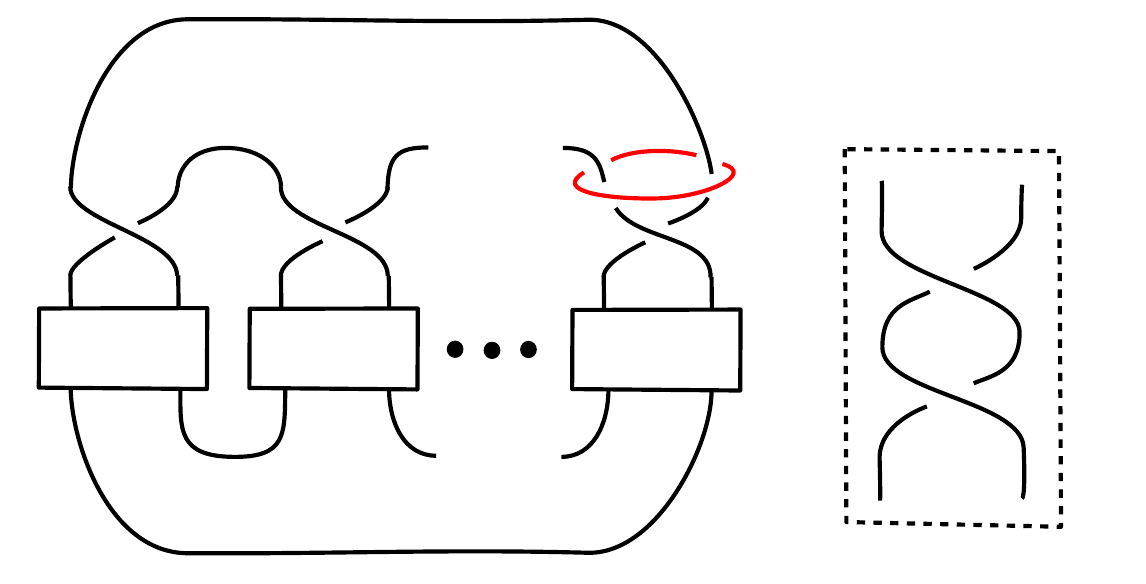}
		\put(9,20){$a_1$}
		\put(27,20){$a_2$}
		\put(53,20){$a_{2n+1}$}
		\put(56,39){$\delta$}
		\put(80,0){$a_i=1$}
	\end{overpic}
	\caption{The knot $K(a_1,\dots,a_{2n+1})$.}\label{fig: alternating2}
\end{figure}

In this subsection, we study differentials for a family of knots
$K(a_1,\dots,a_{2n+1})$, where $a_1,\dots, a_{2n+1}$ are the numbers of full-twists as in Figure \ref{fig: alternating2}. 

\bprop\label{prop: k<=n+1, a_i>0}
Suppose $K=K(a_1,\dots,a_{2n+1})$ with $a_i\ge 0$. Let $$k=\#\{i~|~a_i\ge 1\}.$$ If $k\leq n+1$, then
we have the following.
\begin{enumerate}
	\item $\tau_I(K)=g(K)=n$. 
	\item $K$ has torsion order one (\textit{cf.} Definition \ref{defn: torsion order one}).
\end{enumerate}
\eprop

\bpf[Proof of Theorem \ref{thm: alternating knot}]
When $a_i>0$ for all $i$ we know that $K$ is an alternating knot. It follows from \cite[Corollary 1.8]{baldwin2020concordance} and \cite[Theorem 1.2]{li2019tau} that $\tau_I(K)=\tau(K)$ for all alternating knots. Moreover, by the spectral sequence in \cite{kronheimer2011khovanov}, we know $\dim \sut{\mu}=||\Delta_K(t)||$. Then Proposition \ref{prop: k<=n+1, a_i>0} and Corollary \ref{cor: Dehn surgery for torsion order one} apply.
\epf
% \bcor
% Suppose $K=K(a_1,\dots,a_{2n+1})$ such that $k\leq n+1$ and $a_i>0$ for all $i$, and $r\in\mathbb{Q}\backslash \{0\}$. Then we have
% $$\dim I^{\sharp}(S^3_r(K))={\rm rk}\:\widehat{HF}(S^3_r(K)).$$ 
% \ecor

\begin{figure}[ht]
	\begin{overpic}[width=0.7\textwidth]{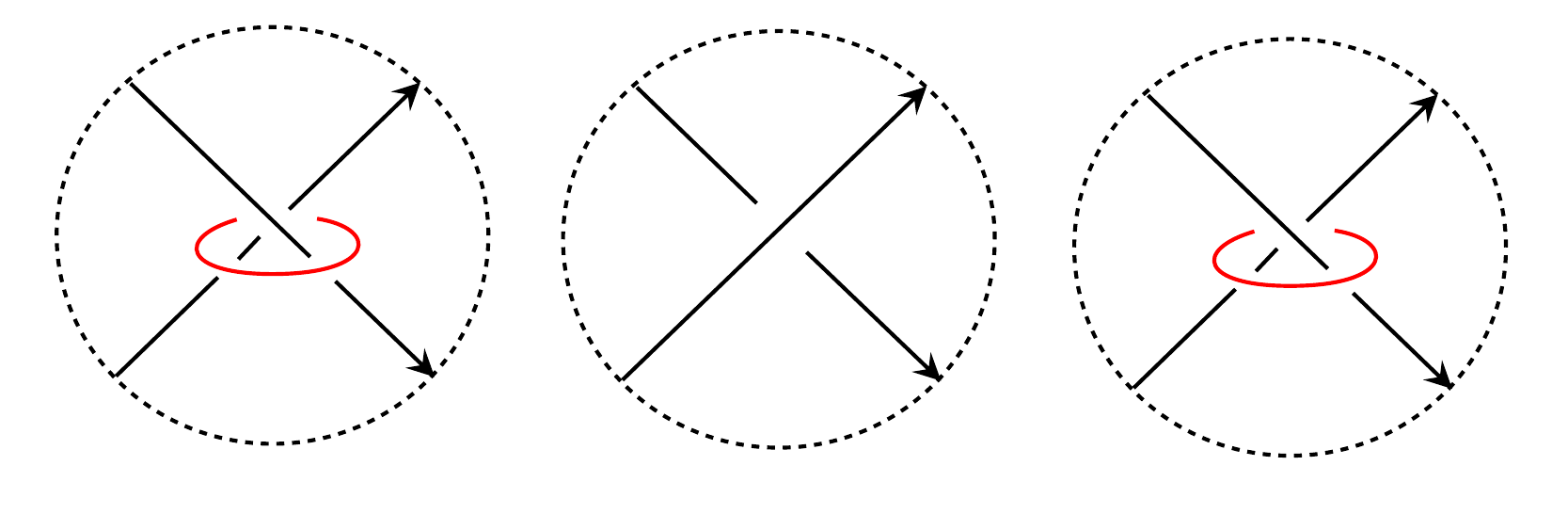}
		\put(10,19){$\delta$}
		\put(16,2){$K_+$}
		\put(48,2){$K_-$}
		\put(88,18){\color{red}$0$}
		\put(80,0){$K_0$}
	\end{overpic}
	\caption{The knots $K_+$, $K_-$ and $K_0$.}\label{fig: skein}
\end{figure}
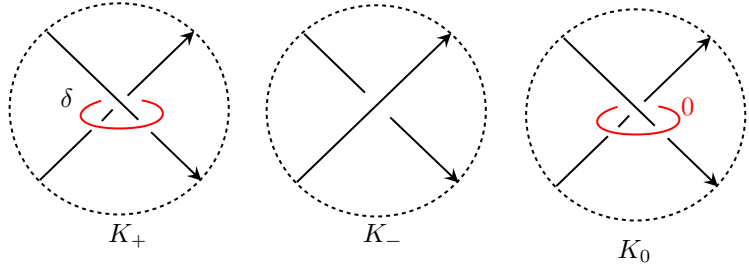

We start with some preparation lemmas.
Suppose $K_+\subset S^3$ is a knot and $\delta$ is a curve circling around a crossing of $K_+$ as in Figure \ref{fig: skein}. Write $K_-\subset S^3_{-1}(\delta)\cong S^3$ and $K_0\subset S^3_{0}(\delta)\cong S^1\times S^2$. We can discuss the tau invariant for the knot $K_0\subset S^1\times S^2$ once we fix an element in $I^{\sharp}(-S^1\times S^2)$. For any $*\in\mathbb{Q}\cup\{\mu\}$, write
$$\sutb{\pm}{*}=\shi(-S^3\backslash N(K_{\pm}),-\Gamma_{*})~{\rm and~}\sutb{0}{*}=\shi(-S^1\times S^2\backslash N(K_{0}),-\Gamma_{*})$$
 The bypass maps are written as $\psi^{\bullet,n_1}_{\pm,n_2}$ for $\bullet\in\{+,-,0\}$. The maps $F_{n}^{\bullet}$, $G_{n}^{\bullet}$ from surgeries along a meridian of $K_+$ are defined similarly. For simplicity, we will write $d_1^\bullet$ for $d_{1,+}^\bullet$.
 
Since $\dim I^{\sharp}(-S^1\times S^2)=2$, we have two effective tau-invariants for the knot $K_0$. To specify choices, first note that there is a surgery exact triangle associated to $\delta$:
\begin{equation}
	\xymatrix{
	I^{\sharp}(-S^3)\ar[rr]^{H_{\delta}}&&I^{\sharp}(-S^3)\ar[dl]^{F_{\delta}}\\
	&I^{\sharp}(-S^1\times S^2)\ar[ul]^{G_{\delta}}&
	}
\end{equation}
Pick $\al_1\neq0\in \im F_{\delta}$ then we can define
$$\tau_{\al_1}(K_0)=\max\{i~|~\exists~x\in\sutbg{0}{n}{i}~s.t.~F_{n}^0(x)=\al_1\}-\frac{n-1}{2}.$$
Note that $\im F_{\delta}$ is $1$-dimensional, so it does not matter what scalar to put on $\al_1$. We pick $$\al_2\in I^{\sharp}(-S^1\times S^2)\backslash \im F_{\delta}$$ such that the value
$$\tau_{\al_2}(K_0)=\max\{i~|~\exists~x\in\sutbg{0}{n}{i}~s.t.~F_{n}^0(x)=\al_2\}-\frac{n-1}{2}$$
takes the maximal value among all possible $\al_2$.

\blem\label{lem: from 0,- to +}
For the knots $K_{+}$, $K_-$, and $K_0$, suppose the following.
\begin{itemize}
	\item[(i)] $\dim\sutb{+}{\mu}=\dim\sutb{-}{\mu}+\dim\sutb{0}{\mu}.$
	\item[(ii)] $\tau_I(K_-)=\tau_{\al_2}(K_0)-1.$
	\item[(iii)]The knots $K_-$ and $K_0$ both have torsion order one.
\end{itemize}
Then $K_+$ has torsion order one.

\elem 
\bpf
The surgery triangle with respect to $\delta$ gives rise to an exact triangle
\begin{equation*}
	\xymatrix{
	\sutb{-}{\mu}\ar[rr]^{H_{\delta,\mu}}&&\sutb{+}{\mu}\ar[dl]^{F_{\delta,\mu}}\\
	&\sutb{0}{\mu}\ar[ul]^{G_{\delta,\mu}}&
	}
\end{equation*}
Condition (i) implies that $G_{\delta,\mu}=0$. Since the surgery maps associated to $\delta$ commute with the differentials $d_1^{\bullet}$ on $\sutb{\bullet}{\mu}$, we have a short exact sequence of chain complexes:
$$0\ra (\sutb{-}{\mu},d^-_1)\to(\sutb{-}{\mu},d^-_1)\to (\sutb{-}{\mu},d^-_1)\to0.$$
The Zigzag lemma gives rise to an exact triangle
\begin{equation}\label{eq: zigzag}
\xymatrix{
	H(\sutb{-}{\mu},d_1^-)\ar[rr]^{H_{\delta,\mu,*}}&&H(\sutb{+}{\mu},d_1^+)\ar[ld]^{F_{\delta,\mu,*}}\\
	&H(\sutb{0}{\mu},d_1^0)\ar[ul]^{\partial_*}&
	}
\end{equation}
From condition (iii) we know that $H(\sutb{-}{\mu},d_1^-)\cong\mathbb{C}$ and $H(\sutb{0}{\mu},d_1^0)\cong\mathbb{C}^2$. So in order to prove the lemma, it suffices to show that $\partial_*\neq 0$. To do this, let $\be_2=G_{\delta}(\al_2)\neq 0\in I^{\sharp}(-S^3)$. Pick $n$ large enough and $x^-\in\sutbg{-}{n}{\tau_I(K_-)}$ such that $F_{n}^-(x^-)=\be$. Take $u^-=\psp{-,n}{\mu}(x^-)$, where the first superscript of $\psp{-,n}{\mu}$ corresponds to the superscript of $\sutb{\pm}{\mu}$, and the first subscript corresponds to the positve and the negative bypasses. It is straightforward to check that $d^-_1(u^-)=0$ and $u^-\notin\im d^-_1$. Then $H(\sutb{-}{\mu},d_1^-)$ is generated by $[u^-]$. 

Pick $x_2^0\in\sutbg{0}{n}{\tau_{\al_2}(K^0)}$ such that $F_{n}^0(x_2^0)=\al_2$. We have the following diagrams in which the triangles are exact and the parallelograms are commutative.
\begin{equation}\label{eq: commutative, al, psi}
\xymatrix{
	\sutb{-}{n}\ar[rrrr]^{H_{\delta,n}}\ar[dd]^{\psp{-,n}{\mu}}&&&&\sutb{+}{n}\ar[dll]^{F_{\delta,n}}\ar[dd]^{\psp{+,n}{\mu}}\\
	&&\sutb{0}{n}\ar[ull]^{G_{\delta,n}}\ar[dd]^<<<<{\psp{0,n}{\mu}}&&\\
	\sutb{-}{\mu}\ar[rrrr]^<<<<<<<<<{H_{\delta,\mu}}&&&&\sutb{+}{\mu}\ar[dll]^{F_{\delta,\mu}}\\
	&&\sutb{0}{\mu}\ar[ull]^{G_{\delta,\mu}}&&
}
\xymatrix{
	\sutb{-}{n}\ar[rrrr]^{H_{\delta,n}}&&&&\sutb{+}{n}\ar[dll]^{F_{\delta,n}}\\
	&&\sutb{0}{n}\ar[ull]^{G_{\delta,n}}&&\\
	\sutb{-}{\mu}\ar[rrrr]^<<<<<<<<<{H_{\delta,\mu}}\ar[uu]^{\psp{-,\mu}{n}}&&&&\sutb{+}{\mu}\ar[dll]^{F_{\delta,\mu}}\ar[uu]^{\psp{\mu}{n}}\\
	&&\sutb{0}{\mu}\ar[ull]^{G_{\delta,\mu}}\ar[uu]_>>>>{\psp{0,\mu}{n}}&&
}
\end{equation}

{\bf Claim}. We have $G_{\delta,n}(x_2^0)=U(x^-)$. 

\bpf[Proof of Claim.]
Suppose not, \textit{i.e.}, $y=G_{\delta,n}(x_2^0)-Ux^-\neq 0$. Note
$$F_{n}^-(y)=F_{n}^-\circ G_{\delta,n}(x_2^0)-F_{n}^-\circ U(x^-)=G_{\delta}(\al)-G_{\delta}(\al)=0.$$
Hence the fact that $K^-$ has torsion order one implies that $y\notin \im U$. As a result, 
$$v=\psp{-,n}{\mu}(y)=\psp{-,n}{\mu}\circ G_{\delta,n}(x^0)\neq 0.$$
From commutativity, we know
$$G_{\delta,\mu}\circ\psp{0,n}{\mu}(x_2^0)\neq0$$
which contradicts the fact that $G_{\delta,\mu}=0$ as implied by Condition (i).
\epf

Now since $Ux^-=G_{\delta,n}(x_2^0)$, we know from commutativity that
$$U\circ H_{\delta,n}(x^-)=H_{\delta,n}\circ U(x^-)=0.$$
Then there exists $u^+\in\sutb{+}{\mu}$ such that $H_{\delta,n}(x^-)=\psp{+,n}{\mu}(u^+).$ Take $u^0=F_{\delta,\mu}(u^+)$, we know from the commutativity that
$$\psp{0,n}{\mu}(u^0)=F_{\delta,n}\circ \psp{+,n}{\mu}(u^+)=0.$$
As a result, $d_1^0(u^0)=0$. Also, we know
\beq
d_1^+(u^+)&=\psp{+,n}{\mu}\circ\psp{+,\mu}{n}(u^+)\\
&=\psp{+,n}{\mu}\circ H_{\delta,n}(x^-)\\
&=H_{\delta,\mu}\circ\psp{-,n}{\mu}(x^-)\\
&=H_{\delta,\mu}(u^-).
\eeq
By the construction of $\partial_*$, we conclude that
$\partial_*([u^0])=[u^-]$
and we conclude the proof of the lemma.

\epf

\blem\label{lem: tau(K_+)=tau_al_1}
For the knots $K_{+}$, $K_-$, and $K_0$, suppose we have the following.
\begin{itemize}
	\item [(i)] All three knots $K_-$, $K_0$ and $K_+$ have torsion order one.
	\item [(ii)] Either $\dim\sutb{+}{\mu}=\dim\sutb{-}{\mu}+\dim\sutb{0}{\mu}$, or $\dim\sutb{0}{\mu}=\dim\sutb{-}{\mu}+\dim\sutb{+}{\mu}.$ 
\end{itemize}
Then $\tau_{\al_1}(K_0)=\tau_I(K_+)$.
\elem

\bpf
From Condition (ii) and the zigzag lemma there exist an exact triangle
\begin{equation}\label{eq: zigzag, 2}
\xymatrix{
	H(\sutb{-}{\mu},d_1^-)\ar[rr]&&H(\sutb{+}{\mu},d_1^+)\ar[ld]^{F_{\delta,\mu,*}}\\
	&H(\sutb{0}{\mu},d_1^0)\ar[ul]&
	}
\end{equation}
Condition (i) implies that $H(\sutb{+}{\mu},d_1^+)\cong\mathbb{C}$ and $F_{\delta,\mu,*}\neq0$. Take $\be_1\in I^{\sharp}(-S^3)$ such that
$F_{\delta}(\be_1)=\al_1.$
Take $x^+\in\sutbg{+}{n}{\tau_I(K_+)}$ such that
$F_{n}^+(x^+)=\be_1.$
We know from the commutativity that
$$F_{n}^0\circ F_{\delta,n}(\be_1)=F_{\delta}\circ F_{n}^+(x^+)=\al_1.$$
Hence by the definition of $\tau$ we know
$\tau_{\al_1}(K_0)\geq \tau_I(K_+).$
Suppose $\tau_{\al_1}(K_0)=\tau_I(K_+)+k$ for some $k>0$.
Now take $v^+=\psp{+,n}{\mu}(x^+)$. We know that $d^+_1(v^+)=0$ and $v^+\notin \im d^+_1$. So $H(\sutb{+}{\mu},d_1^+)\cong\mathbb{C}$ is generated by $[v^+]$. Let $v^0=F_{\delta,\mu}(v^{+})$, we know that 
$F_{\delta,\mu,*}([v^+])=[v^0].$

Pick $x^0_1\in\sutbg{0}{n}{\tau_{\al_1}(K_0)}$ such that $F_{n}^0(x^0_2)=\al_1$, then we know that
$$F_{n}^0\bigg(F_{\delta,n}(x^+)-U^k(x^0_1)\bigg)=\al_1-\al_1=0.$$
Since $K_0$ has torsion order one as in Condition (i), we know that
$$U\bigg(F_{\delta,n}(x^+)-U^k(x^0_1)\bigg)=0.$$
As a result, there exists $w^+\in\sutb{0}{\mu}$ such that
$F_{\delta,n}(x^+)-U^k(x^0_1)=\psp{+,\mu}{n}(w^0).$
As a result, we know that
\beq
v^0&=F_{\delta,\mu}(v^{+})\\
&=\psp{0,n}{\mu}\circ F_{\delta,n}(x^+)\\
&=\psp{0,n}{\mu}\bigg(F_{\delta,n}(x^+)-U^k(x^0_1)\bigg)\\
&=\psp{0,n}{\mu}\psp{+,\mu}{n}(w^0)\\
&=d_1^0(w^0).
\eeq
As a result, we know that
$F_{\delta,\mu,*}([v^+])=[v^0]=0,$ which contradicts the fact that $F_{\delta,\mu,*}$ fits into the exact triangle (\ref{eq: zigzag, 2}) and the fact that $F_{\delta,\mu,*}\neq 0$.
\epf

\blem\label{lem: tau_al_2=tau(K_-)}
For the knots $K_{+}$, $K_-$, and $K_0$, suppose we have the following.
\begin{itemize}
	\item [(i)] The knots $K_-$ and $K_+$ both have torsion order one.
	\item [(ii)] We have
	$\dim\sutb{0}{\mu}=\dim\sutb{-}{\mu}+\dim\sutb{+}{\mu}$
\end{itemize}
Then we have the following.
\begin{enumerate}
	\item $K_0$ has torsion order one.
	\item We have
	$\tau_{\al_2}(K_0)=\tau_I(K_-).$
\end{enumerate}
\elem
\bpf
From Condition (ii) and the Zigzag lemma, we have an exact triangle
\begin{equation}\label{eq: zigzag, 3}
\xymatrix{
	H(\sutb{-}{\mu},d_1^-)\ar[rr]&&H(\sutb{+}{\mu},d_1^+)\ar[ld]^{F_{\delta,\mu,*}}\\
	&H(\sutb{0}{\mu},d_1^0)\ar[ul]^{G_{\delta,\mu,*}}&
	}
\end{equation}
From Condition (i) we know
$H(\sutb{+}{\mu},d_1^+)\cong H(\sutb{+}{\mu},d_1^+)\cong\mathbb{C}$
so
$\dim H(\sutb{0}{\mu},d_1^0)\leq 2.$
On the other hand, we know
$\dim H(\sutb{0}{\mu},d_1^0)\geq \dim I^\sharp(-S^1\times S^2)=2$. As a result, we know
$\dim H(\sutb{0}{\mu},d_1^0)= 2$
which implies that $K_0$ has torsion order one, $F_{\delta,\mu,*}\neq0$, and $G_{\delta,\mu,*}\neq0$.

Note that all hypotheses of Lemma \ref{lem: tau(K_+)=tau_al_1} are satisfied so the argument in the proof of that lemma applies. In particular, we know
$\tau_{\al_1}(K_0)=\tau_I(K_+).$
We can pick $x^+\in\sutbg{+}{n}{\tau_I(K^+)}$ such that $F_{n}^+(x^+)=\be_1$, and take
$$x^0_1=F_{\delta,n}(x^+),~v^0_1=\psp{0,n}{\mu}(x^0_1),~{\rm and~}v^+=\psp{0,n}{\mu}(x^+).$$
We know from the proof of Lemma \ref{lem: tau(K_+)=tau_al_1} that
$F_{\delta,\mu,*}([v^+])=[v^0_1]\neq0.$

Pick $x^0_2\in\sutbg{0}{n}{\tau_{\al_2}(K)}$ such that $F_n^0(x^0_2)=\al_2$ and take $v^0_2=\psp{0,n}{\mu}(x^0_2)$. It is straightforward to check that $d^0_1(v^0_2)=0$ and $v^0_2\notin\im d^0_1.$

{\bf Claim}. The homology $H(\sutb{0}{\mu},d^0_1)$ is generated by $[v^0_1]$ and $[v^0_2]$.
\bpf[Proof of Claim.]
When $v^0_1$ and $v^0_2$ have different gradings in $\sutb{0}{\mu}$, the claim follows immediately from the fact that $\dim H(\sutb{0}{\mu},d^0_1)=2$. If $v^0_1$ and $v^0_2$ have the same grading, then $x^0_1$ and $x^0_2$ have the same grading in $\sutb{0}{n}$, and hence
$\tau_{\al_1}(K_0)=\tau_{\al_2}(K_0).$
In this case, suppose that there exists complex numbers $c_1$, $c_2$, not both zero, and an element $w^0\in\sutb{0}{\mu}$ such that
$$c_1\cdot v^0_1+c_2\cdot v^0_2+d_1^0(w^0)=0.$$
Take $y^0=\psp{0,\mu}{n}(w^0)$, then the above equality is equivalent to
$$\psp{0,n}{\mu}(c_1\cdot x^0_1+c_2\cdot x^0_2+y^0)=0,$$
which implies that there exists $z^0\in\sutbg{0}{n}{\tau_{\al_1}(K_0)+1=\tau_{\al_2}(K_0)+1}$ such that
$$c_1\cdot x^0_1+c_2\cdot x^0_2+y^0+Uz^0=0.$$
From the construction we know $Uy^0=0$ hence $F_{n}^0(y^0)=0$. Since the grading of $z^0$ is strictly larger than both $\tau_{\al_1}(K_0)$ and $\tau_{\al_2}(K_0)$, the choice of $\al_1$ and $\al_2$ implies that $F_n^0(z)=0$. As a result we have
$$c_1\cdot\al_1+c_2\cdot\al_2=F_{n}^0(c_1\cdot x^0_1+c_2\cdot x^0_2+y^0+Uz^0)=0.$$
Thus we must have $c_1=c_2=0$ and hence $[v^0_1]$ and $[v^0_2]$ are linearly independent.
\epf

With the help of the above claim, the proof that $\tau_{\al_2}(K_0)=\tau_I(K_-)$ is similar to the proof of Lemma \ref{lem: tau(K_+)=tau_al_1}.
\epf

\bpf[Proof of Proposition \ref{prop: k<=n+1, a_i>0}]
We use induction on $k$ to prove the following.
\begin{itemize}
	\item [(i)] The knot $K$ has genus $g(K)=n$.
	\item [(ii)] The coefficient of the term $t^i$ in $\Delta_K(t)$ has sign $(-1)^{n-i}$.
	\item [(iii)] The knot has torsion order one.
	\item [(iv)] We have $\tau_I(K)=n$.
\end{itemize}

When $k=0$, the knot is the torus knot $T(2,2n+1)$, so all the above four statements hold. Suppose we have proved the above four statements for $k$. Now we deal with the case $k+1$. Without loss of generality, we can assume that $a_{2n+1}>1$ while $a_{2n}=1$. Write
$$K_{l}=K(a_1,a_2,a_3,\dots,a_{2n},l).$$
We know $K=K_{a_{2n+1}}$. Note that $K_{-1}=K(a_1,a_2,\dots,a_{2n-1})$ and $K_1=K(a_1,\dots,a_{2n}=1,1)$, so the inductive hypothesis applies to both $K_{-1}$ and $K_1$. Let $\delta$ be a curve circling around the crossing corresponding to $a_{2n+1}$ as shown in Figure \ref{fig: alternating2}. Then we can take $K_+=K_1$, $K_-=K_{-1}$ and there is a corresponding $K_0\subset S^1\times S^2$. From \cite{kronheimer2010instanton} we know that
$$\chi_{\mathrm{gr}}(\sutb{\pm}{\mu})=-\Delta_{K_{\pm}}(t){\rm~and~}\chi_{\mathrm{gr}}(\sutb{+}{\mu})-\chi_{\mathrm{gr}}(\sutb{-}{\mu})=\chi_{\mathrm{gr}}(\sutb{0}{\mu}).$$
Also, since $K_{\pm1}$ are both alternating knots, we know that
$$\dim \sutb{\pm}{\mu}=||\chi_{\mathrm{gr}}(\sutb{\pm}{\mu})||,$$
where $||\cdot||$ denotes the sum of the absolute values of coefficients. Statements (i) and (ii) applied to $K_{\pm 1}$ then imply that
\beq
\dim \sutb{0}{\mu}&\geq||\chi_{\mathrm{gr}}(\sutb{0}{\mu})||\\
&=||\chi_{\mathrm{gr}}(\sutb{+}{\mu})-\chi_{\mathrm{gr}}(\sutb{+}{\mu})||\\
&=||\chi_{\mathrm{gr}}(\sutb{+}{\mu})||+||\chi_{\mathrm{gr}}(\sutb{+}{\mu})||\\
&=\dim\sutb{+}{\mu}+\dim\sutb{-}{\mu}
\eeq
Then it follows from the exact triangle
\begin{equation*}
	\xymatrix{
	\sutb{-}{\mu}\ar[rr]^{H_{\delta,\mu}}&&\sutb{+}{\mu}\ar[dl]^{F_{\delta,\mu}}\\
	&\sutb{0}{\mu}\ar[ul]^{G_{\delta,\mu}}&
	}
\end{equation*}
that
$\dim \sutb{0}{\mu}=\dim\sutb{+}{\mu}+\dim\sutb{-}{\mu}.$
Then we can apply Lemma \ref{lem: tau(K_+)=tau_al_1} and Lemma \ref{lem: tau_al_2=tau(K_-)} to conclude that $K_0$ has torsion order one, and $\tau_{\al_1}(K_0)=\tau_{\al_2}(K_0)+1=g.$

Now for any odd $l>0$, we can take $K_+=K_l$, $K_-=K_{l-2}$, and take $K_0$ to be the same knot as the one for $K_1$ and $K_{-1}$. Hence we can apply Lemma \ref{lem: torsion order 1} to inductively conclude all four statements. 

In the statement of Proposition \ref{prop: k<=n+1, a_i>0}, we require that $k\leq n+1$. This extra assumption is due to our strategy is to cancel two crossings when $a_{2n}=1$ and $a_{2n+1}=-1$. In particular, in the proof we need $a_{2n}=1$ throughout the induction so that we have enough information to start with to understand larger $a_{2n+1}$. This means at the very beginning we need at least half of $a_{i}$ to be $1$.
\epf

% \brem
% Here we present Proposition \ref{prop: k<=n+1, a_i>0} as a proof of usefulness. 
% \erem

\section{Twisted Whitehead doubles and splicings}\label{sec: twisted Whitehead doubles}
The techniques in Section \ref{subsec: alternating knots} can also be used to study twisted Whitehead doubles.

\begin{figure}[ht]
	\begin{overpic}[width=0.7\textwidth]{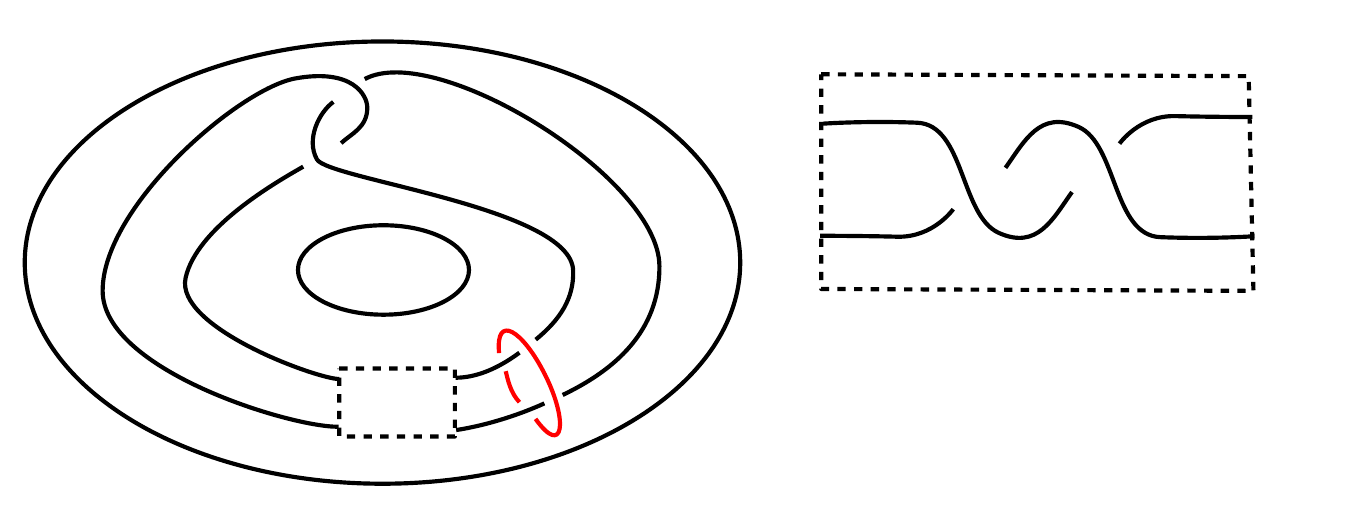}
		\put(10,16){$\widehat{K}$}
		\put(-2,16){$\widehat{V}$}
		\put(28.5,8){$t$}
		\put(73,12){$t=1$}
	\end{overpic}
	\caption{The Whitehead double.}\label{fig: pattern of WHD}
\end{figure}

\bdefn
Suppose $\widehat{V}\subset S^3$ is an unknotted solid torus. Let $\widehat{K}\subset \widehat{V}$ be the knot as shown in Figure \ref{fig: pattern of WHD}. Let $\hat{\mu}$ be a non-separating curve on $\partial \widehat{V}$ bounding a disk in $\widehat{V}$ and let $\hat\lambda$ be a non-separating curve on $\partial \widehat{V}$ bounding a disk in $S^3\backslash \widehat{V}$. Let $J$ be a knot in $S^3$ and $V$ a tubular neighborhood of $J$. Let $\mu$ and $\lambda$ be the meridian and Seifert longitude of $J$, respectively. Let
$$f:\widehat{V}\hookrightarrow S^3$$
be an embedding such that $f(\widehat{V})=V$, $f(\hat{\mu})=\mu$ and $f(\hat{\lambda})=\lambda$. Let $K=f(\widehat{K})$. Then $K$ is called the {\bf positively clasped $t$-twist Whitehead double} of $J$, denoted by $D_t^{+}(J)$.
\edefn

\brem\label{rem: negative clasped}
We can also study the negatively clasped Whitehead doubles. Note that they are the mirror of positively clasped Whitehead doubles as in \cite{Hedden2007WHD}.
\erem

Here are some basic properties of $K$.
\blem[\cite{Hedden2007WHD}]\label{lem: basic properties of WHD}
Suppose $K$ is the positively clasped $t$-twist Whitehead double of $J$. Then we have the following.
\begin{enumerate}
	\item The genus of $K$ is one. 
	\item $\Delta_K(T)=-t\cdot T+(2t+1)-t\cdot T^{-1}.$
\end{enumerate}
\elem

Since the Whitehead doubles all have genus $1$, there are only three nontrivial gradings of its $\khii$ to study. Note that the top and bottom gradings are isomorphic to each other. The following lemma describes the top (and hence the bottom) grading.

\blem\label{lem: top grading of KHI(WHD)}
Suppose $K$ is the positively clasped $t$-twist Whitehead double of $J$. Then
$$\khii(S^3,K,1)\cong \shi(S^3\backslash N(J),\Ga_{-t}).$$
\elem

\begin{figure}[ht]
	\begin{overpic}[width=0.7\textwidth]{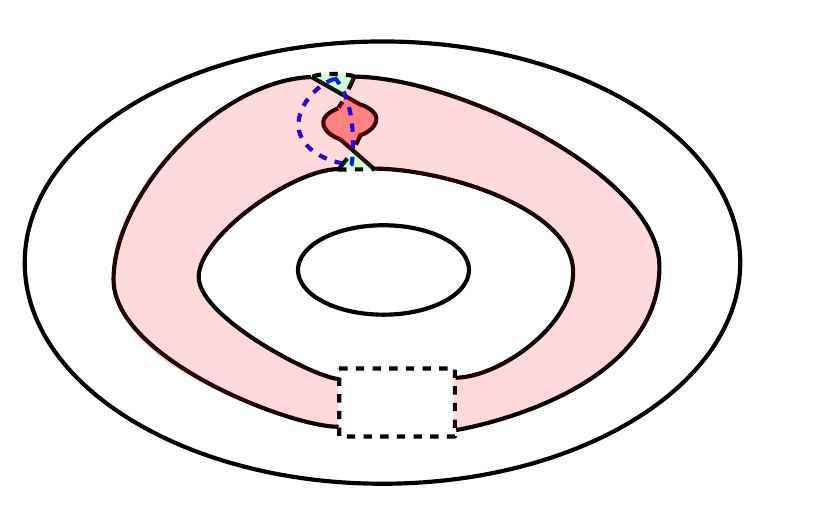}
		\put(10,27){$\widehat{K}$}
		%\put(-2,16){$\widehat{V}$}
		\put(47,13){$t$}
		\put(33,47){$\alpha$}
	\end{overpic}
	\caption{A genus-one Seifert surface $S$ of $\widehat{K}$. Two sides of the Seifert surface are shaded in red and green. The blue curve $\al$ is a curve on $S$ bounding a disk in $\widehat{V}\backslash S$. This disk can be viewed as a product disk in the sutured manifold $(S^3\backslash [-1,1]\times S,\{0\}\times\partial S)$.}\label{fig: Seifert surface of WHD}
\end{figure}

\bpf
A genus-one Seifert surface $S$ of $K$ can be drawn as in Figure \ref{fig: Seifert surface of WHD} (inside $\widehat{V}$). From the proof of \cite[Proposition 7.16]{kronheimer2011knot}, we know that there is an isomorphism
$$\khii(S^3,K,1)\cong \shi(S^3\backslash [-1,1]\times S,\{0\}\times\partial S).$$
As shown in Figure \ref{fig: Seifert surface of WHD}, the sutured manifold $(S^3\backslash [-1,1]\times S,\{0\}\times\partial S)$ admits a product disk $D$ and it is straightforward to check that there is a sutured manifold decomposition
$$(S^3\backslash [-1,1]\times S,\{0\}\times\partial S)\stackrel{D}{\leadsto}(S^3\backslash N(J),\Gamma_{-t}).$$
Note that we can fix the suture as $\Ga_{-t}$ by counting its intersections with $\mu$ and $\lambda$ explicitly. Hence we conclude that
\beq
\khii(S^3,K,1)&\cong \shi(S^3\backslash [-1,1]\times S,\{0\}\times\partial S)\\
&\cong \shi(S^3\backslash N(J),\Ga_{-t}).
\eeq
\epf

Now we compute the tau invariants for the twisted Whitehead doubles.
\blem\label{lem: tau for Whitehead doubles}
Suppose $K$ is the positively clasped $t$-twist Whitehead double of $J$. Then
\begin{equation*}
	\tau_I(K)=\begin{cases}
		1&t<2\cdot \tau_I(J)\\
		0&t\geq 2\cdot \tau_I(J)
	\end{cases}
\end{equation*}
\elem
\bpf
Write
$(K^+_t,\Ga_n)=\shi(-S^3(K_t^+),-\Ga_n),~{\rm and~}(K^+_t,\Ga_n,i)=\shi(-S^3(K_t^+),-\Ga_n,i).$
We take the surgery exact triangle along the curve $\delta$. The maps in the surgery triangle associated to $\delta$ commute with the $2$-handle attachments along the meridian of the knots, so we have the following diagram, in which the triangles are exact and the parallelograms are commutative.
\begin{equation}\label{eq: WHD, 1}
	\xymatrix{
	(D^+_{t+1}(J),\Ga_n)\ar[rrrr]^{H_{\delta,n}}\ar[dd]^{F_{t+1,n}}&&&&(D^+_{t}(J),\Ga_n)\ar[dll]^{F_{\delta,n}}\ar[dd]^{F_{t,n}}\\
	&&(K^+,\Ga_n)\ar[ull]^{G_{\delta,n}}\ar[dd]^<<<<{F_{\times,n}}&&\\
	I^{\sharp}(-S^3)\ar[rrrr]^<<<<<<<<<<<<<<{H_{\delta}}&&&&I^{\sharp}(-S^3)\ar[dll]^{F_{\delta}}\\
	&&I^{\sharp}(-S^1\times S^2)\ar[ull]^{G_{\delta}}&&
}
\end{equation}
Note that the above diagram works for any $n\in\mathbb{Z}$, but we fix a large enough $n\in\intg$. Here $K^+\subset S^1\times S^2$ is obtained from $D_t^{+}(J)$ by performing a $0$-surgery along $\delta$. Note that $g(K^+)=1$. Also, the knots $D_t^{+}(J)$ depend on the companion knot $J$, yet due to the $3$-dimensional light bulb theorem we know that $K^+$ is independent of the companion knot $J$. As a result, we can assume $J=U$ is the unknot to obtain information about $J$. When $J=U$, we know $D^+_{0}(U)$ is the unknot and $D^+_{-1}(U)$ is the right-handed trefoil. As a result, we know from Lemma \ref{lem: tau(K_+)=tau_al_1} and Lemma \ref{lem: tau_al_2=tau(K_-)} that
$$\tau_{\al_1}(K^+)=1,~{\rm and~}\tau_{\al_2}(K^+)=0.$$
Note that we also know from \cite[Section 6]{li2019tau} that
$$(D_{0}^+(U),\Ga_{n},1)=0~{\rm and}~(D_{-1}^+(U),\Ga_{n},1)\cong\mathbb{C}.$$
As a result, we know from the exactness that $(K^+,\Ga_n,1)\cong\mathbb{C}$. Also, from Lemma \ref{lem: symmetry and dimension of Gamma_n} part (2), we know
\begin{equation}\label{eq: WHD, 2}
	\dim (J,\Ga_{n})=\dim (J,\Ga_{-2\cdot\tau_I(J)})+|n+2\cdot\tau_I(J)|.
\end{equation}
As a result, when $t<2\tau_I(J)$, we know from Lemma \ref{lem: top grading of KHI(WHD)} and Lemma \ref{lem: bypass n,n+1,mu} that
$$\dim (D^+_{t+1}(J),\Ga_n,1)=\dim (D^+_{t+1}(J),\Ga_{\mu},1)=\dim (D^+_{t}(J),\Ga_{\mu},1)-1=\dim (D^+_{t}(J),\Ga_{n},1)-1.$$
Hence $F_{\delta,n}$ restricted to $(D^+_{t}(J),\Ga_{n},1)$ is nontrivial. Since $\tau_{\al_1}(K^+)=1$ ({\it cf.} Section \ref{subsec: commutativity for d_1}) and $F_{\delta}$ is injective, we know that
$$F_{t,n}|_{(D^+_{t+1}(J),\Ga_n,1)}\neq 0$$
which implies that
$\tau_I(D^+_{t+1}(J))=1.$

When $t\geq 2\tau_I(J)$, we know similarly that 
$$\dim (D^+_{t+1}(J),\Ga_n,1)=\dim (D^+_{t}(J),\Ga_n,1)+1,$$
so $F_{\delta,n}$ restricted to $(D^+_{t}(J),\Ga_{n},1)$ is trivial and the injectivity of $F_{\delta}$ implies that 
$$F_{t,n}|_{(D^+_{t+1}(J),\Ga_n,1)}=0$$
which means $\tau_I(D^+_{t+1}(J))<1.$ To further refine the $\tau_I$, we can look at the mirrors of such knots. Taking the mirror corresponding to reversing the orientation of the $3$-manifold so we have a different diagram
\begin{equation}\label{eq: WHD, 3}
	\xymatrix{
	(\overline{D^+_{t}(J)},\Ga_n)\ar[rrrr]^{\widebar{H}_{\delta,n}}\ar[dd]^{\widebar{F}_{t,n}}&&&&(\overline{D^+_{t+1}(J)},\Ga_n)\ar[dll]^{\widebar{F}_{\delta,n}}\ar[dd]^{\widebar{F}_{t+1,n}}\\
	&&(\widebar{K}^+,\Ga_n)\ar[ull]^{\widebar{G}_{\delta,n}}\ar[dd]^<<<<{\widebar{F}_{\times,n}}&&\\
	I^{\sharp}(-S^3)\ar[rrrr]^<<<<<<<<<<<<<<{H_{\delta}}&&&&I^{\sharp}(-S^3)\ar[dll]^{F_{\delta}}\\
	&&I^{\sharp}(-S^1\times S^2)\ar[ull]^{G_{\delta}}&&
}
\end{equation}
As above, we can use the case $J=U$ and $t=-1$ to compute that
$$\tau_{\al_1}(K^+)=0,~{\rm and~}\tau_{\al_2}(K^+)=-1.$$
As a result we have $\widebar{F}_{\times,n}\circ \widebar{F}_{\delta,n}$ is trivial and the injectivity of $\widebar{F}_{\delta}$ implies that 
$$\widebar{F}_{t,n}|_{(\overline{D^+_{t+1}(J)},\Ga_n,1)}=0,$$
which means $\tau_I(\overline{D^+_{t+1}(J)})<1$ for any $t\in\intg$. In particular, for $t\geq 2\cdot\tau_{I}(J)$, we must have $\tau_I(D^+_{t+1}(J))=0.$
\epf

\brem
Lemma \ref{lem: tau for Whitehead doubles} answers \cite[Question 1.25]{baldwin2020concordance} affirmatively.
\erem

\begin{figure}[ht]
	\begin{overpic}[width=0.7\textwidth]{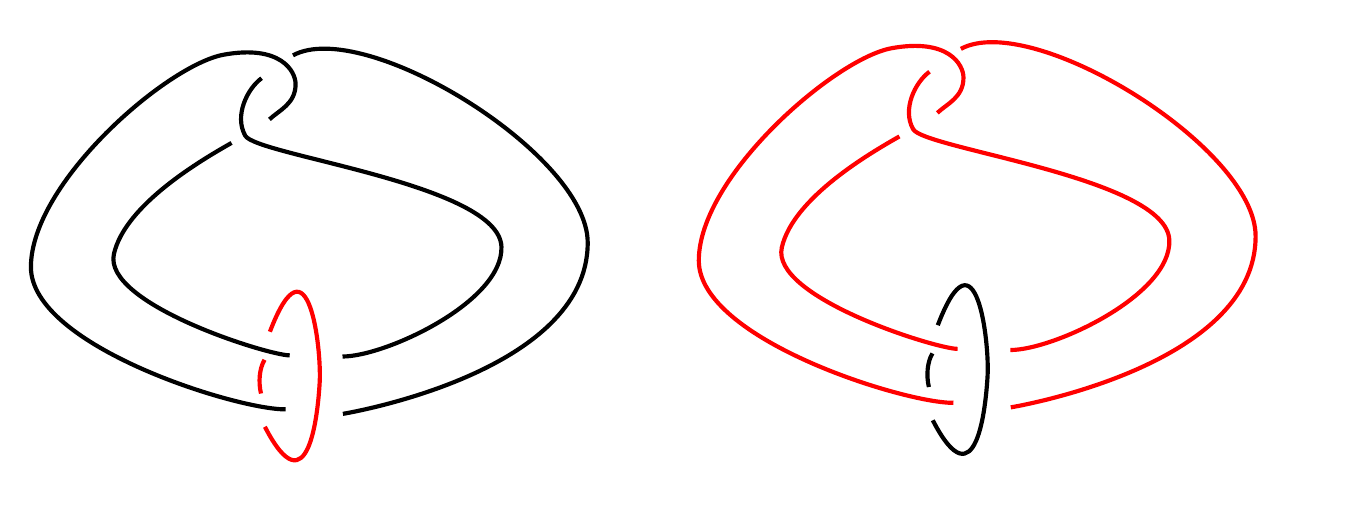}
		\put(10,27){$\widehat{K}$}
		\put(60,27){$L$}
		\put(40,28){$\pm1$}
		\put(18,17){$L$}
		\put(67,17){$\widehat{K}$}
		\put(72,3){$\pm1$}
		%\put(33,47){$\alpha$}
	\end{overpic}
	\caption{The two components of a Whitehead double link can be swapped.}\label{fig: swap two components of WHD}
\end{figure}
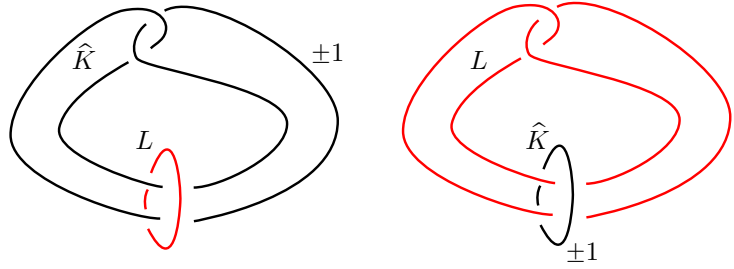

\bpf[Proof of Theorem \ref{thm: Whitehead double}]
Part (1) and part (2) are Lemma \ref{lem: top grading of KHI(WHD)} and Lemma \ref{lem: tau for Whitehead doubles}. Part (3) follows from Lemma \ref{lem: basic properties of WHD}, Lemma \ref{lem: tau for Whitehead doubles}, and Corollary \ref{cor: 1-surgery of genus 1 knot}. 
\epf
\bpf[Proof of Theorem \ref{thm: splicing with twist knot}]
Let $K$ be the positively clasped $0$-twist Whitehead double of $J$. Let $L\subset \partial \widehat{V}$ be a meridian of $\widehat{V}$ as in Figure \ref{fig: swap two components of WHD}. The knot $S^3_{-\frac{1}{n}}(K)$ can be viewed as the splicing of the complements of the knot $J\subset S^3$ and the knot $L\subset S^{3}_{-\frac{1}{n}}(\widehat{K})$. It is well known that the two components of the Whitehead link can be swapped so $S^{3}_{-\frac{1}{n}}(\widehat{K})$ is still $S^3$, while the knot $L\subset S^{3}_{-\frac{1}{n}}(\widehat{K})$ becomes the knot $K_n$. Theorem \ref{thm: Whitehead double} part (3) applies to compute the $(\pm1)$-surgeries of the knot $K$. Then we can apply \cite[Theorem 1.1]{baldwin2020concordance} after knowing $(\pm1)$-surgeries.
\epf

% \bcor\label{cor: splicing with figure 8}
% Suppose $J\subset S^3$ is a knot. Let $Y$ be obtained by gluing the complement of $J$ with the complement of a figure-eight knot so that the gluing map sends the meridian of one knot to the longitude of the other and vice versa. Then
% \begin{equation*}
% 	\dim I^{\sharp}(Y)=2\cdot\dim \shi(S^3\backslash N(J),\Ga_0)+1.
% \end{equation*}
% \ecor

% \bpf
% Let $K$ be the positively clasped $0$-twist Whitehead double of $J$. Let $L\subset \partial \widehat{V}$ be a meridian of $\widehat{V}$ as in Figure \ref{fig: swap two components of WHD}. It is clear that $S^3_{-1}(K)$ can be viewed as the splicing of the complements of the knot $J\subset S^3$ and the knot $L\subset S^{3}_{-1}(\widehat{K})$. It is well known that the two components of the Whitehead link can be swapped so $S^{3}_{-1}(\widehat{K})$ is still $S^3$, while the knot $L\subset S^{3}_{-1}(\widehat{K})$ becomes the figure-eight. Hence Corollary \ref{cor: surgery along WHD} applies.
% \epf

\section{Almost L-space knots}\label{sec: Almost L-space knots}
% \bdefn\label{defn: almost L-space knot}
% A knot $K\subset S^3$ is said to be an {\bf almost L-space knot} if $K$ is not an $L$-space knot and there exists $n\in\intg$, so that
% $$\dim I^{\sharp}(S^3_n(K))\cong\mathbb{C}^{|n|+2}.$$
% \edefn

In this section, we study almost $L$-space knots; see (\ref{eq: defn almost l}) for the definition. We adopt the following notations from \cite[Definition 5.2]{LY2021large}: Let $K\subset S^3$ be a knot (such that $(p,q)=(1,0))$. Define 
$$T_{n,i}=\sutg{n}{i+\frac{n-1}{2}}~{\rm and~}B_{n,i}=\sutg{n}{i-\frac{n-1}{2}}.$$
By Lemma \ref{lem: structure of Gamma_n} part (2), we know that when $n\geq 2g(K)+1$,
$$T_{n,i}\cong T_{n+1,i}~{\rm and~}B_{n,i}\cong B_{n+1,i}.$$
We can rewrite the bypass exact triangles in Lemma \ref{lem: bypass n+1,n,2n+1/2} using $T_{n,i}$ and $B_{n,i}$ as follows.
\blem\label{lem: bypasses for T and B}
Adopting the notations as above, we have the following two bypass exact triangles:
\begin{center}
\begin{minipage}{0.4\textwidth}
	\begin{equation*}
\xymatrix{
\sutg{\frac{2n-1}{2}}{i}\ar[rr]^{\psp{\frac{2n-1}{2}}{n}}&&T_{n,i}\ar[dl]^{\psp{n}{n-1}}\\
&B_{n-1,i-1}\ar[ul]^{\psp{n-1}{\frac{2n-1}{2}}}&
}	
\end{equation*}
\end{minipage}
\begin{minipage}{0.4\textwidth}
	\begin{equation*}
\xymatrix{
\sutg{\frac{2n-1}{2}}{i}\ar[rr]^{\psm{\frac{2n-1}{2}}{n}}&&B_{n,i}\ar[dl]^{\psm{n}{n-1}}\\
&T_{n-1,i+1}\ar[ul]^{\psm{n-1}{\frac{2n-1}{2}}}&
}	
\end{equation*}
\end{minipage}
\end{center}

Furthermore, we have the following.
\begin{enumerate}
	\item (\cite[Theorem 1.12]{LY2020} or from the large surgery formula) When $n\geq 2g(K)+1$, we have
	$$I^{\sharp}(-S^3_{-n}(K))\cong\bigoplus_{|i|\leq g(K)}\sutg{\frac{2n-1}{2}}{i}\oplus \mathbb{C}^{n-2g-1}.$$
	\item (\cite[Proposition 5.5]{LY2021large}) We have
	$$\psp{n-1}{n-2}\circ\psm{n}{n-1}=\psm{n-1}{n-2}\circ\psp{n}{n-1}=0.$$
\end{enumerate}
\elem

\bthm\label{thm: almost L-space knots}
Suppose $K\subset S^3$ is an almost L-space knot. Then we have the following.
\begin{enumerate}
	\item If $g(K)\geq 2$, then
	$\dim \khii(S^3,K,i)\leq 1$
	for any $i\in\intg$ such that $|i|>1$. Furthermore, the knot $K$ is fibered and strongly quasi-positive.
	\item If $g(K)=1$, then either $K$ is the figure-eight or $\tau_I(K)=1$ and
	\begin{equation*}
		\khii(S^3,K,i)\cong\begin{cases}
			0&|i|>1\\
			\mathbb{C}^2&|i|=1\\
			\mathbb{C}~{\rm or}~\mathbb{C}^3&i=0\\
		\end{cases}
	\end{equation*}

	\item If $g(K)=2$, then
	\begin{equation*}
		\khii(S^3,K,i)\cong\begin{cases}
			0&|i|>2\\
			\mathbb{C}&|i|=2\\
			\mathbb{C}~{\rm or}~\mathbb{C}^2&|i|=1\\
			\mathbb{C}~{\rm or}~\mathbb{C}^3&i=0\\
		\end{cases}
	\end{equation*}

\end{enumerate}
\ethm
\bpf
Suppose $n\in\mathbb{N}_+$ such that $\dim I^\sharp(S_n^3(K))=n$. From \cite[Section 2.2]{baldwin2020concordance}, we have the following exact triangle
\begin{equation*}
	\xymatrix{
	I^\sharp(S_n^3(K))\ar[rr]&&I^\sharp(S_{n+1}^3(K))\ar[dl]\\
	&I^\sharp(S^3)\ar[ul]&
	}
\end{equation*}
From the fact that $I^\sharp(S^3)\cong \mathbb{C}$, we know either
$\dim I^{\sharp}(S^3_{n+1}(K))=n+1$, which implies that $K$ is an instanton L-space knot and hence a contradiction, or 
$\dim I^{\sharp}(S^3_{n+1}(K))=n+3$. Hence, by induction we can assume that $n\geq 2g+1$. Suppose $\widebar{K}$ is the mirror of $K$. We have
$I^{\sharp}(-S^3_{-n}(\widebar{K}))\cong I^{\sharp}(S^3_n(K))$. From now on, all sutured instanton homologies are for the mirror knot. Applying Lemma \ref{lem: bypasses for T and B}, we know that
$$\bigoplus_{|i|\leq g} \sutg{\frac{2n-1}{2}}{i}\cong\mathbb{C}^{2g+3}.$$
From Lemma \ref{lem: symmetry and dimension of Gamma_n} part (1), we know that
$$\sutg{\frac{2n-1}{2}}{i}\cong \sutg{\frac{2n-1}{2}}{-i}.$$
From \cite[Proposition 1.21]{LY2021}, we know that
$$\chi(\sutg{\frac{2n-1}{2}}{i})=\chi(I^{\sharp}(-S^3))=1.$$
As a result, we conclude that
$$\sutg{\frac{2n-1}{2}}{i}\cong \mathbb{C} {\rm ~when} ~0<|i|\leq g ~{\rm and}~
\sutg{\frac{2n-1}{2}}{0}\cong \mathbb{C}^3.$$

When $g(K)\geq 2$, the argument for the fact that
$\dim \khii(-S^3,\widebar{K},i)\leq 1$
in the proof of \cite[Theorem 5.14]{LY2021large}
applies verbatim for $|i|>1$. In particular, for $i=g$, we can apply \cite[Lemma 5.7]{LY2021large} (with $m=g$) to conclude that (for the knot $\widebar{K}$)
$B_{n,g-1}=0.$
By Lemma \ref{lem: symmetry and dimension of Gamma_n} part (1), we know that
$$\sutg{n}{1-g+\frac{n-1}{2}}\cong \sutg{n}{g-1-\frac{n-1}{2}}=B_{n,g-1}=0.$$
From \cite[Section 5]{li2019direct}, we know that
$$\khii(-S^3,\widebar{K},1-g)\cong \sutg{n}{1-g+\frac{n-1}{2}}=0,$$
which, implies that
$\tau_I(\widebar{K})=-g$
as in Definition \ref{defn: tau}.
By Lemma \ref{lem: basic properties of tau}, we know
$\tau_I(K)=g.$

From \cite[Proposition 7.16]{kronheimer2011knot} and \cite[Proposition 4.1]{kronheimer2010instanton}, the fact that $$\dim \khii(-S^3,\widebar{K},g(K))\leq 1$$
implies $K$ is fibered. The fibration gives rise to a partial open book decomposition and hence a contact structure $\xi$ on $S^3$. Note that $K$ is strongly quasi-positive if and only if $\xi$ is tight on $S^3$. We can perturb $K$ so that $K$ is Lengendrian in $(S^3,\xi)$ and, furthermore, the knot complement $S^3\backslash N(K)$ is obtained by removing a standard tight contact neighborhood of $K$ from $S^3$. Let $\xi^\p$ be the restriction of $\xi$ on ${S^3\backslash N(K)}$. Hence $\partial (S^3\backslash N(K))$ is convex with dividing set described by the suture $\Ga_m$ for some integer $m\in\intg$. We can perform suitable stabilizations to make $m\geq2g(K)+1$. In \cite{baldwin2016contact}, Baldwin-Sivek defined a contact invariant
$\theta(\xi^\p)\in \sut{m}.$
By the proof of \cite[Theorem 1.17]{baldwin2018khovanov} and (the conclusion of) \cite[Theorem 1.18]{baldwin2018khovanov}, we know that
$$\theta(\xi^\p)\neq 0\in \sutg{m}{g}\cong \mathbb{C}.$$
We can attach a contact $2$-handle along the meridian of $K$ to $(S^3,\Ga_m)$ so that the sutured manifold becomes $S^3(1)$, which is a $3$-ball with a connected simple closed curve as the suture. After gluing, the contact structure $\xi^\p$ on $(S^3\backslash N(K),\Ga_m)$ becomes the restriction of $\xi$ on ${S^3(1)}$. So by \cite[Theorem 1.2]{baldwin2016contact}, we have
$$F_m(\theta(\xi^\p))=\theta(\xi|_{S^3(1)})\in \shi(-S^3(1))= I^{\sharp}(-S^3).$$
where $F_m$ is the map associated to the contact $2$-handle attachment. Note that, by \cite[Proposition 3.17]{li2019tau} and the fact that $\tau_I(K)=g$, we know that
$$\theta(\xi|_{S^3(1)})\neq0$$
which implies that $\xi$ is tight on $S^3$ by \cite[Theorem 1.3]{baldwin2016contact}. Hence we conclude that $K$ is strongly quasi-positive.

To prove the arguments when $g(K)\le 2$, we need to unpack the proof of \cite[Lemma 5.7 and Lemma 5.8]{LY2021large}. First, assume $g(K)=1$. From the above discussions, we can pick $n\geq 6$ and have
\begin{equation*}
	\sutg{\frac{2n-1}{2}}{i}\cong\sutg{\frac{2n-3}{2}}{i}\cong\sutg{\frac{2n-5}{2}}{i}\cong
	\begin{cases}
		\mathbb{C}&|i|=1\\
		\mathbb{C}^3&i=0.
	\end{cases}
\end{equation*}

From Lemma \ref{lem: bypass n,n+1,mu} and the definition of $T_{n}$, we know that
$$T_{n,1}=\sutg{n}{1+\frac{n-1}{2}}\cong \khii(-S^3,\widebar{K},1).$$
Assume 
$T_{n,1}\cong \khii(-S^3,\widebar{K},1)\cong \mathbb{C}^k.$
Lemma \ref{lem: bypasses for T and B} leads to the following diagram where the vertical and horizontal sequences are exact:
\begin{equation*}
	\xymatrix{
	&\sutg{\frac{2n-1}{2}}{1}\cong\mathbb{C}\ar[d]&\\
	&T_{n,1}\cong\mathbb{C}^k\ar[d]^{\psi_{-,n-1}^{n,1}}&\\
	\sutg{\frac{2n-3}{2}}{0}\cong\mathbb{C}^3\ar[r]^{\quad\quad\psi_{+,n-1}^{2n-3,0}}&B_{n-1,0}\ar[r]^{\psi_{+,n-2}^{n-1,0}\quad\quad}&T_{n-2,1}\cong\mathbb{C}^k
	}
\end{equation*}
where the second superscript of the bypass map indicates the grading. Note that, since $\sutg{\frac{2n-1}{2}}{1}\cong\mathbb{C}$, the map $\psi_{-,n-1}^{n,1}$ is either injective or surjective. 

{\bf Genus 1, case 1} $\psi_{-,n-1}^{n,1}$ is surjective. Then by the exactness
$B_{n-1,0}\cong \mathbb{C}^{k-1}.$
From Lemma \ref{lem: bypasses for T and B} part (3), we know that $\psi_{+,n-2}^{n-1,0}=0$ so from the exactness we know that
$3=k-1+k$, which means $k=2.$ Thus, we know that
$\khii(-S^3,\widebar{K},\pm1)\cong\mathbb{C}^k=\mathbb{C}^2.$
Applying Lemma \ref{lem: bypass n,n+1,mu}, we know that
$$\dim \khii(-S^3,\widebar{K},0)\leq\dim \sutg{n}{1+\frac{n-1}{2}}+\dim\sutg{n-1}{0+\frac{n-2}{2}}.$$
From the definitions of $T_{n,i}$ and $B_{n,i}$ and Lemma \ref{lem: symmetry and dimension of Gamma_n}, we know that
$$\dim \sutg{n}{1+\frac{n-1}{2}}=\dim T_{n,1}=k=2,\text{ and }$$
$$\dim \sutg{n-1}{0+\frac{n-2}{2}}=\dim \sutg{n-1}{-\frac{n-2}{2}}=\dim B_{n-1,0}=k-1=1.$$
Hence $\dim \khii(-S^3,\widebar{K},0)\leq3.$
From the Euler characteristic result in \cite[Theorem 1.1]{kronheimer2010instanton}, we know that $\dim \khii(-S^3,\widebar{K},0)$ is odd, so it must be either $1$ or $3$.

It remains to show that $\tau_I(K)=1$ for all such knots. By \cite[Section 5]{li2019direct}, if $n\geq3$ then
$$\khii(-S^3,\widebar{K},0)\cong\sutg{n}{0+\frac{n-1}{2}}\cong \mathbb{C}^{k-1}\sutg{n}{-1+\frac{n-1}{2}}\cong\mathbb{C}.$$
Note that the last isomorphism follows from Lemma \ref{lem: structure of Gamma_n}. Also, there is an exact triangle by Lemma \ref{lem: the U map} part (2).
\begin{equation*}
	\xymatrix{
	\khii(-S^3,\widebar{K},0)\ar[rr]^{U}&&\khii(-S^3,\widebar{K},-1)\ar[dl]\\
	&\khii(-S^3,\widebar{K},-1)\cong \mathbb{C}^2\ar[lu]&
	}
\end{equation*}
Hence we know
$$U|_{\khii(-S^3,\widebar{K},0)}=0.$$
and hence by the definition of $\tau_I$ and Lemma \ref{lem: basic properties of tau} we know
$\tau_I(K)=-\tau_I(\widebar{K})=1.$

{\bf Genus 1, case 2} $\psi_{-,n-1}^{n,1}$ is injective. 
Then by the exactness
$B_{n-1,0}\cong \mathbb{C}^{k+1}.$
Lemma \ref{lem: bypasses for T and B} implies another exact triangle
\begin{equation*}
	\xymatrix{
	&\sutg{\frac{2n-5}{2}}{1}\cong\mathbb{C}\ar[d]\\
	B_{n-1,0}\cong\mathbb{C}^{k+1}\ar[r]^{\psi_{+,n-2}^{n-1,0}}&T_{n-2,1}\cong\mathbb{C}^{k}\ar[d]^{\psi_{-,n-3}^{n-2,1}}\ar[r]&\sutg{\frac{2n-3}{2}}{0}\cong\mathbb{C}^3\\
	&B_{n-3,0}\cong \mathbb{C}^{k+1}
	}
\end{equation*}
The vertical exactness implies that $\psi_{-,n-3}^{n-2,1}$ is injective and hence from Lemma \ref{lem: bypasses for T and B} part (3), we know that $\psi_{+,n-2}^{n-1,0}$ is zero. Hence from the horizontal exactness we know
$k+1+k=3$, which means $ k=1.$

From \cite[Proposition 4.1]{kronheimer2010instanton} we know $K$ is fibered. It is well known that there are only two genus-one fibered knots in $S^3$, namely the trefoil and the figure-eight, among which the trefoil is an L-space knot. Hence $K$ is the figure-eight.

Finally we study the case of $g(K)=2$. First, as in the proof of part (1), since $\khii(-S^3,\widebar{K},2)\neq0$, we can apply \cite[Lemma 5.7]{LY2021large} directly with $m=2$ and conclude that for $n\geq 9$, we have
$\khii(-S^3,\widebar{K},2)\cong T_{n,2}\cong\mathbb{C}$ and $B_{n-1,1}=0.$
Note that we have
$$\sutg{n+1}{-1+\frac{n-1}{2}}\cong B_{n+1,1}\cong B_{n-1,1}=0.$$
and hence from Lemma \ref{lem: bypass n,n+1,mu} and Lemma \ref{lem: symmetry and dimension of Gamma_n} we know
$$\khii(-S^3,\widebar{K},-1)=\sutg{\mu}{-1}\cong\sutg{n}{0+\frac{n-1}{2}}\cong \sutg{n}{0-\frac{n-1}{2}}=B_{n,0}\cong B_{n-1,0}.$$
Then the argument above for genus-one almost L-space knots applies verbatim and we can conclude the following two cases:

{\bf Genus 2, case 1} $T_{n,1}\cong\mathbb{C}$ and $B_{n-1,0}\cong\mathbb{C}^2$.

{\bf Genus 2, case 2} $T_{n,1}\cong\mathbb{C}^2$ and $B_{n-1,0}\cong\mathbb{C}$.
Note that in both cases from Lemma \ref{lem: bypass n,n+1,mu} and Lemma \ref{lem: symmetry and dimension of Gamma_n} we know that
\beq
\dim \khii(-S^3,\widebar{K},0)&\leq\dim \sutg{n}{1+\frac{n-1}{2}}+\dim\sutg{n+1}{0+\frac{n}{2}}\\
&=\dim T_{n,1}+\dim\sutg{n+1}{0-\frac{n}{2}}\\
&=\dim T_{n,1}+\dim B_{n+1,0}\\
&=\dim T_{n,1}+\dim B_{n-1,0}\\
&=3.
\eeq
Hence we conclude the proof of part (3).
\epf

\brem
For genus-two almost L-space knots, we know
$$\dim \khii(S^3,K,2)=1~{\rm and ~}\dim \khii(S^3,K,1)=1{\rm~or~}2.$$
Recent techniques developed in \cite{baldwin21t25,BLSY21} can show that $\dim \khii(S^3,K,2)=1$ implies that $K=T_{2,\pm 5}$, while the case $\dim \khii(S^3,K,1)=2$ is still open.
\erem

The techniques in proving the above lemma can also lead to the following.

\bcor\label{cor: 1-surgery of genus 1 knot}
Suppose $K$ is a genus-one knot such that
$I^{\sharp}(S^3_1(K))=2d+1,$
then either
\begin{enumerate}
	\item $\dim \khii(S^3,K,1)=d+1$ and $\tau_I(K)=1$, or 
	\item $\dim \khii(S^3,K,1)=d$ and $\tau_I(K)\leq 0$.
\end{enumerate}
\ecor
\bpf
Note that $g(K)=1$ so by \cite[Section 1.1 and Theorem 1.1]{baldwin2020concordance}, we know that
$\dim I^{\sharp}(S^3_3(K))=2d+3.$
From Lemma \ref{lem: bypasses for T and B} we know that if $n\geq 7$,
$$\dim \bigoplus_{-1\leq i \leq 1} \sutg{\frac{2n-1}{2}}{i}=2d+3.$$
Lemma \ref{lem: bypasses for T and B} implies a triangle
\begin{equation*}
	\xymatrix{
	\sutg{\frac{2n-1}{2}}{1}\ar[rr]&&\sutg{n}{2+\frac{n-1}{2}}\ar[dl]\\
	&\sutg{n-1}{1-\frac{n-2}{2}}\ar[lu].&
	}
\end{equation*}
From Lemma \ref{lem: structure of Gamma_n} with $Y=S^3$, we know that
$$\sutg{n}{2+\frac{n-1}{2}}=0~{\rm and~}\sutg{n-1}{1-\frac{n-2}{2}}\cong\mathbb{C}$$
so we conclude that
$$\dim \sutg{\frac{2n-1}{2}}{1}=\dim \sutg{\frac{2n-1}{2}}{-1}=1.$$
As a result, we have
$$\dim \sutg{\frac{2n-1}{2}}{0}=2d+1.$$
The argument in the proof of Theorem \ref{thm: almost L-space knots} for the case $g=1$ applies. The original setup $\dim \sutg{\frac{2n-1}{2}}{0}\cong\mathbb{C}^3$ is the case $d=1$. So as in that proof, we have two cases

{\bf Case 1}. $\khii(-S^3,\widebar{K},1)\cong T_{n,1}\cong \mathbb{C}^{d+1}$, $B_{n-1,0}\cong\mathbb{C}^{d}$, and $\tau_I(K)=1$.

{\bf Case 2}. $\khii(-S^3,\widebar{K},1)\cong T_{n,1}\cong \mathbb{C}^{d}$, $B_{n-1,0}\cong\mathbb{C}^{d+1}$. As in the proof of Theorem \ref{thm: almost L-space knots}, we know
$$\khii(-S^3,\widebar{K},0)\cong \sutg{n-1}{0+\frac{n-1}{2}}\cong B_{n-1,0}\cong \mathbb{C}^{d+1}$$
and
$$\khii(-S^3,\widebar{K},-1)\cong \sutg{n-1}{-1+\frac{n-1}{2}}\cong\mathbb{C}.$$
From the exact triangle in Lemma \ref{lem: the U map} part (2)
\begin{equation*}
	\xymatrix{
	\khii(-S^3,\widebar{K},0)\cong\mathbb{C}^{d+1}\ar[rr]^{U}&&\khii(-S^3,\widebar{K},-1)\cong\mathbb{C}\ar[dl]\\
	&\khii(-S^3,\widebar{K},-1)\cong \mathbb{C}^{d}\ar[lu]&
	}
\end{equation*}
the map $U:\khii(-S^3,\widebar{K},0)\ra \khii(-S^3,\widebar{K},-1)$
is surjective and hence $\tau_I(\widebar{K})\geq 0$ which implies that $\tau_I(K)\leq0.$
\epf

\bcor\label{cor: dim >=5}
Suppose $K\subset S^3$ is a knot with $g(K)\geq 2$. Then
$$\dim I^{\sharp}(S^3_1(K))\geq 5.$$
\ecor
\bpf
Suppose the contrary, that $\dim I^{\sharp}(S^3_1(K))\leq 3$. Then there are two cases: $K$ is either an instanton L-space knot or an almost L-space knot. If $K$ is an instanton L-space knot we can apply the main result in \cite{lidman2020framed} (or \cite{LY2021large}) and conclude that $\dim I^{\sharp}(S^3_1(K))\geq 5$ directly. If $K$ is an almost L-space then from Theorem \ref{thm: almost L-space knots} and \cite[Theorem 1.2]{li2019tau}, we know that the invariant $\nu^\sharp(K)$ in \cite{baldwin2020concordance} satisfies
$$\nu^{\sharp}(K)\ge 2\tau^\sharp(K)-1=2\tau_I(K)-1=2g(K)-1\geq 3.$$
Then from \cite[Theorem 1.1]{baldwin2020concordance} we know that $\dim I^{\sharp}(S^3_1(K))\geq 5$.
\epf

\bcor\label{cor: 15n43522}
Suppose $K=15n_{43522}$, then we have $\tau_I(K)=0$ and
\begin{equation*}
	\khii(S^3,K,i)\cong\begin{cases}
		0&|i|>1\\
		\mathbb{C}^2&|i|=1\\
		\mathbb{C}^5&|i|=0
	\end{cases}
\end{equation*}
\ecor

\bpf
From \cite{BS2022nearly}, we know that $g(K)=1$, $\Delta_K(t)=2t-3+2t^{-1}$, and $\dim_\mathbb{Q}\widehat{HFK}(S^3,K,1)=2$. From \cite[Corollary 1.4]{LY2022nearly}, we know that
$\dim \khii(S^3,K,1)=2$. From the Euler characteristic result in \cite[Theorem 1.1]{kronheimer2010instanton}, we know that
$\dim \khii(S^3,K,1)$ is either $3$ or $5$. If it is $3$-dimensional, from \cite[Proposition 6.8]{LY2021large}, we know that either
$$\dim I^{\sharp}(S^3_1(K))=3~{\rm or~}\dim I^{\sharp}(S^3_{-1}(K))=3.$$
However this contradicts Corollary \ref{cor: dim >=5} and the facts that
$$S^3_1(K)\cong \pm S^3_{-1}(9_{42})~{\rm and~}S^3_{-1}(K)\cong \pm S^3_{-1}(8_{20})$$
as in the proof of \cite[Proposition A1]{BS2022nearly}. As a result, we must have
$$\dim \khii(S^3,K,1)=5.$$
\epf

\bpf[Proof of Theorem \ref{thm: genus-one nearly fibred}]
We know from \cite{LY2022nearly,BS2022nearly} that up to mirror $K$ must be one of the following knots:
$$5_2,~D^{\pm}_2(J),~15n_{43522}, P(-3,3,2n+1).$$
Since $K=5_2$ is an alternating knot it follows from \cite{kronheimer2011khovanov} that
$$\dim KHI(S^3,K)=||\Delta_{K}(t)||=7.$$
For $K=D_2^+(J)$ or $P(-3,3,2n+1)$, we know that
$$\Delta_{K}(t)=-2t+5-2t^{-1}.$$
From \cite[Theorem 1.1]{kronheimer2010instanton} we know that
$$\dim KHI(S^3,K,0)\geq 5.$$
From the proof of \cite[Proposition 6.3]{LY2021large}, we know that
$$\dim KHI(S^3,K,0)\leq 5.$$
As a result, we have
$$\dim KHI(S^3,K,0)=5.$$

For $K=D_2^-(J)$ or $15n_{43522}$, we know that
$$\Delta_{K}(t)=2t-3+2t^{-1}.$$
From the above argument we know that
$$\dim KHI(S^3,K,0)=3{\rm ~or}~5.$$
If $\dim KHI(S^3,K,0)=3$ then \cite[Proposition 6.8]{LY2021large} and Corollary \ref{cor: 1-surgery of genus 1 knot} imply that
$$\tau_{I}(K)=\pm1,$$
which contradicts Corollary \ref{cor: 15n43522} and Lemma \ref{lem: tau for Whitehead doubles}.
\epf

\bpf[Proof of Theorem \ref{thm: almost L space knots, intro}]
Part (1) follows from Theorem \ref{thm: almost L-space knots}. We prove part (2) as follows. From Theorem \ref{thm: almost L-space knots} part (2), when $K$ is a genus-one almost L-space knot, we have either $\khii(S^3,K,1)\cong\mathbb{C}$ so that $K$ is the figure eight, or $\khii(S^3,K,1)\cong\mathbb{C}^2$. In the latter case, we know from Theorem \ref{thm: genus-one nearly fibred} that $K=\widebar{5}_2$, is indeed an almost L-space knot again by \cite[Theorem 1.20]{LY2021large}.
\epf

\bpf[Proof of Corollary \ref{cor: 1 surgery having dimension 3}]
If $\dim I^{\sharp}(S^3_1(K))=3$, we know that either $K$ is an L-space knot or an almost L-space knot. From Corollary \ref{cor: dim >=5} we know that $g(K)=1$. Then the corollary follows from Theorem \ref{thm: almost L space knots, intro} part (2).
\epf

%————End from here————
%\newpage
\section{Conflict of interest}
On behalf of all authors, the corresponding author states that there is no conflict of interest.

\section{Data availability}
On behalf of all authors, the corresponding author states that there is no associated data.

\bibliographystyle{alpha}

\begin{thebibliography}{ABDS22}

\bibitem[ABDS22]{alfieri2020framed}
Antonio Alfieri, John~A. Baldwin, Irving Dai, and Steven Sivek.
\newblock Instanton {F}loer homology of almost-rational plumbings.
\newblock {\em Geom. Topol.}, 26(5):2237--2294, 2022.

\bibitem[BHS21]{baldwin21t25}
John~A. Baldwin, Ying Hu, and Steven Sivek.
\newblock {Khovanov homology and the cinquefoil}.
\newblock {\em ArXiv:2105.12102, v1}, 2021.

\bibitem[BLSY24]{BLSY21}
John~A. Baldwin, Zhenkun Li, Steven Sivek, and Fan Ye.
\newblock Small {D}ehn surgery and {${\rm SU}(2)$}.
\newblock {\em Geom. Topol.}, 28(4):1891--1922, 2024.

\bibitem[BS15]{baldwin2015naturality}
John~A. Baldwin and Steven Sivek.
\newblock Naturality in sutured monopole and instanton homology.
\newblock {\em J. Differ. Geom.}, 100(3):395--480, 2015.

\bibitem[BS16a]{baldwin2016contact}
John~A. Baldwin and Steven Sivek.
\newblock A contact invariant in sutured monopole homology.
\newblock {\em Forum Math. Sigma}, 4:e12, 82, 2016.

\bibitem[BS16b]{baldwin2016instanton}
John~A. Baldwin and Steven Sivek.
\newblock Instanton {F}loer homology and contact structures.
\newblock {\em Selecta Math. (N.S.)}, 22(2):939--978, 2016.

\bibitem[BS21]{baldwin2020concordance}
John~A. Baldwin and Steven Sivek.
\newblock Framed instanton homology and concordance.
\newblock {\em J. Topol.}, 14(4):1113--1175, 2021.

\bibitem[BS22a]{baldwin22knot52}
John~A. Baldwin and Steven Sivek.
\newblock Characterizing slopes for $5_2$.
\newblock {\em ArXiv:2209.09805, v1}, 2022.

\bibitem[BS22b]{BS2022nearly}
John~A. Baldwin and Steven Sivek.
\newblock {Floer homology and non-fibered knot detection}.
\newblock {\em ArXiv:2208.03307, v1}, 2022.

\bibitem[BS22c]{baldwin2018khovanov}
John~A. Baldwin and Steven Sivek.
\newblock Khovanov homology detects the trefoils.
\newblock {\em Duke Math. J.}, 171(4):885--956, 2022.

\bibitem[BS23]{baldwin2019lspace}
John~A. Baldwin and Steven Sivek.
\newblock Instantons and {L}-space surgeries.
\newblock {\em J. Eur. Math. Soc. (JEMS)}, 25(10):4033--4122, 2023.

\bibitem[Don02]{donaldson2002floer}
S.~K. Donaldson.
\newblock {\em Floer homology groups in {Y}ang-{M}ills theory}, volume 147 of
  {\em Cambridge Tracts in Mathematics}.
\newblock Cambridge University Press, Cambridge, 2002.
\newblock With the assistance of M. Furuta and D. Kotschick.

\bibitem[GL23]{li2019decomposition}
Sudipta Ghosh and Zhenkun Li.
\newblock Decomposing sutured monopole and instanton {F}loer homologies.
\newblock {\em Selecta Math. (N.S.)}, 29(3):Paper No. 40, 60, 2023.

\bibitem[GLW24]{li2019tau}
Sudipta Ghosh, Zhenkun Li, and C.-M.~Michael Wong.
\newblock On the tau invariants in instanton and monopole {F}loer theories.
\newblock {\em J. Topol.}, 17(2):Paper No. e12346, 53, 2024.

\bibitem[Han20]{Hanselman20surgery}
Jonathan Hanselman.
\newblock {Heegaard Floer homology and cosmetic surgeries in $S ^3$}.
\newblock {\em ArXiv: 1906.06773, v2}, 2020.

\bibitem[Hed07]{Hedden2007WHD}
Matthew Hedden.
\newblock {Knot Floer homology of Whitehead doubles}.
\newblock {\em Geom. Topol.}, 11(4):2277--2338, 2007.

\bibitem[Juh06]{juhasz2006holomorphic}
Andr\'{a}s Juh\'{a}sz.
\newblock Holomorphic discs and sutured manifolds.
\newblock {\em Algebr. Geom. Topol.}, 6:1429--1457, 2006.

\bibitem[KM10a]{kronheimer2010instanton}
Peter~B. Kronheimer and Tomasz~S. Mrowka.
\newblock Instanton {F}loer homology and the {A}lexander polynomial.
\newblock {\em Algebr. Geom. Topol.}, 10(3):1715--1738, 2010.

\bibitem[KM10b]{kronheimer2010knots}
Peter~B. Kronheimer and Tomasz~S. Mrowka.
\newblock Knots, sutures, and excision.
\newblock {\em J. Differ. Geom.}, 84(2):301--364, 2010.

\bibitem[KM11a]{kronheimer2011khovanov}
Peter~B. Kronheimer and Tomasz~S. Mrowka.
\newblock Khovanov homology is an unknot-detector.
\newblock {\em Publ. Math. Inst. Hautes \'{E}tudes Sci.}, 113:97--208, 2011.

\bibitem[KM11b]{kronheimer2011knot}
Peter~B. Kronheimer and Tomasz~S. Mrowka.
\newblock Knot homology groups from instantons.
\newblock {\em J. Topol.}, 4(4):835--918, 2011.

\bibitem[Li20]{li2018contact}
Zhenkun Li.
\newblock Contact structures, excisions and sutured monopole {F}loer homology.
\newblock {\em Algebr. Geom. Topol.}, 20(5):2553--2588, 2020.

\bibitem[Li21a]{li2018gluing}
Zhenkun Li.
\newblock Gluing maps and cobordism maps in sutured monopole and instanton
  {F}loer theories.
\newblock {\em Algebr. Geom. Topol.}, 21(6):3019--3071, 2021.

\bibitem[Li21b]{li2019direct}
Zhenkun Li.
\newblock Knot homologies in monopole and instanton theories via sutures.
\newblock {\em J. Symplectic Geom.}, 19(6):1339--1420, 2021.

\bibitem[Lim10]{Lim2009}
Yuhan Lim.
\newblock {Instanton homology and the Alexander polynomial}.
\newblock {\em Proc. Amer. Math. Soc.}, 138(10):3759--3768, 2010.

\bibitem[LPCS22]{lidman2020framed}
Tye Lidman, Juanita Pinz\'{o}n-Caicedo, and Christopher Scaduto.
\newblock Framed instanton homology of surgeries on {L}-space knots.
\newblock {\em Indiana Univ. Math. J.}, 71(3):1317--1347, 2022.

\bibitem[LY21]{LY2021large}
Zhenkun Li and Fan Ye.
\newblock {SU(2) representations and a large surgery formula}.
\newblock {\em ArXiv:2107.11005, v1}, 2021.

\bibitem[LY22a]{LY2020}
Zhenkun Li and Fan Ye.
\newblock {Instanton Floer homology, sutures, and {H}eegaard diagrams}.
\newblock {\em J. Topol.}, 15(1):39--107, 2022.

\bibitem[LY22b]{LY2022integral1}
Zhenkun Li and Fan Ye.
\newblock {Knot surgery formulae for instanton Floer homology I: the main
  theorem}.
\newblock {\em ArXiv:2206.10077, v2}, 2022.

\bibitem[LY22c]{LY2022nearly}
Zhenkun Li and Fan Ye.
\newblock {Seifert surface complements of nearly fibered knots}.
\newblock {\em ArXiv:2208.05382, v1}, 2022.

\bibitem[LY23a]{LY2021enhanced}
Zhenkun Li and Fan Ye.
\newblock {An enhanced Euler characteristic of sutured instanton homology}.
\newblock {\em Int. Math. Res. Not. IMRN}, page rnad066, 04 2023.

\bibitem[LY23b]{LY2021}
Zhenkun Li and Fan Ye.
\newblock Instanton {F}loer homology, sutures, and {E}uler characteristics.
\newblock {\em Quantum Topol.}, 14(2):201--284, 2023.

\bibitem[OS03]{Ozsvath2003}
Peter~S. Ozsv{\'{a}}th and Zolt{\'{a}}n Szab{\'{o}}.
\newblock {Heegaard Floer homology and alternating knots}.
\newblock {\em Geom. Topol.}, 7(1):225--254, 2003.

\bibitem[OS04]{ozsvath2004holomorphicknot}
Peter~S. Ozsv\'{a}th and Zolt\'{a}n Szab\'{o}.
\newblock Holomorphic disks and knot invariants.
\newblock {\em Adv. Math.}, 186(1):58--116, 2004.

\bibitem[OS08]{Ozsvath2008integral}
Peter~S. Ozsv{\'{a}}th and Zolt{\'{a}}n Szab{\'{o}}.
\newblock {Knot {F}loer homology and integer surgeries}.
\newblock {\em Algebr. Geom. Topol.}, 8(1):101--153, 2008.

\bibitem[OS11]{Ozsvath2011rational}
Peter~S. Ozsv{\'{a}}th and Zolt{\'{a}}n Szab{\'{o}}.
\newblock {Knot {F}loer homology and rational surgeries}.
\newblock {\em Algebr. Geom. Topol.}, 11(1):1--68, 2011.

\bibitem[Pet13]{Petkova2009thin}
Ina Petkova.
\newblock {Cables of thin knots and bordered Heegaard Floer homology}.
\newblock {\em Quantum Topol.}, 4(4):377--409, 2013.

\bibitem[Sca15]{scaduto2015instanton}
Christopher Scaduto.
\newblock Instantons and odd {K}hovanov homology.
\newblock {\em J. Topol.}, 8(3):744--810, 2015.

\bibitem[XZ23]{yixiesu2}
Yi~Xie and Boyu Zhang.
\newblock On meridian-traceless su(2)–representations of link groups.
\newblock {\em Advances in Mathematics}, 418:108947, 2023.

\end{thebibliography}

\end{document}